\pdfoutput=1
\RequirePackage{ifpdf}
\ifpdf 
\documentclass[pdftex]{sigma}
\else
\documentclass{sigma}
\fi

\usepackage{ytableau}

\numberwithin{equation}{section}

\newtheorem{Theorem}{Theorem}[section]
\newtheorem{Corollary}[Theorem]{Corollary}
\newtheorem{Lemma}[Theorem]{Lemma}
\newtheorem{Proposition}[Theorem]{Proposition}
 { \theoremstyle{definition}
\newtheorem{Definition}[Theorem]{Definition}

\newtheorem{Problem}[Theorem]{Problem}

\newtheorem{Remark}[Theorem]{Remark} }

\newcommand{\ce}{\mathcal{E}}
\newcommand{\ct}{\mathcal{T}}
\newcommand{\cc}{\boldsymbol{c}}
\newcommand{\bg}{\boldsymbol{g}}
\newcommand{\bG}{\boldsymbol{G}}
\newcommand{\bn}{\boldsymbol{n}}
\newcommand{\II}{{\rm I\hspace{-.2mm}I}}
\newcommand{\IIo}{\hspace{0.4mm}\mathring{\rm{ I\hspace{-.2mm} I}}{\hspace{.0mm}}}

\newcommand{\IVnn}{{{\bf\rm I\hspace{-.2mm} V}{\hspace{.2mm}}}{}}

\newcommand{\hd }{\hat{D}}

\renewcommand{\=}{\stackrel \Lambda =}
\newcommand{\csdot}{\hspace{-0.75mm} \cdot \hspace{-0.75mm}}

\DeclareMathOperator{\Id}{Id}

\begin{document}
\allowdisplaybreaks

\newcommand{\arXivNumber}{2405.07692}

\renewcommand{\PaperNumber}{002}

\FirstPageHeading

\ShortArticleName{Holography of Higher Codimension Submanifolds: Riemannian and Conformal}

\ArticleName{Holography of Higher Codimension Submanifolds:\\ Riemannian and Conformal}

\Author{Samuel BLITZ and Josef \v{S}ILHAN}

\AuthorNameForHeading{S.~Blitz and J.~\v{S}ilhan}

\Address{Department of Mathematics and Statistics, Masaryk University,\\
Building 08, Kotl\'a\v{r}sk\'a 2, Czech Republic}
\Email{\href{mailto:blitz@math.muni.cz}{blitz@math.muni.cz}, \href{mailto:silhan@math.muni.cz}{silhan@math.muni.cz}}

\ArticleDates{Received June 01, 2024, in final form December 16, 2024; Published online January 04, 2025}

\Abstract{We provide a natural generalization to submanifolds of the holographic method used to extract higher-order local invariants of both Riemannian and conformal embeddings, some of which depend on a choice of parallelization of the normal bundle. Qualitatively new behavior is observed in the higher-codimension case, giving rise to new invariants that obstruct the order-by-order construction of unit defining maps. In the conformal setting, a novel invariant (that vanishes in codimension 1) is realized as the leading transverse-order term appearing in a holographically-constructed Willmore invariant. Using these same tools, we also investigate the formal solutions to extension problems off of an embedded submanifold.}

\Keywords{Riemannian geometry; conformal geometry; submanifold embeddings; holo\-graphy}

\Classification{53C18; 53A55; 53C21; 58J32}

\section{Introduction}
The local geometry of codimension $k$ embedded submanifolds $\Lambda^{d-k}$ in Riemannian $d$-manifolds $\bigl(M^d,g\bigr)$ has been studied extensively for over a century, with special focus being given to curves (codimension $d-1$ submanifolds) and hypersurfaces (codimension $1$ submanifolds). These cases stand out as describing particularly useful systems: curves can describe the trajectory of a particle in space and hypersurfaces can provide boundary data for elliptic PDEs. On the other hand, \textit{conformal manifolds} -- smooth manifolds $M^d$ endowed with an equivalence class of Riemannian metrics $\cc := \big\{\Omega^2 g\big\}$ for every $\Omega \in C^\infty_+ M$ -- have also been studied since the early 20th century. Early investigations include work by Cartan~\cite{CartanConformal}, Thomas~\cite{ThomasConformal} and Fialkow~\cite{Fialkow}, among others. In particular, \textit{conformal invariants} -- those quantities which are equivariant under the rescaling group action -- are central to any description of conformal geometry. In the last thirty years, a~method of Thomas~\cite{BEG} has been rediscovered that allows for the construction of such invariants in the same way that the Ricci calculus allows for the construction of invariants of a~Riemannian manifold: this is called the tractor calculus, which will be discussed in more detail in Section~\ref{conf-geom-sec}. (Another distinct but equally powerful method of producing conformal invariants by embedding conformal manifolds into higher dimensional manifolds solving particular problems was first introduced by Fefferman and Graham in~\cite{FGbook}.) A confluence of submanifold geometry and conformal geometry has gained popularity in the last few decades~\cite{Forms,WillB,GR,GZscatt,MondinoNguyen,SZ2019,Chang}, having been driven in part by physics applications~\cite{HS, Maldacena}. Specifically, conformal hypersurface geometry has taken center stage, with a focus on an extrinsic characterization of the embedding $\Sigma^{d-1} \hookrightarrow \bigl(M^d,\cc\bigr)$ itself~\cite{Will2,Will1}. Recent work conducted by Curry, Gover and Snell~\cite{RodSnellCurry} also investigated conformal submanifold geometry using tractor methods, however their approach differed significantly from the methods applied here. Furthermore, Poincar\'e--Einstein embedding approaches have also been used to construct conformal submanifold invariants, see for example~\cite{CGKTW} -- this approach avoids tractor calculus entirely.

One popular method to study any submanifold embedding structure, whether Riemannian, conformal, or otherwise, is called \textit{holography}. Given a submanifold embedding into some geometric structure $\Lambda \hookrightarrow (M,\mathcal{S})$, one can study the invariants of the embedding by studying solutions to geometrically-natural PDEs relating data on $\Lambda$ to data on $M$. In the study of Riemannian hypersurface embeddings $\Sigma \hookrightarrow (M,g)$, this method uses the existence of a unique \textit{unit defining function} for $\Sigma$, i.e., a function $s \in C^\infty M$ such that $\Sigma = \{p \in M \mid s(p) = 0\}$ and that $|{\rm d}s|^2_g = 1$. It is this second condition -- the eikonal equation -- that distinguishes this function from any other defining function, which need only have that $|{\rm d}s|^2_g > 0$. As $s$ is uniquely determined (up to orientation) by the embedding, its jets can be used to construct invariants of the embedding structure.

In the case of a two-sided conformal hypersurface embeddings $\Sigma \hookrightarrow (M,\cc)$, a similar method bears fruit. In that case, one seeks a double equivalence class (or density) $\sigma = [g;s] = \bigl[\Omega^2 g ; \Omega s\bigr]$ of metrics $g \in \cc$ and functions $s \in C^\infty M$ such that for any representative $[g;s]$, we have that on $M/\Sigma$, the metric $g^o := g/s^2$ has constant scalar curvature. (More detail is provided on these double equivalence classes in Section~\ref{conf-geom-sec}.) Finding such a double equivalence class solves the so-called ``singular Yamabe problem'', which can be shown~\cite{Loewner} to be equivalent to a conformal eikonal equation, given~by
\[
|{\rm d}s|^2_g - \frac{2s}{d} \biggl(\Delta^g s + \frac{1}{2(d-1)} \mathrm{Sc}^g s\biggr) = 1 .
\]
While a one-sided global solution to this problem always exists~\cite{ACF,Aviles,Maz}, it turns out that one cannot (in general) solve this to arbitrarily high order~\cite{ACF,Goal}. Nonetheless, one may still find a~unique \textit{unit defining density} up to order $d+1$ in $\sigma$ that is entirely determined by the embedding and the conformal structure. As it is unique, its jets below this order will therefore characterize the conformal embedding. A study of these jets was carried out in~\cite{BGW1}.

The goal of this article is to investigate the holographic method of unit defining functions as it applies to higher-codimension submanifold embeddings in both the Riemannian and conformal settings. In general, we find that this method does not yield true submanifold invariants, but rather invariants that depend on a choice of parallelization of the normal bundle on $\Lambda$. In the following subsection, we introduce notations and conventions that will be used throughout this article. In Section~\ref{riem-frames}, we first discuss when and how a canonical frame for the normal bundle of a Riemannian submanifold embedding can be constructed. In Section~\ref{unit-def-funcs}, we proceed with the holographic program in the Riemannian setting, culminating in an optimal construction in Theorem~\ref{ho-Riem}. In Section~\ref{conf-geom-sec}, we give a brief review of conformal geometry, before considering the frame construction in a conformal submanifold embedding in Section~\ref{conf-frames}. In Section~\ref{conf-densities-sec}, we investigate the use of the holographic method to study conformal submanifold embeddings, leading to a optimal (to second order) construction in Theorem~\ref{ho-conformal-construction}. While we find that we cannot proceed to particularly high order, we are still able to construct an interesting parallelization-dependent conformal submanifold invariant in the special case of a conformally-embedded surface. This is detailed in Section~\ref{Willmore-section} -- see equation~\eqref{surface-willmore}. Finally, in Section~\ref{ext-prob-sec}, we examine the submanifold extension problems, which are a natural problems to study given a holographically-constructed defining map. The results are contained in Theorems~\ref{sym-ext-prop}, \ref{ho-extension} and \ref{conformal-ext-solution}.

\subsection{Notations and conventions}
In this article, $M$ denotes a $d$-dimensional (with $d \geq 3$) smooth manifold endowed with a metric tensor $g$ or conformal class of metrics $\cc$, both of which will be positive definite. Furthermore, $\Lambda$~denotes a $(d-k)$-dimensional smooth submanifold of $M$, and $0 < k < d$. Throughout the entirety of this article, we work only locally and so will assume that trivializations are always available for any section of any bundle. In the Riemannian setting $(M,g)$, we denote the Levi-Civita connection by $\nabla$. Our convention for the curvature tensor is given~by
\[
R(x,y)z = (\nabla_x \nabla_y - \nabla_y \nabla_x)z - \nabla_{[x,y]} z .
\]
Here, $x$, $y$, $z$ are vector fields on $M$ and $[x,y]$ is the Lie bracket of $x$ and $y$. In general, for any bundle $\mathcal{V}M$ over $M$, we will denote sections of this bundle by $\Gamma(\mathcal{V} M)$. That is, $x,y,z \in \Gamma(TM)$. When necessary, we will refer to arbitrary tensor products of a bundle (and its dual) using capital Greek letter, so that, for example, $T^\Phi M$ is some unspecified tensor product of $TM$ and~$T^* M$.

For ease of calculation, we will use Penrose's abstract index notation, labelling such indices with letters from the beginning of the Latin alphabet $a$, $b$, $c$, etc., so that, for example, the~Riemann curvature tensor may be expressed as
\[
R_{ab}{}^c{}_d x^d = [\nabla_a, \nabla_b] x^c .
\]
Using the metric and its inverse, indices can be raised and lowered, as usual; the Einstein summation convention is used to indicate contraction. (Sometimes, when the notation is unambiguous, we will use dot product notation to represent contraction; for example, when $\omega$ is some 1-form, we will sometimes write $\nabla^a \omega_a =: \nabla \cdot \omega$.) In this language, we have $\mathrm{Ric}_{ab} := R_{ca}{}^c{}_b$ and its corresponding scalar is $\mathrm{Sc} := \mathrm{Ric}_a^a$. When $d > 2$, a particular trace-correction of the Ricci tensor, known as the Schouten tensor, arises frequently in conformal considerations,
\[
P_{ab} := \frac{1}{d-2} \left(\mathrm{Ric}_{ab} - \frac{1}{2(d-1)} \mathrm{Sc}\, g_{ab} \right) .
\]
The trace of this tensor is denoted by $J := P_a^a$. Written in terms of the Riemann tensor and the Schouten tensor, we may express the Weyl tensor as
\[
W_{abcd} = R_{abcd} - g_{ac} P_{bd} + g_{bc} P_{ad} + g_{ad} P_{bc} - g_{bd} P_{ac} .
\]
In dimensions $d > 2$, the Cotton tensor is given~by the covariant curl of the Schouten tensor,
\[
C_{abc} = \nabla_a P_{bc} - \nabla_b P_{ac} .
\]
When $d > 3$, this can be reexpressed in terms of the divergence of the Weyl tensor,
\[
C_{abc} = \frac{1}{d-3} \nabla^d W_{dcab} .
\]

In the context of submanifold embeddings, we will use Greek indices to represent tensors on normal frame bundles. Objects that are along the submanifold will be denoted by overbars~$\bar{\bullet}$, so that, for example, $\bar{\nabla}$ represents the induced Levi-Civita connection on $(\Lambda,\bar{g})$, where $\bar{g}$ is the induced metric from $\Lambda \hookrightarrow (M,g)$. In general for any index notation, we will use round brackets~$(\cdots)$ to indicate unit-normalized symmetrization and square brackets $[\cdots]$ to indicate unit-normalized antisymmetrization. Furthermore, we will use a subscript $\circ$ to indicate the tracefree part of a particular symmetric index configuration. However, because Latin and Greek indices represent different types of objects, we may sometimes write $(\alpha a b \beta)$ to specifically mean the symmetrization over $\alpha$ and $\beta$, and similarly for antisymmetrization. To be explicit, whichever type of index is closest to the opening and closing brackets is the type of index that is (anti)symmetrized over. On the other hand, when excluding an index from (anti)symmetrization is necessary, we will use vertical bars, so that $(\alpha | \beta | \gamma)$ refers to the symmetrization over~$\alpha$ and~$\gamma$. Note that Greek indices will always be subscripts, so repeated Greek indices will denote summation in the usual way.

\section[Orthonormal frames for normal bundles of Riemannian submanifold embeddings]{Orthonormal frames for normal bundles \\ of Riemannian submanifold embeddings} \label{riem-frames}

We begin with the case of smooth Riemannian submanifold embeddings $\iota \colon \Lambda^{d-k} \hookrightarrow (M,g)$ with $d \geq 3$. For notational convenience, we will often view $\Lambda$ as a subset of $M$. We now review classical submanifold geometry.

\subsection{Classical submanifold geometry}

 At a point $p \in \Lambda$, every tangent space of the bulk space along $\Lambda$ can be decomposed according~to
\begin{align} \label{TM-decomp}
T_p M = T_p \Lambda \oplus N_p \Lambda ,
\end{align}
where $N_p \Lambda$ is the normal space at $p \in \Lambda$. The normal bundle over $\Lambda$ is defined as $N \Lambda := \sqcup_{p \in \Lambda} N_p \Lambda$.

As the pullback map is well defined, we can construct the \textit{intrinsic metric}, also sometimes called the \textit{first fundamental form}, via the pullback: $\bar{g} := \iota^* g$. Note that this metric is compatible with the decomposition above, and given an orthornormal frame $\{n_{\alpha}\}$ for $N \Lambda$, we may express this invariant in index notation according to $\bar{g}_{ab} = g_{ab}|_{\Lambda} - n_{a \alpha} n_{b \alpha}$. In the notation $n^a_{\alpha}$ (and $n_{a \alpha}$), the Latin index is an abstract index whereas the Greek letter indexes some particular element of the basis -- this will always be the case. Closely related to the first fundamental form is the projection operator given~by $\bar{g}^a_b := \delta^a_b - n^a_{\alpha} n_{b \alpha}$.

With a metric in hand, one may construct the Levi-Civita connection on $\Lambda$, denoted $\bar{\nabla}$. A~classical result of submanifold geometry is that this connection is directly obtainable from the Levi-Civita connection $M$ and the projection operator,
\[
\bar{\nabla}_{\bar v} \bar{u}^a = \bar{g}^a_b \nabla_{\iota_* \bar{v}} (\iota_* \bar{u})^b ,
\]
where $\bar{v},\bar{u} \in \Gamma(T \Lambda)$ and $\iota_* \bar{v}$, $\iota_* \bar{u}$ are their respective pushforwards. Note that we will often write~$\nabla^\top$ for the pullback of the bulk Levi-Civita connection -- this way, we may use abstract index notation and we do not need to carry around arbitrary vector fields. In a similar way, we use $\top$ both as an operator and as a superscript to denote the projection to the submanifold, e.g., we will write $P^\top$ to represent the projection of the Schouten tensor of $(M,g)$ to $\Lambda$.

If we instead project $\nabla_{\iota_* \bar{v}} (\iota_* \bar{u})$ along a unit normal vector, we obtain a classical result: the~\textit{second fundamental form}. Indeed, relative to an orthonormal frame on $N \Lambda$, we define
\[
\bar{u}^a \bar{v}^b \II_{ab \alpha} := -n_{a \alpha} \nabla_{\iota_* \bar{v}} (\iota_* \bar{u})^a .
\]
A short calculation shows that this can be expressed in a simpler way according to
\[
\II_{ab \alpha} = \bar{g}_{bc} \nabla^\top_a n^c_{\alpha} \in \Gamma\bigl(\odot^2 T^* \Lambda\bigr) .
\]
It is worth pointing out that in our calculations, we view $\II_{ab \alpha}$ as a tuple of 2-tensors depending on an orthonormal frame, as per the above section space. However when an orthonormal frame is not prescribed, the second fundamental form is still well defined, taking values in the normal bundle $N \Lambda$.
Accordingly, one may express the Levi-Civita connection of $M$ restricted to $\Lambda$ via the following identity:
\[
\nabla_{\iota_* \bar{v}} (\iota_* \bar{u}) = \bar{\nabla}_{\bar{v}} \bar{u} - \II_{\alpha}(\bar{u},\bar{v}) n_{\alpha} .
\]

The first and second fundamental forms are sufficient to completely describe an embedded hypersurface in $\mathbb{E}^d$ up to Euclidean motions. This is because, for $k=1$, the normal bundle is one-dimensional and hence (up to orientation) each normal space has a canonical basis given~by $n \in N_p \Lambda$, where $g(n,n) = 1$. (Note that as we work locally for the entirety of this article, we will neglect any issues that arise from orientations.) This rigidity effectively removes all geometric structure from the normal bundle. Unfortunately, this rigidity is lost entirely when we consider higher-codimension submanifolds. Indeed, while orthonormality can certainly be retained, there is (a priori) no uniquely preferred basis for $N \Lambda$, even locally. In what follows, we will find that there is sometimes a preferred basis up to constant ${\rm SO}(k)$ rotations -- but in general, it is not that clear that any preferred basis can be constructed. We will now provide a construction for the final piece of data required to describe embedded submanifolds up to Euclidean motions and also provide notation and language that will be necessary for later considerations.

Just as one can induce a connection on the tangent bundle to $\Lambda$, one may also induce on the normal bundle a connection
\[
D \colon \ \Gamma(N \Lambda) \rightarrow \Gamma(T^* \Lambda \otimes N \Lambda)
\]
given~by
\[
D_{\bar v} \xi = {\perp} (\nabla_{\iota_* \bar{v}} \xi) ,
\]
where $\perp$ is the orthogonal projection to $N \Lambda$, $\bar{v} \in \Gamma(T \Lambda)$ is any vector field along $\Lambda$ and $\xi \in \Gamma(N \Lambda)$ is any normal vector field. Note that this connection is compatible with the metric on the normal bundle given~by $\eta = g|_{\Lambda} - \iota^* g$. Then, for any frame $\{n_{\alpha}\}$ for $N \Lambda$, the connection coefficients are defined by
\[
D_a n^b_{\beta} = \omega_{\gamma a \beta} n^b_{\gamma} .
\]
Given any orthonormal frame $\{n_{\alpha}\}_{\alpha = 1}^k$ for the normal bundle $N \Lambda$, $\omega$ is given precisely by the so-called \textit{normal fundamental forms}~\cite{Spivak} (sometimes called the H{\'a}j{\'i}{\v{c}}ek one-form in the physics literature~\cite{hajicek1975, hajicek1973}):
\[
\omega_{\alpha a \beta} = \beta_{a \alpha \beta} := n^b_{\alpha} \nabla^\top_a n_{b \beta} .
\]
Note that we use $\beta$ only for orthonormal frames.

Associated to the normal bundle is the orthonormal frame bundle $F(N \Lambda)$, which is a principal ${\rm O}(k)$-bundle. So in particular, any two orthonormal frame fields are related to one another by the action of a section $m \in C^\infty(\Lambda, {\rm O}(k))$. Furthermore, the connection coefficients (and thus the normal fundamental forms) transform according to
\begin{align} \label{beta-transf}
\beta \mapsto m \beta m^{-1} + m \mathrm{d} \bigl(m^{-1}\bigr) ,
\end{align}
where adjacency is understood to imply matrix multiplication.
Note that for constant rotations, i.e., when $m$ is a constant matrix, $\beta$ transforms tensorially.

Because $\beta$ is not covariant under the group action, it is not a true invariant of the normal bundle (in that it does not transform linearly under the group action). Rather, its curvature 2-form
\begin{align} \label{norm-curv}
\mathcal{R} = \mathrm{d} \beta + [\beta, \beta]
\end{align}
does transform linearly under the group action and hence \textit{is} a true invariant of the normal bundle. In terms of the connection itself, the curvature 2-form, hereafter known as the \textit{normal curvature}, obeys the standard relation
\[
\mathcal{R}(\bar{u},\bar{v}) \xi = D_{\bar{u}} (D_{\bar{v}} \xi) - D_{\bar{v}} (D_{\bar{u}} \xi) - D_{[\bar{u},\bar{v}]} \xi
\]
for any $\bar{u},\bar{v} \in \Gamma(T \Lambda)$ and $\xi \in \Gamma(N \Lambda)$. The normal curvature is the final ingredient necessary to completely specify an embedded submanifold in $\mathbb{E}^d$ up to Euclidean motions.

To summarize, the three classical invariants of submanifold embeddings are $\bar{g}$, $\II$, and $\mathcal{R}$. It~is these tensors that appear in classical submanifold identities as described by Spivak~\cite{Spivak}, and given below for later use. First is the Gauss equation, relating the tangential component of the Riemann curvature tensor to that of the submanifold $\Lambda$,
\begin{align} \label{gauss-equation}
R_{abcd}^\top = \bar{R}_{abcd} + \II_{ad \alpha} \II_{bc \alpha} - \II_{bd \alpha} \II_{ac \alpha} ,
\end{align}
where $R^\top := \top (R)$, the projection of the Riemann curvature tensor to the submanifold. Similar notation will be used throughout.
Next is the Codazzi--Mainardi equation, relating another component of the Riemann curvature tensor to extrinsic invariants,
\begin{align} \label{codazzi}
R_{abc \alpha}^\top \= 2\bar{\nabla}_{[a} \II_{b]c \alpha} + 2\II_{c[b \beta} \beta_{a] \alpha \beta} .
\end{align}
Note that here, the right-hand side does not look explicitly tensorial -- however, a short calculation shows that it is indeed independent of the choice of orthonormal basis $\{n_{\alpha}\}$. Finally, we have the Ricci equation (which is trivial for hypersurfaces), relating another projection of the Riemann curvature to the normal curvature and the second fundamental form,
\[
R_{ab \alpha \beta}^\top = \mathcal{R}_{ab \alpha \beta} + \II^c_{a \beta} \II_{cb \alpha} - \II^c_{b \beta} \II_{ca \alpha} .
\]

As noted above, all of these equations are covariant under the action of a section $m \in C^\infty(\Lambda,{\rm O}(k))$. It would be desirable, in principle, to promote the normal fundamental forms to submanifold invariants so that covariance would be more explicitly seen. However, as $\beta$ does not transform covariantly under non-constant sections $m$, one might wish to restrict the structure group just to those constant sections by employing a geometric constraint on the normal frame bundle -- a procedure known in the physics literature as ``gauge-fixing''. Alas, in general, there is no preferred frame for the normal bundle fixed by a differential constraint. This is known as the Gribov ambiguity~\cite{gribov1978}. Furthermore, even global minimizing conditions of energy functionals depending on the normal connection are not necessarily unique~\cite{vanBaal1992}. In certain situations, canonical frames can be identified, at least up to global gauge transformations.

\subsection{Rotation minimizing frames}

In the simplest case, where the normal curvature vanishes, it is easy to show that a canonical frame (up to constant rotations) for the normal bundle can be constructed. We capture this result in the following proposition.
\begin{Proposition} \label{riem-canonical-frame}
Let $\Lambda^{d-k} \hookrightarrow (M,g)$ be a smooth submanifold embedding with $1 \leq k \leq d-1$. Then the normal curvature $\mathcal{R}$ vanishes if and only if there exists in a neighborhood around every point $p \in \Lambda$ an orthonormal frame for the normal bundle $N \Lambda$ with vanishing normal fundamental form $\beta$. Furthermore, this orthonormal frame is unique $($up to constants in $C^\infty(\Lambda,{\rm O}(k)))$.
\end{Proposition}
\begin{proof}
It follows trivially from equation~(\ref{norm-curv}) that if there exists an orthonormal frame such that $\beta = 0$, then $\mathcal{R} = 0$. Now suppose that $\mathcal{R} = 0$. We wish to show that there exists a~frame such that $\beta = 0$. This follows from a~standard result: namely, if the bundle $N \Lambda \rightarrow \Lambda$ has a~connection $D$ with vanishing curvature, then there always exists locally a frame $\{n_{\alpha}\}$ consisting of parallel sections, namely $D n_{\alpha} = 0$. But then by the Gram--Schmidt procedure, it follows that (along $\Lambda$) there exists an orthonormal frame $\{n_{\alpha}\}$ consisting of parallel sections, and so we have that for that choice of frame, $\beta = 0$.

Now suppose $\mathcal{R} = 0$ and that two frames $\{n_{\alpha}\}$ and $\{n'_{\alpha}\}$ have vanishing normal fundamental forms. Because they are both orthonormal frames, there exists a (possibly non-constant) $m \in C^\infty(\Lambda, {\rm O}(k))$ such that $n'_{\alpha} = m_{\alpha \beta} n_{\beta}$. It then follows from equation~(\ref{beta-transf}) that
\[
0 = m_{\alpha \gamma} \nabla_a m_{\beta \gamma} .
\]
But clearly this identity only holds if $\nabla_a m_{\beta \gamma} = 0$. Thus we have that $m$ is a constant section of~${{\rm O}(k) \rightarrow \Lambda}$.
\end{proof}

It is thus useful to think of this canonical orthonormal frame as, in some sense, ``maximally parallel''. From a geometrical perspective, the associated $\beta$ for a given orthonormal frame $\{n_{\alpha}\}$ can be thought of as a generalization of the torsion of a spacecurve, insofar as torsion describes how one normal vector rotates into another normal vector while moving along a spacecurve. In~fact, for $\Lambda^1 \hookrightarrow \mathbb{R}^3$, it is straightforward to show that in the Frenet frame $\{T, N, B\}$ with torsion~$\tau$, we have that
\begin{align*}
\beta_{a B N} = \tau .
\end{align*}
It is because of this analogy that those orthonormal frames for the normal bundle that fix $\beta = 0$ are sometimes called rotation minimizing~\cite{RMFs} or Bishop frames~\cite{bishop}. We provide a formal definition below:
\begin{Definition}
Let $\Lambda^{d-k} \hookrightarrow (M,g)$ be a Riemannian submanifold embedding and let $\{n_{\alpha}\}$ be an orthonormal frame for $N \Lambda$. Then we say that $\{n_{\alpha}\}$ is a \textit{rotation minimizing frame} when, for each $\alpha \in \{1, \ldots, k\}$ and $\bar{v} \in \Gamma(T \Lambda)$, we have that
\[
\nabla_{\iota_* \bar{v}} n_{\alpha} \in \Gamma(T \Lambda) .
\]
\end{Definition}

By construction, rotation minimizing frames are precisely the canonical frames described in Proposition~\ref{riem-canonical-frame}, as captured in the following lemma:
\begin{Lemma} \label{rot-min-frame}
Let $\Lambda^{d-k} \hookrightarrow (M,g)$ have a rotation minimizing frame $\{n_\alpha\}$. Then modulo the action of constant sections of ${\rm O}(k) \rightarrow \Lambda$, this frame is the canonical frame described in Proposition~$\ref{riem-canonical-frame}$.
\end{Lemma}
\begin{proof}
Let $\{n_{\alpha}\}$ be a rotation minimizing frame for $N \Lambda$. From the definition, we have that for any $\bar{v} \in \Gamma(T \Lambda)$, $\nabla_{\iota_* \bar{v}} n_{\alpha} \in \Gamma(T \Lambda)$. But by the orthogonal decomposition $T M|_{\Lambda} = T \Lambda \oplus N \Lambda$, it~follows that $n_{a \beta} \nabla_{\iota_* \bar{v}} n^a_{\alpha} = 0$. But this implies that $\bar{v}^a \beta_{a \beta \alpha} = 0$ for all $\alpha,\beta \in \{1, \ldots, k\}$ and any choice of $\bar{v}$. This is precisely the frame described in Proposition~\ref{riem-canonical-frame}.
\end{proof}

We thus have the following result for spacecurves.
\begin{Corollary}
Let $\Lambda^1 \hookrightarrow \bigl(M^d,g\bigr)$ be a curve embedded in a Riemannian manifold. Then there exists a unique $($up to constant rotations$)$ rotation minimizing frame for $N \Lambda$.
\end{Corollary}
\begin{proof}
Because $\Lambda$ is 1-dimensional and the normal curvature is a $\operatorname{End}(N \Lambda)$-valued 2-form, we have that $\mathcal{R}$ necessarily vanishes. Hence from Proposition~\ref{riem-canonical-frame} and Lemma~\ref{rot-min-frame}, there exists a~unique rotation minimizing frame.
\end{proof}

\subsection{Gauge-fixed frames} \label{gauge-fixed}
When the normal curvature is non-vanishing, rotation minimizing frames do not exist. However, that does not necessarily preclude the existence of a preferred frame determined by some condition, differential or otherwise.

We begin by considering the simplest non-trivial case, where $k=2$. In that setting, we can directly parametrize sections of ${\rm SO}(2) \rightarrow \Lambda$ by a single function $\theta \in C^\infty S^1$ so that if $m \in C^\infty(\Lambda, {\rm SO}(2))$, then we have that
\[
m_{\alpha \beta} = \cos \theta \delta_{\alpha \beta} - \sin \theta \epsilon_{\alpha \beta} .
\]
Given a orthonormal frame $\{n_{\alpha} \}$, the normal fundamental form is specified entirely by a single 1-form $\beta_a$. Under the group action, we find that this function transforms according to
\[
\beta_a \mapsto \beta_a + \bar{\nabla}_a \theta .
\]

Because this gauge transformation is sufficiently simple, we may be able to find a unique solution to some differential gauge-fixing condition. Consider the Coulomb gauge $\bar{\nabla} \csdot \beta = 0$. Fixing this gauge is equivalent to finding $\theta$ such that
\[
\bar{\Delta} \theta = - \bar{\nabla}^a \beta_a .
\]
In general, any global solution to this equation will require a global frame and hence requires that $F(N \Lambda)$ is trivial. In that case, if $\Lambda$ is compact and $\int_{\Lambda} \bar{\nabla} \csdot \beta = 0$, then there exists $\theta$ that solves the Poisson equation globally and allows the construction of a preferred frame (up to a~constant) $\big\{m_{\alpha \beta} n^a_{\beta}\big\}$. Trivially, then, if $\Lambda$ is closed such a preferred frame exists.

More generally, for $k \geq 2$, Uhlenbeck showed~\cite{uhlenbeck1982} that when the normal curvature is sufficiently small in a trivializing neighborhood $U \subset \Lambda$ around a point $p$, there exists a gauge transformation $m \in C^\infty(U, {\rm O}(k))$ such that the basis $\big\{m_{\alpha \beta} n^a_{\beta}\big\}$ has normal fundamental forms that satisfy $\bar{\nabla} \csdot \beta_{\alpha \beta} = 0$ and $u \csdot \beta = 0$ on $\partial U$, where $u$ is the normal to the boundary $\partial U$.

While the existence of a ``clean'' gauge choice can sometimes be taken for granted, uniqueness is much harder to demonstrate -- and indeed, it is not always present (see the discussion about the Gribov ambiguity above). Regardless, we can still study the extrinsic geometry of submanifold embeddings equipped with an orthonormal basis for the normal bundle. When a global orthonormal frame exists for $N \Lambda$, it is called a \textit{parallelization} of $N \Lambda$. Since we work locally, we can always assume that an orthonormal frame is specified for $N \Lambda$. We call such an embedding equipped with an orthonormal frame for the normal bundle a \textit{parallelized $($Riemannian$)$ submanifold embedding}, specified by the pair $(\Lambda \hookrightarrow (M,g), \{n_{\alpha}\})$.

\section[Defining maps for parallelized Riemannian submanifold embeddings]{Defining maps for parallelized Riemannian submanifold\\ embeddings} \label{unit-def-funcs}

Given a Riemannian submanifold embedding $\iota \colon \Lambda^{d-k} \hookrightarrow \bigl(M^d,g\bigr)$, we are interested in studying the embedding by examining the jets of a geometrically-adapted \textit{defining map} $s_{\alpha} \colon M \rightarrow \mathbb{R}^k$, which is a $k$-tuple of functions such that for any $p \in \Lambda$, we have that $s_{\alpha}(p) = 0$ for every $\alpha \in \{1, \ldots, k\}$ and ${\rm d}s_{1}(p) \wedge \cdots \wedge {\rm d}s_{k}(p) \neq 0$. That such families exist is well known, at least locally: for every $p \in \Lambda$, there exists a neighborhood $U \subset M$ such that $U \cap \Lambda$ is the zero locus of some defining map $s_{\alpha} \colon U \rightarrow \mathbb{R}^k$. However, there are many such defining maps, and in order for the jets of such a defining map to contain information about the embedding $\iota$, it must be determined canonically by the geometry. To realize this canonical defining map, we consider the \text{Gram matrix} $G_{\alpha \beta} := g^{-1}({\rm d}s_{\alpha},{\rm d}s_{\beta})$ of a given defining map $s_{\alpha}$. In the case of hypersurfaces, we canonically choose $G = \Id$, which is equivalent to the condition that $|{\rm d}s|^2 = 1$. Similarly, in the case of a higher-codimension submanifold we would like to achieve $G= \Id$. Because we will work locally, it is sufficient for our purposes to find a defining map $s_{\alpha}$ such that $G_{\alpha \beta} = \delta_{\alpha \beta} + \mathcal{O}(s^\infty)$, meaning we wish to find a defining map such that for any $m \in \mathbb{Z}_{\geq 1}$, we have that $G_{\alpha \beta} = \delta_{\alpha \beta} + \mathcal{O}(s^m)$. We call such a defining map a \textit{unit defining map}. Going forward, we will abuse notation by also denoting $g^{-1}({\rm d}s_{\alpha},\cdot)$ by $n_{\alpha}$. Note that, when clarification is necessary, we will decorate $n_{\alpha}$ with abstract indices (such as $n^a_{\alpha}$ and $n_{a \alpha}$) to indicate the whether we are referring to the vector or its metric dual.

\subsection{First order} \label{fo-Riem}
Given the parallelization $\{n_{\alpha}\}$ of $N \Lambda$ and any defining map $\hat{s}_{\alpha}$ for $\Lambda$, there exists a linear combination of components of $\hat{s}_{\alpha}$ denoted by $s_{\alpha}$ such that $g^{-1}({\rm d}s_{\alpha},\cdot) \= n_{\alpha}$. This linear combination provides a unique $s_{\alpha}$ modulo terms that are at least second order in a defining function. The behavior of ${\rm d}s_{\alpha}$ away from $\Lambda$, however, is harder to control. Our goal is to obtain a canonical family of unit defining functions adapted to our parallelization satisfying $g^{-1}({\rm d}s_{\alpha},{\rm d}s_{\beta}) = \delta_{\alpha \beta}$. However, this will not be possible in general. Instead, we will choose $s_{\alpha}$ to, in some sense, get as close to the identity Gram matrix condition as possible.
We begin with an iterative procedure. Because the parallelization of $N \Lambda$ is defined to be orthonormal, we can certainly obtain $G_{\alpha \beta} \= \delta_{\alpha \beta}$, so we assume this going forward. This constraint defines a unique family of defining maps that are aligned to the parallelization of $N \Lambda$ and whose elements have a difference given~by \smash{$s'_{\alpha} - s_{\alpha} = \mathcal{O}\bigl(s^2\bigr)$}. For any of these defining maps, it thus follows that
\begin{align} \label{0-order}
G_{\alpha \beta} = \delta_{\alpha \beta} + F^{(1)}_{\alpha \beta \gamma} s_{\gamma} ,
\end{align}
where $F^{(1)}_{\alpha \beta \gamma}$ represents an array of smooth functions on $M$. Our goal will be to pare down this family of defining maps by controlling \smash{$F^{(1)}$}.

However, note that equation~(\ref{0-order}) does not uniquely specify \smash{$F^{(1)}_{\alpha \beta \gamma}$} as an array of functions on $C^\infty M$. To see this, consider the mapping
\[
F^{(1)}_{\alpha \beta \gamma} \mapsto F^{(1)}_{\alpha \beta \gamma} + F^{(1,1)}_{\alpha \beta \gamma \delta} s_{\delta} ,
\]
where \smash{$F^{(1,1)}$} is any array of functions so that \smash{$F^{(1,1)}_{\alpha \beta \gamma \delta} = F^{(1,1)}_{\alpha \beta [\gamma \delta]}$}. This mapping clearly preserves equation~(\ref{0-order}), and so we see that \smash{$F^{(1)}$} is only uniquely determined along $\Lambda$. Nonetheless, we \textit{may} choose \smash{$F^{(1)}$} in a collar neighborhood of $\Lambda$ in a non-canonical way. Going forward, we will assume such a choice has been made; later, we will show that this choice actually yields a unique result.

Because of the symmetry of $G_{\alpha \beta}$, we have that \smash{$F^{(1)}_{\alpha \beta \gamma} = F^{(1)}_{(\alpha \beta) \gamma}$}, and thus as a representation of $S_k$, we can see that the tensor structure of \smash{$F^{(1)}$} in the normal bundle takes the form
\[
\raisebox{3pt}{\scalebox{0.4}{\ydiagram{2,1}}} \; \oplus \; \scalebox{0.4}{\ydiagram{3}} .
\]
For what follows, we will abuse notation and write, for example, that $F^{(1)} \in C^\infty (M,\raisebox{3pt}{\scalebox{0.2}{\ydiagram{2,1}}} \; \oplus \; \raisebox{2pt}{\scalebox{0.2}{\ydiagram{3}}})$.

Now observe that under the mapping
\[
s_{\alpha} \mapsto \tilde{s}_{\alpha} := s_{\alpha} + A^{(1)}_{\alpha \gamma_1 \gamma_2} s_{\gamma_1} s_{\gamma_2} ,
\]
where $A^{(1)}$ is a list of smooth functions on $M$, we have that the zeroth order part of the Gram matrix condition given in equation~(\ref{0-order}) is preserved while changing $F^{(1)}$. Indeed, by inspection we may consider only those \smash{$A^{(1)}$} satisfying \smash{$A^{(1)}_{\alpha \gamma_1 \gamma_2} = A^{(1)}_{\alpha (\gamma_1 \gamma_2)}$}, and so $A^{(1)} \in C^\infty (M , \raisebox{3pt}{\scalebox{0.2}{\ydiagram{2,1}}} \; \oplus \; \raisebox{2pt}{\scalebox{0.2}{\ydiagram{3}}})$. Evidently, it may be possible to choose \smash{$A^{(1)}$} so that it cancels \smash{$F^{(1)}$}.

Differentiating the above display and using the definition of the Gram matrix for the new defining map (combined with equation~\eqref{0-order}), we obtain
\begin{align}
\tilde{G}_{\alpha \beta} &{}= \delta_{\alpha \beta}+ \bigl(F^{(1)}_{\alpha \beta \gamma} + 4 A^{(1)}_{(\alpha \beta) \gamma}\bigr) s_{\gamma} \nonumber
\\
&\quad{}+ \bigl(4 A^{(1)}_{\alpha (\rho \gamma_1)} A^{(1)}_{\beta (\rho \gamma_2)} + 2 n_{a (\alpha} \nabla^a A^{(1)}_{\beta) (\gamma_1 \gamma_2)} + 4 A^{(1)}_{\alpha (\rho \gamma_1)} F^{(1)}_{\beta \rho \gamma_2} \bigr) s_{\gamma_1} s_{\gamma_2}
\nonumber\\
&\quad{}+ \bigl(4A^{(1)}_{\alpha (\rho \gamma_1)} \nabla_{n_{\rho}} A^{(1)}_{\beta \gamma_2 \gamma_3} + 4 A^{(1)}_{\alpha (\omega_1 \gamma_2)} A^{(1)}_{\beta (\omega_2 \gamma_2)} F^{(1)}_{\omega_1 \omega_2 \gamma_3} \bigr) s_{\gamma_1} s_{\gamma_2} s_{\gamma_3} \nonumber\\
&\quad{}+ \bigl(\nabla^a A^{(1)}_{\alpha (\gamma_1 \gamma_2)}\bigr) \bigl(\nabla_a A^{(1)}_{\beta (\gamma_3 \gamma_4)} \bigr) s_{\gamma_1} s_{\gamma_2} s_{\gamma_3} s_{\gamma_4} ,
\label{1-order-correction}
\end{align}
where the above is symmetric in $\alpha \beta$ (even when the round brackets are not displayed).
Observe that in order to fix $A^{(1)}$, we must express it in terms of $F^{(1)}$. To that end, we fix an ansatz: \smash{$A^{(1)}_{\alpha \beta \gamma} = a F^{(1)}_{\alpha \beta \gamma} + b F^{(1)}_{\gamma \alpha \beta} + c F^{(1)}_{\beta \gamma \alpha}$} and then solve explicitly; these are the only three terms~possible~as \smash{$F^{(1)}_{\alpha \beta \gamma} = F^{(1)}_{(\alpha \beta) \gamma}$}. We find that
\[
A^{(1)}_{\alpha \beta \gamma} = -\frac{1}{4}\bigl(F^{(1)}_{\alpha \beta \gamma} + F^{(1)}_{\gamma \alpha \beta} - F^{(1)}_{\beta \gamma \alpha}\bigr)
\]
sets the first order term in equation~\eqref{1-order-correction} to zero identically, and thus we have found $s_{\alpha}$ such that
\[
G_{\alpha \beta} = \delta_{\alpha \beta} + \mathcal{O}\bigl(s^2\bigr) ,
\]
fixing $s_{\alpha}$ to second order in a way dependent on the original choice of extension of $F^{(1)}$.

With this result, we find a useful identity for the normal fundamental form $\beta$ associated with the parallelization of $N \Lambda$:
\begin{align}
n^b_{\alpha} \nabla_a n_{b \beta} \={}& n^b_{[\alpha} \nabla_a n_{b \beta]} + n^b_{(\alpha} \nabla_a n_{b \beta)}
\= n^b_{[\alpha} \bigl(\nabla_a^\top + n_{a \gamma} n^c_{\gamma} \nabla_c \bigr) n_{b \beta]} + \frac{1}{2} \nabla_a n^b_{(\alpha} n_{b \beta)} \nonumber \\
\={}& \beta_{a \alpha \beta} + n_{a \gamma} n^c_{\gamma} n^b_{[\alpha} \nabla_c n_{b \beta]}
\= \beta_{a \alpha \beta} + n_{a \gamma} n^c_{\gamma} n^b_{[\alpha} \nabla_b n_{c \beta]} \nonumber \\
\={}& \beta_{a \alpha \beta} - n_{a \gamma} n^b_{[\alpha} n^c_{\beta]} \nabla_b n_{c \gamma}
\= \beta_{a \alpha \beta} . \label{ndn-Riem}
\end{align}
In the above computation, we used that when $n_{\alpha} = {\rm d} s_{\alpha}$, we evidently have that $\nabla_{[a} n_{b]\beta} = 0$. Observe from this calculation that \smash{$[n_{\alpha},n_{\beta}]^a \= 2 \beta^a_{\alpha \beta}$}. Further, it thus follows that
\begin{align} \label{nabla-n}
\nabla_a n_{b \beta} \= \II_{ab \beta} + n_{b \alpha} \beta_{a \alpha \beta} + n_{a \alpha} \beta_{b \alpha \beta} .
\end{align}

\subsection{Second order} 
Section~\ref{fo-Riem} established that for a given parallelized Riemannian submanifold embedding, there exists a unique family of defining maps related by $s'_{\alpha} - s_{\alpha} = \mathcal{O}\bigl(s^3\bigr)$ such that
\begin{align} \label{1-order}
G_{\alpha \beta} = \delta_{\alpha \beta} + F^{(2)}_{\alpha \beta \gamma_1 \gamma_2} s_{\gamma_1} s_{\gamma_2} .
\end{align}
As in that section, our goal is to pare down this family by controlling $F^{(2)}$. From equation~(\ref{1-order-correction}), $F^{(2)}$ is uniquely determined in a collar neighborhood of $\Lambda$ given the initial choice of $F^{(1)}$. Our goal will be to adjust $s_{\alpha}$ so that we may remove as much of $F^{(2)}$ as possible.

Using the same reasoning as in Section~\ref{fo-Riem}, we have that $F^{(2)} \in C^\infty (M, \raisebox{3pt}{\scalebox{0.2}{\ydiagram{2,2}}} \; \oplus \raisebox{3pt}{\scalebox{0.2}{\ydiagram{3,1}}} \; \oplus \; \raisebox{1.5pt}{\scalebox{0.2}{\ydiagram{4}}})$. So, we attempt to fix the family of defining functions to cancel $F^{(2)}$ by applying the map
\[
s_{\alpha} \mapsto \tilde{s}_{\alpha} := s_{\alpha} + A^{(2)}_{\alpha \gamma_1 \gamma_2 \gamma_3} s_{\gamma_1} s_{\gamma_2} s_{\gamma_3} ,
\]
which clearly preserves the structure of the Gram matrix condition~(\ref{1-order}) while possibly changing~$F^{(2)}$. However, note that $A^{(2)} \in C^\infty (M, \raisebox{3pt}{\scalebox{0.2}{\ydiagram{3,1}}} \; \oplus \; \raisebox{1.5pt}{\scalebox{0.2}{\ydiagram{4}}})$. Because $A^{(2)}$ and $F^{(2)}$ can live in different representations of $S_k$, it is not always possible to construct $A^{(2)}$ such that $F^{(2)} \mapsto 0$. Nonetheless, one can always find $A^{(2)}$ so that
\[
F^{(2)}_{\alpha \beta \gamma_1 \gamma_2} \mapsto \tilde{F}^{(2)}_{(\alpha \beta)(\gamma_1 \gamma_2)} \in C^\infty (M,\raisebox{3pt}{\scalebox{0.2}{\ydiagram{2,2}}}).
\]

From the above construction, for a given parallelization of $N \Lambda$, there always exists at least one defining map such that
\begin{align*} 
G_{\alpha \beta} = \delta_{\alpha \beta} + F^{(2)}_{\alpha \beta \gamma_1 \gamma_2} s_{\gamma_1} s_{\gamma_2} ,
\end{align*}
where $F^{(2)} \in C^\infty (M,\raisebox{3pt}{\scalebox{0.2}{\ydiagram{2,2}}})$ is entirely determined by the parallelization of $N \Lambda$ and the choice of~$F^{(1)}$. We call defining maps that are in this infinite family \textit{associated defining maps}. Along~$\Lambda$, $F^{(2)}|_{\Lambda}$ is a first obstruction to producing a holonomic basis for the extension of the normal bundle. In fact, $F^{(2)}|_{\Lambda}$ can be computed in terms of tensors specified by the parallelization. We now carry out this computation, and going forward we will work with the family of associated defining maps to the parallelization.

Now, observe that
\[
P_{\scalebox{0.2}{\ydiagram{2,2}}}\, n^a_{\gamma_1} n^b_{\gamma_2} \nabla_a \nabla_b G_{\alpha \beta} \= 2 F^{(2)}_{\alpha \beta \gamma_1 \gamma_2} ,
\]
where $P_{\scalebox{0.2}{\ydiagram{2,2}}}$ is the projector onto the $\raisebox{3pt}{{\scalebox{0.2}{\ydiagram{2,2}}}}$ component of a rank-4 tensor in $S^k$. In calculations that will follow, it is useful to have explicit formulas for relevant projection operators. To that end, we find that
\begin{align*}
P_{\scalebox{0.2}{\ydiagram{2,2}}}\, F_{\alpha \beta \gamma \delta} = \frac{1}{12} \bigl(&F_{\alpha \beta \gamma \delta} - F_{\gamma \beta \alpha \delta} - F_{\alpha \delta \gamma \beta} + F_{\gamma \delta \alpha \beta} + F_{\beta \alpha \gamma \delta} - F_{\gamma \alpha \beta \delta} - F_{\beta \delta \gamma \alpha} + F_{\gamma \delta \beta \alpha} \\
&+ F_{\alpha \beta \delta \gamma}- F_{\delta \beta \alpha \gamma} - F_{\alpha \gamma \delta \beta} + F_{\delta \gamma \alpha \beta} + F_{\beta \alpha \delta \gamma} - F_{\delta \alpha \beta \gamma} - F_{\beta \gamma \delta \alpha} + F_{\delta \gamma \beta \alpha}\bigr) ,
\end{align*}
and
\begin{align*}
P_{\scalebox{0.2}{\ydiagram{3,1}}}\, F_{\alpha \beta \gamma \delta} = \frac{1}{8}\bigl(&F_{\alpha \beta \gamma \delta} - F_{\delta \beta \gamma \alpha} + F_{\beta \alpha \gamma \delta} - F_{\delta \alpha \gamma \beta} + F_{\gamma \beta \alpha \delta} - F_{\delta \beta \alpha \gamma} + F_{\alpha \gamma \beta \delta} - F_{\delta \gamma \beta \alpha}
\\&+ F_{\beta \gamma \alpha \delta}- F_{\delta \gamma \alpha \beta} + F_{\gamma \alpha \beta \delta} - F_{\delta \alpha \beta \gamma} \bigr) .
\end{align*}
We can now compute directly
\begin{align*}
P_{\scalebox{0.2}{\ydiagram{2,2}}}\, n^a_{\gamma_1} n^b_{\gamma_2} \nabla_a \nabla_b n^c_{\alpha} n_{c \beta} ={}& P_{\scalebox{0.2}{\ydiagram{2,2}}}\, n^{(a}_{\gamma_1} n^{b)}_{\gamma_2} \bigl[ 2 (\nabla n_{\alpha})^c_{(a} (\nabla n_{\beta})_{b) c} + 2 n^c_{(\alpha|} \nabla_a \nabla_b n_{c| \beta)} \bigr] \\
\={}& -2 \beta_{c \alpha (\gamma_1 } \beta^c_{\gamma_2) \beta} - \frac{2}{3} R_{\gamma_1 (\alpha \beta) \gamma_2} .
\end{align*}
Putting this together, we have that
\begin{align} \label{Riem-2o-obstruction}
F^{(2)}_{\alpha \beta \gamma_1 \gamma_2} \= - \beta_{c \alpha (\gamma_1 } \beta^c_{\gamma_2) \beta} - \frac{1}{3} R_{\gamma_1 (\alpha \beta) \gamma_2} .
\end{align}
It is important to note that $F^{(2)}$ is (as yet) dependent on the choice of extension of $F^{(1)}|_{\Lambda}$. However, what this computation verifies is that this choice of extension yields a unique obstruction $F^{(2)}|_{\Lambda}$ that depends only on the parallelization of $N \Lambda$ and not on some choice of how the parallelization or $F^{(1)}|_{\Lambda}$ extends off $\Lambda$.

Thus, the obstruction in equation~(\ref{Riem-2o-obstruction}) is a true invariant of the parallelized Riemannian submanifold embedding. As such, it may be of independent interest to study those structures for which this obstruction vanishes. While there is indeed an obstruction, we may proceed regardless.

\subsection{Higher orders} \label{ho-Riem}
We now seek to generalize the result of the previous subsection to all orders. As noted, at second order $F^{(2)}|_{\Lambda}$ is unique, despite the choice made of extension of $F^{(1)}$ off $\Lambda$. We would like to show that this phenomenon continues. This result is captured by the following theorem.

\begin{Theorem} \label{general-ho-ext}
Let $\Lambda^{d-k} \hookrightarrow \bigl(M^d,g\bigr)$ be parallelized by $\{n_{\alpha}\}$. Then, in a collar neighborhood of~$\Lambda$, there exists formally a unique defining map determined by $\{n_{\alpha}\}$ such that
\begin{align} \label{window-expansion}
G_{\alpha \beta} = \delta_{\alpha \beta} + F^{(2)}_{\alpha \beta \gamma_1 \gamma_2} s_{\gamma_1} s_{\gamma_2} + F^{(3)}_{\alpha \beta \gamma_1 \gamma_2 \gamma_3} s_{\gamma_1} s_{\gamma_2} s_{\gamma_3} + \cdots,
\end{align}
where each $F^{(i)}$ is canonically determined to all orders and lives in a generalized window tensor representation of $S^k$.
\end{Theorem}
For clarity, a generalized window tensor representation of $S^k$ corresponds to a Young diagram with a first row of length at least two, and second row of length two.
\begin{proof}
We first show that one may find a defining map that satisfies equation~(\ref{window-expansion}) via induction. To do so, we begin by fixing an arbitrary extension of $F^{(1)}|_{\Lambda}$ as described in Section~\ref{fo-Riem}. Given such an extension, the base case is handled by the previous subsection, so that there exists a defining map so that \smash{$G = \delta + \mathcal{O}\bigl(s^2\bigr)$}. Now, suppose that there exists a defining map such that
\[
G_{\alpha \beta} = \delta_{\alpha \beta} + F^{(2)}_{\alpha \beta \gamma_1 \gamma_2} s_{\gamma_1} s_{\gamma_2} + \cdots + F^{(m)}_{\alpha \beta \gamma_1 \cdots \gamma_m} s_{\gamma_1} \cdots s_{\gamma_m} + F^{(m+1)}_{\alpha \beta \gamma_1 \cdots \gamma_{m+1}} s_{\gamma_1} \cdots s_{\gamma_{m+1}} ,
\]
where for each $2 \leq i \leq m$, $F^{(i)}$ occupies a generalized window tensor representation. Now let
\[
s_{\alpha} \mapsto \tilde{s}_{\alpha} := s_{\alpha} + A^{(m+1)}_{\alpha \gamma_1 \cdots \gamma_{m+2}} s_{\gamma_1} \cdots s_{\gamma_{m+2}} .
\]
Then
\[
\tilde{G}_{\alpha \beta} = G_{\alpha \beta} + 2(m+1) A^{(m+1)}_{(\alpha \beta) \gamma_1 \cdots \gamma_{m+1}} s_{\gamma_1} \cdots s_{\gamma_{m+1}} + \mathcal{O}\bigl(s^{m+2}\bigr) .
\]
From the representation theory of the normal bundle, we find that $A^{(m+1)}$ may thus be chosen in a collar neighborhood of $\Lambda$ so that $F^{(m+1)}$ occupies the generalized window tensor representation. This completes the induction, and thus given an arbitrary extension of $F^{(1)}|_{\Lambda}$ off $\Lambda$, we have uniquely determined a defining map $s_{\alpha}$ that satisfies equation~(\ref{window-expansion}).

Next, we must show that this defining map is independent of the choice of $F^{(1)}$. Let $F^{(1)}$ and $F^{(1)'}$ be distinct extensions of $F^{(1)}|_{\Lambda}$. Then we may produce unique defining maps $s$ and $s'$ respectively associated to each of these choices as described above. Now because the differentials of the defining maps necessarily agree on the submanifold, we must have that
\[
s'_{\alpha} = s_{\alpha} + M_{\alpha \gamma_1 \gamma_2} s_{\gamma_1} s_{\gamma_2}
\]
for some $M_{\alpha \gamma_1 \gamma_2} \in C^\infty M \times (\raisebox{3pt}{\scalebox{0.2}{\ydiagram{2,1}}} \; \oplus \; \raisebox{2pt}{\scalebox{0.2}{\ydiagram{3}}})$. We may now compare the Gram matrices,
\[
G'_{\alpha \beta} = G_{\alpha \beta} + 2 M_{(\alpha \beta) \gamma_1} s_{\gamma_1} + \mathcal{O}\bigl(s^2\bigr) .
\]
However, observe that $G'_{\alpha \beta}$ and $G_{\alpha \beta}$ must be the identity to second order, and so~\smash{$M_{(\alpha \beta) \gamma} \= 0$}. So because \smash{$M_{\alpha \beta \gamma} \= M_{[\alpha \beta] \gamma}$},
\[
M_{\alpha \beta \gamma} \= -M_{\beta \alpha \gamma} \= -M_{\beta \gamma \alpha} \= M_{\gamma \beta \alpha} \= M_{\gamma \alpha \beta} \= -M_{\alpha \gamma \beta} \= -M_{\alpha \beta \gamma} .
\]
Thus, we have that $M_{\alpha \beta \gamma} \= 0$.

So this means that
\[
s'_{\alpha} = s_{\alpha} + M_{\alpha \gamma_1 \gamma_2 \gamma_3} s_{\gamma_1} s_{\gamma_2} s_{\gamma_3} .
\]
By the same considerations, we have that
\[
G'_{\alpha \beta} - G_{\alpha \beta} = 4 M_{(\alpha \beta) \gamma_1 \gamma_2} s_{\gamma_1} s_{\gamma_2} .
\]
However, at second order the left-hand side necessarily occupies a window tensor representation~-- but $M$ is symmetric in three indices, and thus cannot contribute to this tensor representation. Therefore, it must vanish along $\Lambda$, and hence $s' - s = \mathcal{O}\bigl(s^4\bigr)$. By a similar argument, we may apply induction to find that $s'_{\alpha} = s_{\alpha}$ to arbitrarily high order. Thus, the defining map is independent of the choice of the extension of $F^{(1)}|_{\Lambda}$.
\end{proof}

As the defining map constructed above is unique and fixed for a given parallelization, we call any such defining maps \textit{canonical defining maps}. Note the canonical defining map for a given parallelization is also an associated defining map.

\begin{Remark} \label{curve-invariants}
Note that for one-dimensional submanifolds, a geometrically-determined canonical parallelization exists, i.e., the rotation minimizing frame. As such, all of the obstructions produced above (while still generically non-vanishing) depend only on the embedding and not on a choice of parallelization.
\end{Remark}

\begin{Remark}
A related construction can be performed using Fermi normal coordinates -- for explicit details, see~\cite{Gray,GrayBook}. Clearly, given a Fermi normal coordinatization, one can use our construction to build a canonical defining map. However, because our construction necessarily requires that the Gram matrix coefficients have certain symmetries, it is not clear that the curves generated by these defining maps are geodesics. On the other hand, given a (canonical) defining map, it is not clear how this determines a coordinate system in general.
\end{Remark}

When the submanifold embedding is particularly geometrically simple, we may do better. To see this, we first need a technical lemma.
\begin{Lemma} \label{lower-order-derivs}
Suppose $\Lambda^{d-k} \hookrightarrow (M,g)$ has vanishing normal curvature, that $(M,g)$ is flat, and that there exists a defining map satisfying $\beta = 0$ and $G_{\alpha \beta} = \delta_{\alpha \beta} + \mathcal{O}(s^m)$ for some $m \geq 2$. Then
\[
n^{a_1}_{\gamma_1} \cdots n^{a_j}_{\gamma_j} \nabla_{a_1} \cdots \nabla_{a_j} n^b_{\alpha} \= 0
\]
for any $1 \leq j \leq m-1$.
\end{Lemma}

\begin{proof}
We prove by induction on $j$. For that, we first examine the case where $j = 1$. We must~show that \smash{$n^a_{\gamma} \nabla_a n^b_{\alpha} \= 0$}. However, because $m \geq 2$, equation~(\ref{nabla-n}) holds, and thus this condition is equivalent to the condition that $\beta = 0$, which is satisfied by hypothesis. Note that this entirely covers the case where $m = 2$, so going forward we assume that $m \geq 3$.

Now fix some $1 \leq j \leq m-2$, and suppose that the equation in the lemma holds for all $1 \leq \ell \leq j$. We wish to show that this implies that the equation holds for $j+1$. We compute
\begin{align*}
n^{a_1}_{\gamma_1} \cdots n^{a_{j+1}}_{\gamma_{j+1}} \nabla_{a_1} \cdots \nabla_{a_{j+1}} n^b_\alpha \=&-n^{a_1}_{\gamma_1} \cdots n^{a_{j}}_{\gamma_{j}} \bigl(\nabla_{a_1} n^{a_{j+1}}_{\gamma_{j+1}}\bigr) \nabla_{a_2} \cdots \nabla_{a_{j+1}} n^b_{\alpha} \\&{}+ n^{a_1}_{\gamma_1} \cdots n^{a_{j}}_{\gamma_{j}} \nabla_{a_1} n^{a_{j+1}}_{\gamma_{j+1}} \nabla_{a_2} \cdots \nabla_{a_{j+1}} n^b_\alpha \\
\={}& \cdots \\
\={}& n^{a_1}_{\gamma_1} \cdots n^{a_{j}}_{\gamma_{j}} \nabla_{a_1} \cdots \nabla_{a_j} n^{a_{j+1}}_{\gamma_{j+1}} \nabla_{a_{j+1}} n^b_{\alpha} \\
\={}& n^{a_1}_{\gamma_1} \cdots n^{a_{j}}_{\gamma_{j}} \nabla_{a_1} \cdots \nabla_{a_j} n^{a_{j+1}}_{\gamma_{j+1}} \nabla^b n_{a_{j+1} \alpha} .
\end{align*}
The inductive hypothesis is used repeatedly in the above computation when moving normal vectors inside of covariant derivatives. Furthermore, note that we can always rearrange covariant derivatives because $(M,g)$ is flat. However, as a result of this, we know that the expression must be symmetric in $\big\{\gamma_1, \ldots, \gamma_{j+1}\big\}$. Now we may swap the Greek indices in the expression $n^{a_{j+1}}_{\gamma_{j+1}} \nabla^b n_{a_{j+1} \alpha}$ at the expense of a Gram matrix -- however, that Gram matrix is only differentiated $j+1$ times, which because $j+1 < m$, vanishes along $\Lambda$ by hypothesis. Thus, we~have
\begin{align*}
n^{a_1}_{\gamma_1} \cdots n^{a_{j+1}}_{\gamma_{j+1}} \nabla_{a_1} \cdots \nabla_{a_{j+1}} n^b_\alpha \= -n^{a_1}_{\gamma_1} \cdots n^{a_{j}}_{\gamma_{j}} \nabla_{a_1} \cdots \nabla_{a_j} n^{a_{j+1}}_{\alpha} \nabla^b n_{a_{j+1} \gamma_{j+1}} .
\end{align*}
Repeating this process of moving~$n$'s out, rearranging covariant derivatives, and pushing different~$n$'s in, we arrive at
\[
n^{a_1}_{\gamma_1} \cdots n^{a_{j+1}}_{\gamma_{j+1}} \nabla_{a_1} \cdots \nabla_{a_{j+1}} n^b_\alpha \= -n^{a_1}_{\gamma_1} \cdots n^{a_{j-1}}_{\gamma_{j-1}} n^{a_{j+1}}_{\alpha} \nabla_{a_1} \cdots \nabla_{a_{j-1}} \nabla_{a_{j+1}} n^{a_{j}}_{\gamma_{j}} \nabla^b n_{a_{j} \gamma_{j+1}} ,
\]
which vanishes because it is symmetric in $\gamma_j \gamma_{j+1}$ and hence can be replaced with $j+1$ derivatives of the Gram matrix.
\end{proof}

We are now equipped to produce a unique unit defining map.
\begin{Corollary}
Let $\Lambda^{d-k} \hookrightarrow \bigl(M^d,g\bigr)$ have vanishing normal curvature and let $(M,g)$ be flat. Then there exists a unique unit defining map $($up to constant rotations$)$.
\end{Corollary}
\begin{proof}
Observe from Lemma~\ref{rot-min-frame} that we may choose a parallelization such that $\beta = 0$, as the submanifold embedding has vanishing normal curvature. This parallelization is unique up to constant rotations. Given this parallelization, we may use Theorem~\ref{general-ho-ext} to determine the canonical defining map for this parallelization. Then it suffices to show that each obstruction~$F^{(i)}$ vanishes for this parallelization. Evidently, $F^{(2)}|_{\Lambda}$ vanishes because $\beta = 0$ and $(M,g)$ is flat. We would like to show that each subsequent obstruction vanishes along $\Lambda$. To do so, we will induct on $i$.

Suppose that $i \geq 2$ and that $G_{\alpha \beta} = \delta_{\alpha \beta} + \mathcal{O}\bigl(s^{i+1}\bigr)$. Then, from Theorem~\ref{general-ho-ext}, $F^{(i+1)}$ occupies a generalized window tensor representation of $S^k$, and thus it can be computed by projecting appropriately the expression \smash{$\tilde{n}^{a_1}_{(\gamma_1} \cdots \tilde{n}^{a_{i+1}}_{\gamma_{i+1})} \nabla_{a_1} \cdots \nabla_{a_{i+1}} \tilde{G}_{\alpha \beta}$}. By Lemma~\ref{lower-order-derivs} and the flatness condition on $(M,g)$, the only term that may appear in $F^{(i+1)}$ is of the form
\[
\tilde{n}^{a_1}_{(\gamma_1} \cdots \tilde{n}^{a_{i+1}}_{\gamma_{i+1})} n^b_{\alpha} \nabla_{a_1} \cdots \nabla_{a_{i+1}} \nabla_b s_{\beta} .
\]
But as this expression is symmetric in Latin indices, it is necessarily symmetric on all the Greek indices except $\beta$. Thus, the generalized window projection of this expression necessarily vanishes, as this expression is symmetric in $i+2$ Greek indices. Thus, $F^{(i+1)}|_{\Lambda} = 0$, and hence $G_{\alpha \beta} = \delta_{\alpha \beta} + \mathcal{O}\bigl(s^{i+2}\bigr)$, completing the induction.
\end{proof}

Notably, the unique unit defining map is that which produces a parallelization with vanishing~$\beta$ and is as described in Theorem~\ref{general-ho-ext}.

\section{Conformal geometry}\label{conf-geom-sec}
We now take a detour to explain the conformal calculus that will be used in the following sections. Much of what is summarized here can be found in more detail in any of~\cite{BEG,CapGover1999,CapGover2003,GOpet,GW}.

Consider a conformal manifold $\bigl(M^d,\cc\bigr)$ with $d \geq 3$. We may view this conformal manifold as a smooth ray subbundle $\mathcal{Q} \subset \odot^2 T^* M$ coordinatized by $\bigl(p,t^2 g\bigr)$ for $p \in M$, $t > 0$, and $g \in \cc$ a~metric on $M$. Indeed, this ray subbundle is a principal $\mathbb{R}_+$ bundle and hence, for each $w \in \mathbb{R}$, induces a natural associated line bundle, called a \textit{density bundle of weight $w$}, associated to $\mathcal{Q}$ via the irreducible representation $\mathbb{R}_+ \ni t \mapsto t^{-w} \in \operatorname{End}(\mathbb{R})$. Such a density bundle is denoted~$\ce M[w]$. Sections of this bundle can then be viewed as functions $\varphi$ on $\mathcal{Q}$ with the property that
\[\varphi\bigl(p, \Omega^2 g\bigr) = \Omega^w \varphi(p,g) .\]
We call such a function a (\textit{conformal}) \textit{density}. For practical purposes, it is often more convenient to view a density as a double equivalence class of functions, $\varphi = [g;f] = \bigl[\Omega^2 g; \Omega^w f\bigr]$, where $f$ is some function in $C^\infty M$. Note that we can tensor a density bundle of weight $w$ with any bundle~$\mathcal{V} M$ over $M$ to obtain a weighted bundle denoted by $\mathcal{V} M[w]$.

Of particular importance is the so-called \textit{conformal metric} $\bg$ which is the tautological section of $\odot^2 T^* M[2]$ given~by $\bg = [g;g]$. In the conformal setting, the conformal metric can be used for raising and lowering indices on tensor-valued densities.

A feature of conformal geometry is that, on the tangent bundle $TM$, there is no connection that respects the conformal structure. Indeed, even on densities, the naive exterior derivative only respects the conformal structure on weight-0 densities. On the other hand, for any $g \in \cc$, there exists a \textit{true scale} $0<\tau \in \Gamma(\ce M[1])$ which is a weight 1 positive density that determines~$g$ via~$g = \bg/\tau^2$. Then, for any true scale $\tau$, we can define a connection on weight-$w$ densities via $\nabla^\tau := \tau^w \circ d \circ \tau^{-w}$. A straightforward calculation shows that for any two true scales $\tau$ and $\tau' = \Omega \tau$ and any $\varphi \in \Gamma(\ce M[w])$, we have that
\[\bigl(\nabla^{\tau'} - \nabla^{\tau}\bigr) \varphi = -w \Upsilon \varphi ,\]
where $\Upsilon := d \log \Omega$. Nonetheless, this family of connections is insufficient for generating invariants of a conformal structure, so instead we turn to tractors. The \textit{standard tractor bundle} over $(M,\cc)$ is a rank $d+2$ vector bundle $\mathcal{T} M$, which for any choice of metric $g \in \cc$ is canonically isomorphic to
\[\mathcal{T}M \stackrel{g}{\cong} \ce M[1] \oplus TM[-1] \oplus \ce M[-1] .\]
We call this isomorphism a \textit{choice of splitting}. In a choice of splitting, we may thus decompose a section of the standard tractor bundle $T \in \Gamma(\ct M)$ according to
\[T^A \stackrel{g}{=} \bigl(\tau^+, \tau^a, \tau^-\bigr) .\]
Any choice of splitting can be related to another according to
\[T^A \stackrel{\Omega^2 g}{=} \biggl(\tau^+, \tau^a + \Upsilon^a \tau^+, \tau^- - \Upsilon_a \tau^a - \frac{1}{2} \Upsilon_a \Upsilon^a \tau^+\biggr) ,\]
where $\Upsilon^a := \bg^{ab} \Upsilon_a$. Observe that $\tau^+$ is independent of choice of splitting and hence is conformally invariant.

Note that for any $T,U \in \Gamma(\ct M)$, there exists an inner product that is independent of the choice of splitting given~by
(in a choice of splitting)
\[T \csdot U \stackrel{g}{=} \bg_{ab} \tau^a \mu^b + \tau^+ \mu^- + \tau^- \mu^+ ,\]
where $U^A \stackrel{g}{=} \bigl(\mu^+, \mu^a, \mu^-\bigr)$. Indeed, this inner product defines the \textit{tractor metric} $h \in \Gamma\bigl( \odot^2 \ct^* M\bigr)$ expressible in a choice of splitting as
\[h_{AB} \stackrel{g}{=} \begin{pmatrix}
0 & 0 & 1 \\
0 & \bg_{ab} & 0 \\
1 & 0 & 0
\end{pmatrix} .\]
Furthermore, this tractor metric provides an isomorphism between the standard tractor bundle and its dual
\[\ct^* M \stackrel{g}{\cong} \ce M[1] \oplus T^* M[1] \oplus \ce M[-1] .\]

The standard tractor bundle also appears in the construction of other useful bundles. By~tensoring the standard tractor bundle with the weight-$1$ density bundle, we obtain $\ct M[1]$, which has a canonically defined section called the \textit{canonical tractor}, given in any choice of splitting by $X^A \stackrel{g}{=} (0,0,1) \in \Gamma(\ct M[1])$. More generally, we shall use the terminology \textit{tractor bundle} to refer to any tensor product of copies of the standard tractor bundle and density bundles.

Notably, there exists a well-defined and canonical \textit{standard tractor connection} on $\mathcal{T}M$,
\begin{align*}
\nabla^\ct \colon \ \Gamma(\ct M) &\rightarrow \Gamma(T^* M \otimes \ct M) \\
T^B &\mapsto \nabla^\ct_a T^B \stackrel{g}{=} \begin{pmatrix} \nabla_a \tau^+ - \tau_a \\ \nabla_a \tau^b + \bg_a^b \tau^- + (P^g)_a^b \tau^+ \\ \nabla_a \tau^- - P^g_{ab} \tau^b \end{pmatrix} .
\end{align*}
It is this connection that ensures that the 1-form-valued tractor above is independent of the choice of splitting (and hence is a section of a tractor bundle). Note that this connection can be extended to tractor bundle. However, this connection does \textit{not} map to $\mathcal{T}^* M \times \mathcal{T}M$. Instead, we can use the tractor connection to build a tractor-valued operator known as the \textit{Thomas-$D$ operator}, given in a choice of splitting by
\begin{align*}
D_A \colon \ \Gamma\bigl(\ct^\Phi M[w]\bigr) &\rightarrow \Gamma\bigl(\ct^* M \otimes \ct^\Phi M[w-1]\bigr), \\
T^\Phi &\mapsto D_A T^\Phi \stackrel{g}{=} \begin{pmatrix}
(d+2w-2)w T^\Phi \\
(d+2w-2) \nabla^\ct_a T^\Phi \\
-(\Delta^\ct + w J^g) T^\Phi
\end{pmatrix} ,
\end{align*}
where $\Delta^\ct = \bg^{ab} \nabla^\ct_a \nabla^\ct_b$ is the tractor Laplacian. From the Thomas-$D$ operator, we can define another useful operator on $\ct^\Phi M[w]$ -- the \textit{hatted} Thomas-$D$ operator
\[\hd_A = \frac{1}{d+2w-2} D_A ,\]
which is well defined so long as $w \neq 1 - \frac{d}{2}$. While the Thomas-$D$ operator is not a derivation, it has the useful feature that $D_A \circ D^A = 0$. Furthermore, its failure to be a derivation is well controlled and given~by the following formula~\cite{Will1}:
\begin{gather*}
\hd^A (T_1 \cdots T_\ell) - \sum_{i=1}^\ell T_1 \cdots \bigl(\hd^A T_i\bigr) \cdots T_\ell \\
\qquad{}= - \frac{2 X^A}{d+2 \sum_{i=1}^\ell w_i -2} \sum_{1 \leq i < j \leq \ell} T_1 \cdots \bigl(\hd^B T_i\bigr) \cdots \bigl(\hd_B T_j\bigr) \cdots T_\ell ,
\end{gather*}
where for each $1 \leq i \leq \ell$, \smash{$T_i \in \Gamma\bigl(\ct^\Phi M[w_i]\bigr)$}, and the weight constraints $d+2\sum_{i=1}^{\ell} w_i - 2 \neq 0$ and \smash{$\big\{w_i \neq 1 - \frac{d}{2}\big\}_{i=1}^{\ell}$} are satisfied.
While it is clear that this operator is not a covariant derivative, it often plays the role of one. Indeed, the commutator of two hatted Thomas-$D$ operators produces a curvature-like tractor called the $W$-tractor, according to
\[\bigl[\hd_A, \hd_B\bigr] T^D = W_{AB}{}^D{}_E T^E + \frac{4}{d+2w-4} X_{[A} W_{B] C}{}^D{}_E \hd^C T^E ,\]
where such a commutator is defined. The $W$-tractor has Riemann-like symmetries, is trace-free, and satisfies
\[X^A W_{ABCD} = 0 = \hd^A W_{ABCD} .\]
This tractor is composed of the Weyl tensor, the Cotton tensor, and a generalized notion of the Bach tensor. See~\cite{BGW1, GOmin} for more details. Note that we will use the same letter for both the $W$-tractor and the Weyl tensor, but it will be clear from context (i.e., which indices are being used) which object is appearing in a given formula.

Because all of the aforementioned tractor objects respect the conformal structure of $(M,\cc)$, any combination of these tractors and tractor-valued operators can be used to produce a conformally invariant quantity. It is this feature that makes the tractor calculus analogous to the Ricci calculus of Riemannian geometry, and we will utilize this advantage heavily in the sections that follow.

\section[Orthonormal frames for normal bundles of conformal submanifold embeddings]{Orthonormal frames for normal bundles\\ of conformal submanifold embeddings} \label{conf-frames}
We now turn to the case of conformal submanifold embedings $\Lambda^{d-k} \hookrightarrow \bigl(M^d,\cc\bigr)$. As the tangent bundle decomposition for a Riemannian submanifold embedding in equation~(\ref{TM-decomp}) only depends on the orthogonality of vectors rather than their lengths, a conformal structure admits an identical decomposition. Furthermore, given any representative $g \in \cc$ and an orthonormal frame $\{(n^g)^a_{\alpha}\}$ with respect to $g$ of $N \Lambda$, we can promote such an orthonormal frame to an orthonormal frame-valued density (OFVD) $\{\bn^a_{\alpha}\} := \{[g ; (n^g)^a_{\alpha}]\}$ with weight $-1$. Then it follows that $\bg(\bn_{\alpha}, \bn_{\beta}) = \delta_{\alpha \beta}$.

Interestingly, normal fundamental forms are conformally-invariant. Indeed, given an OFVD and a metric representative $g \in \cc$, by direct computation of the conformal transformation of the Levi-Civita connection and the orthonormal frame-valued density, it follows that $\beta_{a \alpha \beta} \in \Gamma(T^* \Sigma[0] \times \raisebox{3pt}{\scalebox{0.2}{\ydiagram{1,1}}})$. So, the normal curvature is also invariant. Thus all of the results of Section~\ref{riem-frames} directly generalize to the equivalent conformal structures without additional consideration. Unfortunately, the results of Section~\ref{unit-def-funcs} do not generalize directly in the conformal setting. Instead, additional obstructions and other troublesome features arise.

\begin{Remark}
The conformal invariance of the normal fundamental form may seem to contradict the observation that, for a curve $\Lambda^1 \hookrightarrow \mathbb{R}^3$, the normal fundamental form of the Frenet frame picks out the torsion. However, this observation is naive. The normal fundamental form is conformally invariant when the orthonormal frame vectors are density valued. However, rescaling of the metric generically changes the plane of osculation for a given curve and hence rotates the Frenet frame. Thus the failure of a spacecurve's torsion to be conformally-invariant is a~consequence of the changing orientation of the Frenet frame.
\end{Remark}

\section[Defining densities for parallelized conformal submanifold embeddings]{Defining densities for parallelized conformal\\ submanifold embeddings} \label{conf-densities-sec}

As in the Riemannian setting, we would like to better study the conformal submanifold embedding $\Lambda^{d-k} \hookrightarrow \bigl(M^d,\cc\bigr)$ by studying the jets of geometrically-determined \textit{defining densities} -- densities of weight 1 that, in any metric representative, are represented by a defining map for $\Lambda$. To do so, we will attempt to construct the conformal analog of associated defining maps. We begin by picking any $g \in \cc$ and constructing the canonical defining map $s_{\alpha}^g$ for the provided OFVD in that metric representative. We can then promote this associated defining map to a defining density $\sigma_{\alpha}$ of weight 1. (Note that we are abusing language by calling an element of $\Gamma(\ce M[1] \times \raisebox{2pt}{\scalebox{0.2}{\ydiagram{1}}})$ a defining density.) Had we picked another metric conformally related to the first, the difference between the two densities would be at least second order in a defining density. For that reason, the initial choice of arbitrary metric to begin the construction is irrelevant for what follows (as we will change higher order terms anyway). That is, for a given OFVD, modulo terms of second order in a defining density, there exists a unique defining density $\sigma_{\alpha} = \bigl[g; s^g_{\alpha}\bigr] \in \Gamma(\ce M[1])$ such that \smash{$\mathrm{d} \sigma_{\alpha} \= \bg(\bn_{\alpha},\cdot)$}. (Note that $\mathrm{d} \sigma_{\alpha}$ can be defined given a choice of scale and depends on this choice. However, it is easy to see that its restriction to $\Lambda$ is independent of this choice, and hence can be regarded as a section of $T^*M[1]|_{\Lambda}$.) Our goal will be to maximally improve this density so that we can extract invariants of the parallelized conformal submanifold embedding.

In general, the Gram matrix of a defining map is not conformally invariant away from $\Lambda$. Even in the codimension 1 case, for two metric representatives $g, \Omega^2 g \in \cc$, we find that
\[\bigl|{\rm d} s^{\Omega^2 g}\bigr|^2_{\Omega^2 g} = |{\rm d} (\Omega s^g)|^2_{\Omega^2 g} = |\Omega (n^g_a + s \Upsilon_a)|^2_{\Omega^2 g} = |n^g_a + s^g \Upsilon_a|^2_g \neq |{\rm d}s^g|^2_g .\]
From this computation, we have shown explicitly that $G^{\Omega^2 g} \neq G^g$ except on the submanifold $\Lambda$. Going forward, we will be dropping the superscripts denoting metric representatives when the context is clear.

The natural analog of the (Riemannian) Gram matrix for a defining map to the conformal case is the \textit{conformal Gram matrix} of a defining density. Given a defining density for a~submanifold~$\sigma_\alpha$, we can build a family of \textit{scale tractors},
\[\big\{N^A_{\alpha} := \hd^A \sigma_{\alpha} \big\}_{\alpha = 1}^k .\]
The conformal Gram matrix is defined according to
\[\bG_{\alpha \beta} := h_{AB} N^A_{\alpha} N^B_{\beta} .\]
The conformal Gram matrix is conformally invariant by construction. Furthermore, upon restriction to the submanifold $\Lambda$, in any representative $g \in \cc$ the conformal Gram matrix agrees with the Riemannian Gram matrix of the representative of $\sigma_{\alpha}$ given~by $[g; s_{\alpha}]$. Indeed, for $g \in \cc$ any metric representative and $\sigma_{\alpha} = [g;s_{\alpha}]$ a defining density for $\Lambda^{d-k} \hookrightarrow (M,\cc)$, the components of $\bG_{\alpha \beta}$ are given~by
\[\bG_{\alpha \beta} = [g; G_{\alpha \beta}] := \biggl[g; (\nabla^a s_{\alpha})(\nabla_a s_{\beta}) - \frac{2}{d} s_{(\alpha} (\Delta^g + J^g) s_{\beta)}\biggr] .\]
Note that we are overloading notation: going forward we use $G_{\alpha \beta}$ to refer to a representative of the conformal Gram matrix, rather than the Riemannian Gram matrix. However note that every representative is the same, as the conformal Gram matrix has weight 0 -- we merely mean that we refer not to the double equivalence class $[g;G_{\alpha \beta}]$ but just $G_{\alpha \beta}$.

In the Riemannian setting, we found that, for a given parallelization of $N \Lambda$, a unique family of canonical defining maps could be obtained to arbitrary order. Further, when the background was flat and the normal curvature vanished, we found that this family consists entirely of unit defining maps. We proceed now with a similar program in the conformal setting.

\subsection{First order} \label{conf-first-order}
Given an initial choice of defining density built from an associated defining map for $\Lambda \hookrightarrow (M,\cc)$, it is easy to see that any representative of the conformal Gram matrix agrees with the Riemannian Gram matrix for the same associated defining map along $\Lambda$. It follows then that
\[\bG_{\alpha \beta} = \delta_{\alpha \beta} + F^{(1)}_{\alpha \beta \gamma} \sigma_{\gamma} \]
for some $F^{(1)} \in \Gamma(\ce M[-1] \times (\raisebox{3pt}{\scalebox{0.2}{\ydiagram{2,1}}} \; \oplus \; \raisebox{2pt}{\scalebox{0.2}{\ydiagram{3}}}))$, similar to the Riemannian setting. As in that setting, we will pick \textit{some} $F^{(1)}$ that extends $F^{(1)}|_{\Lambda}$, with the goal of showing that this choice is immaterial. Note, however, both the trace and trace-free parts of these representations will play a role in what follows.

We will now improve $\sigma$ at second order to cancel $F^{(1)}$ along $\Lambda$. Let
\[\tilde{\sigma}_{\alpha} = \sigma_{\alpha} + A^{(1)}_{\alpha \gamma_1 \gamma_2} \sigma_{\gamma_1} \sigma_{\gamma_2} ,\]
where $A^{(1)} \in \Gamma(\ce M[-1] \times (\raisebox{3pt}{\scalebox{0.2}{\ydiagram{2,1}}} \; \oplus \; \raisebox{2pt}{\scalebox{0.2}{\ydiagram{3}}}))$. Going forward, we will implicitly work in a fixed scale $g \in \cc$. To begin, we compute
\begin{align*}
&\tilde{n}_{\alpha} = n_{\alpha} + 2 A^{(1)}_{\alpha(\gamma_1 \gamma_2)} n_{\gamma_1} s_{\gamma_2} + \mathcal{O}\bigl(s^2\bigr), \\
& \Delta \tilde{s}_{\alpha} = \Delta s_{\alpha} + 2 A^{(1)}_{\alpha \gamma \gamma} + \mathcal{O}(s) .
\end{align*}
Thus, we have that
\begin{align*}
\tilde{G}_{\alpha \beta} &{}= \tilde{n}^a_{\alpha} \tilde{n}_{a \beta} - \frac{2}{d} \tilde{s}_{(\alpha} (\Delta^g + J^g) \tilde{s}_{\beta)} \\
&{}= n^a_{\alpha} n_{a \beta} + 4 n^a_{\gamma_1} n_{a (\alpha} A^{(1)}_{\beta)(\gamma_1 \gamma_2)} s_{\gamma_2} - \frac{2}{d} s_{(\alpha}(\Delta^g + J^g)s_{\beta)} - \frac{4}{d} s_{(\alpha} A^{(1)}_{\beta) \gamma \gamma} + \mathcal{O}\bigl(s^2\bigr) \\
&{}= G_{\alpha \beta} + 4 A^{(1)}_{(\alpha \beta) \gamma} s_{\gamma} - \frac{4}{d} s_{(\alpha} A^{(1)}_{\beta) \gamma \gamma} + \mathcal{O}\bigl(s^2\bigr)\\
&{}= \delta_{\alpha \beta} + \biggl(F^{(1)}_{\alpha \beta \omega} + 4 A^{(1)}_{(\alpha \beta) \omega} - \frac{4}{d} \delta_{\omega (\alpha} A^{(1)}_{\beta) \gamma \gamma}\biggr) s_{\omega} + \mathcal{O}\bigl(s^2\bigr) .
\end{align*}
Indeed, we see that a trace of $A^{(1)}$ appears here. Thus we must be careful to check that all irreducible components of $F^{(1)}$ can be cancelled. In this case, we have that
\[F^{(1)} = \bigl(\raisebox{3pt}{\scalebox{0.4}{\ydiagram{2,1}}}_{\,\circ} \oplus \raisebox{0pt}{\scalebox{0.4}{\ydiagram{1}}}\bigr) \; \oplus \; (\scalebox{0.4}{\ydiagram{3}}_{\,\circ} \oplus \raisebox{0pt}{\scalebox{0.4}{\ydiagram{1}}}) ,\]
and similarly for $A^{(1)}$. Note that the trace-free parts of $F^{(1)}$ can be cancelled for the same reason they can be cancelled in the Riemannian case. To force the traceful parts to vanish, we demand that
\begin{align*}
&0= F^{(1)}_{\alpha \alpha \omega} + 4 A^{(1)}_{\alpha \alpha \omega} - \frac{4}{d} A^{(1)}_{\omega \alpha \alpha} \;\;\text{ and } \\
&0= F^{(1)}_{\omega \alpha \alpha} + 2 A^{(1)}_{\omega \alpha \alpha} + 2 A^{(1)}_{\alpha \alpha \omega} - \frac{2}{d} A^{(1)}_{\omega \alpha \alpha} - \frac{2k}{d} A^{(1)}_{\omega \alpha \alpha} \\
& \hphantom{0}{}= F^{(1)}_{\omega \alpha \alpha} + 2 A^{(1)}_{\alpha \alpha \omega} + \frac{2(d-k-1)}{d} A^{(1)}_{\omega \alpha \alpha} .
\end{align*}
The determinant of this linear system \smash{\big(with $A^{(1)}_{\alpha \alpha \omega}$ and $A^{(1)}_{\omega \alpha \alpha}$ as variables\big)} is $8(d-k)/d$, and thus for all $k \neq d$, we can solve for the traces of $A^{(1)}$ such that they cancel the traces of $F^{(1)}$. We have thus shown the intermediate proposition:
\begin{Proposition}
Let $\Lambda^{d-k} \hookrightarrow \bigl(M^d,\cc\bigr)$ be a conformal submanifold embedding with a parallelization given in a choice of metric representative by $\{n_{\alpha}\}$. Then there exists a defining map $\sigma_{\alpha} \in \Gamma(\ce M[1] \times \raisebox{2pt}{\scalebox{0.2}{\ydiagram{1}}})$ unique up to order $\sigma^3$ such that $\mathrm{d} \sigma_{\alpha}|_{\Lambda} = [g; n_{a \alpha}]$ and its conformal Gram matrix satisfies
\[\bG_{\alpha \beta} = \delta_{\alpha \beta} + \mathcal{O}\bigl(\sigma^2\bigr) .\]
\end{Proposition}

Given such a defining map that satisfies $\bG_{\alpha \beta} = \delta_{\alpha \beta} + \mathcal{O}\bigl(\sigma^2\bigr)$, it will be useful to evaluate~$N^A_{\alpha}|_{\Lambda}$. Doing so requires the evaluation of $\nabla \csdot n_{\alpha}|_{\Lambda}$ in a choice of metric representative. First, note that, unlike in the Riemannian case, we have that
\[n^a_{\alpha} n_{a \beta} = \delta_{\alpha \beta} + \frac{2}{d} s_{(\alpha} \Delta s_{\beta)} + \mathcal{O}\bigl(s^2\bigr) .\]
Hence, following (with a slight modification) the computation leading to equation~(\ref{ndn-Riem}), we find that
\[n^b_{\alpha} \nabla_a n_{b \beta} \= \beta_{a \alpha \beta} + \frac{1}{d} n_{a \alpha} \nabla \csdot n_{\beta} .\]
Consequently, we have that
\begin{align*}
\nabla_a n_{b \beta} \= \II_{ab \beta} + n_{b \alpha} \beta_{a \alpha \beta} + n_{a \alpha} \beta_{b \alpha \beta} + \frac{1}{d}n_{a \alpha} n_{b \alpha} \nabla \csdot n_{\beta} .
\end{align*}
It thus follows that
\[\nabla \csdot n_{\beta} \=  d  H_{\beta} ,\]
where $H_{\beta} = \frac{1}{d-k} \II^a_{a\beta}$ is the mean curvature, and so
\[N^A_{\alpha} \stackrel{\Lambda,g}{=} \begin{pmatrix}
0\\ \, n^a_{\alpha}
\\
-H_{\alpha}
\end{pmatrix} .\]
We may also compute the tractor equivalent of the bulk formula for the normal fundamental form, given~by
\[B^A_{\alpha \beta} := N_{[\alpha} \csdot \hd N^A_{\beta]} \in \Gamma(\ct M[-1] \times \raisebox{2.5pt}{\scalebox{0.2}{\ydiagram{1,1}}} ) .\]
By direct computation, this tractor has the property that
\[X_A B^A_{\alpha \beta} = \hd_A B^A_{\alpha \beta} = 0 .\]
This tractor necessarily takes the form
\[B^A_{\alpha \beta} \stackrel{g}{=} \begin{pmatrix}
0\\
\, n^b_{[\alpha} \nabla_b n^a_{\beta]} - \frac{1}{d} n^a_{[\alpha} \nabla \csdot n_{\beta]} + \mathcal{O}(s)\\
\ast
\end{pmatrix}
\stackrel{\Lambda,g}{=}
\begin{pmatrix}
0\\
\, \beta^a_{\alpha \beta} \\
\ast
\end{pmatrix}
 .\]
Finally, it is useful to express the projecting part of the tractor equivalent of $\nabla_a n_b$,
\begin{align*}
P_{AB \alpha} :={}& \hd_A N_{B \alpha}\\
 ={}& \begin{pmatrix}
0 & 0 & 0 \\
0 & \nabla_a n_{b \alpha} + s_{\alpha}P_{ab} + g_{ab} \rho_{\alpha} & \ast \\
0 & \ast & \ast
\end{pmatrix}
 \= \begin{pmatrix}
0 & 0 & 0 \\
0 & \IIo_{ab \alpha} + 2 n_{(a \beta} \beta_{b) \beta \alpha} & \ast \\
0 & \ast & \ast
\end{pmatrix} ,
\end{align*}
where $\rho_{\alpha} := -\frac{1}{d}(\Delta s_{\alpha} + J s_{\alpha})$ and \smash{$\IIo_{\alpha} := \II_{\alpha} - \bar{g} H_{\alpha} \in \Gamma\bigl(\odot^2_{\circ} T^* \Lambda[1]\bigr)$} is the trace-free component of the second fundamental form.
With this tractor in hand, we define a scalar-valued density called the \textit{submanifold rigidity density} $K_{\alpha \beta} := P_{AB \alpha} P^{AB}_{\beta}$ which is ubiquitous. A short calculation shows that
\[K_{\alpha \beta} \= \IIo^2_{\alpha \beta} - 2 \beta_{b \gamma (\alpha} \beta^b_{\beta) \gamma} .\]
We now proceed to higher orders.

\subsection{Second order}

As in the Riemannian case, we now proceed to adjust $\sigma_{\alpha}$ at third order to eliminate as much from the second order term $\bG_{\alpha \beta}$ as possible. Pulling from the Riemannian case, we expect to fail here as well. Suppose that $\sigma_{\alpha}$ has
\[\bG_{\alpha \beta} = \delta_{\alpha \beta} + F^{(2)}_{\alpha \beta \gamma_1 \gamma_2} \sigma_{\gamma_1} \sigma_{\gamma_2} .\]
As a reminder, $F^{(2)}$ may depend on the initial choice made of $F^{(1)}$. Now, we can decompose~$F^{(2)}$ into irreducible representations so that
\[
F^{(2)} \in
\raisebox{2.5pt}{\scalebox{0.4}{\ydiagram{2,2}}}_{\,\circ}
 \; \oplus \;
 \raisebox{2.5pt}{\scalebox{0.4}{\ydiagram{3,1}}}_{\,\circ}
 \; \oplus \;
 \raisebox{0pt}{\scalebox{0.4}{\ydiagram{4}}}_{\,\circ}
\; \oplus \;
 3\, \raisebox{1pt}{\scalebox{0.4}{\ydiagram{2}}}_{\,\circ}
\; \oplus \;
 \raisebox{2.5pt}{\scalebox{0.4}{\ydiagram{1,1}}}
 \; \oplus \;
2.
\]
On the other hand, we can define $\tilde{\sigma}_{\alpha} = \sigma_{\alpha} + A^{(2)}_{\alpha \gamma_1 \gamma_2 \gamma_3 } \sigma_{\gamma_1} \sigma_{\gamma_2} \sigma_{\gamma_3}$. Then we can decompose $A^{(2)}$ in a similar fashion, yielding
\[A^{(2)} \in
 \raisebox{2.5pt}{\scalebox{0.4}{\ydiagram{3,1}}}_{\,\circ}
 \; \oplus \;
 \raisebox{0pt}{\scalebox{0.4}{\ydiagram{4}}}_{\,\circ}
\; \oplus \;
 2\, \raisebox{1pt}{\scalebox{0.4}{\ydiagram{2}}}_{\,\circ}
\; \oplus \;
 \raisebox{2.5pt}{\scalebox{0.4}{\ydiagram{1,1}}}
 \; \oplus \;
1.
\]
 So naively one might predict there will be obstructions of the form
\[
F^{(2)} \in
\raisebox{2.5pt}{\scalebox{0.4}{\ydiagram{2,2}}}_{\,\circ}
 \; \oplus \;
 \, \raisebox{1pt}{\scalebox{0.4}{\ydiagram{2}}}_{\,\circ}
 \; \oplus \;
1.
\]
To check this naive understanding, we will explicitly compute for the last time in this case. As before, we have that
\begin{align*}
&\tilde{n}_{\alpha}= n_{\alpha} + 3 A^{(2)}_{\alpha (\gamma_1 \gamma_2 \gamma_3)} s_{\gamma_1} s_{\gamma_2} n_{\gamma_3} + \mathcal{O}\bigl(s^3\bigr), \\
&\Delta \tilde{s}_{\alpha}= \Delta s_{\alpha} + 6 A^{(2)}_{\alpha \gamma_1 \gamma \gamma} s_{\gamma_1} + \mathcal{O}\bigl(s^2\bigr) .
\end{align*}
Plugging these into the formula for $G_{\alpha \beta}$, we have
\begin{align*}
\tilde{G}_{\alpha \beta} &{}= \tilde{n}^a_{\alpha} \tilde{n}_{a \beta} - \frac{2}{d} \tilde{s}_{(\alpha} (\Delta^g + J^g) \tilde{s}_{\beta)} \\
&{}= n^a_{\alpha} n_{a \beta} - \frac{2}{d} s_{(\alpha} (\Delta^g + J^g) s_{\beta)} + \biggl(6 A^{(2)}_{(\alpha \beta) \omega_1 \omega_2} - \frac{12}{d} \delta_{\omega_1 (\alpha} A^{(2)}_{\beta) \gamma \gamma \omega_2}\biggr) s_{\omega_1} s_{\omega_2} \\
&{}= \delta_{\alpha \beta} + \biggl(F^{(2)}_{\alpha \beta \omega_1 \omega_2} + 6 A^{(2)}_{(\alpha \beta) \omega_1 \omega_2} - \frac{12}{d} \delta_{\omega_1 (\alpha} A^{(2)}_{\beta) \gamma \gamma \omega_2}\biggr) s_{\omega_1} s_{\omega_2} .
\end{align*}
Note that there is implicit symmetrization over $\omega_1$ and $\omega_2$ in these displays. Just as in the Riemannian case, it is straightforward to see that the $ \raisebox{2.5pt}{\scalebox{0.2}{\ydiagram{3,1}}}_{\,\circ}$ and $\raisebox{1.5pt}{\scalebox{0.2}{\ydiagram{4}}}_{\,\circ} $ components of $F^{(2)}$ can be cancelled directly, and it is impossible to cancel the $\raisebox{2.5pt}{\scalebox{0.2}{\ydiagram{2,2}}}_{\,\circ}$ component of $F^{(2)}$.
We now consider the remaining components of $F^{(2)}$. If the trace-free parts of this coefficient are to vanish, we must demand that
\begin{align*}
&0= F^{(2)}_{\rho \rho (\omega_1 \omega_2)_\circ} + 6 A^{(2)}_{\rho \rho (\omega_1 \omega_2)_\circ} - \frac{12}{d} A^{(2)}_{(\omega_1 \omega_2)_\circ \rho \rho}, \\
&0
= F^{(2)}_{(\omega_1 \omega_2)_\circ \rho \rho} + \frac{6(d-2)}{d} A^{(2)}_{(\omega_1 \omega_2)_\circ \rho \rho}, \\
&0
= F^{(2)}_{\rho (\omega_1 \omega_2)_\circ \rho} + 3 A^{(2)}_{\rho \rho (\omega_1 \omega_2)_\circ} + \frac{3(d-k-2)}{d} A^{(2)}_{(\omega_1 \omega_2)_\circ \rho \rho}, \\
&0
= F^{(2)}_{\rho [\omega_1 \omega_2] \rho} + \frac{3(d-k-2)}{d} A^{(2)}_{[\omega_1 \omega_2] \rho \rho}, \\
&0
= F^{(2)}_{\rho \rho \gamma \gamma} + \frac{6(d-2)}{d} A^{(2)}_{\rho \rho \gamma \gamma}, \\
&0
= F^{(2)}_{\rho \gamma \gamma \rho} + \frac{6(d-k-1)}{d} A^{(2)}_{\rho \rho \gamma \gamma} .
\end{align*}
To more easily visualize this system of equations, we denote \smash{$x_1 = A^{(2)}_{\rho \rho (\omega_1 \omega_2)_\circ}$, $x_2 = A^{(2)}_{(\omega_1 \omega_2)_\circ \rho \rho}$}, \smash{$x_3 = A^{(2)}_{[\omega_1 \omega_2] \rho \rho}$}, and \smash{$x_4 = A^{(2)}_{\gamma \gamma \rho \rho}$}. Then the matrix describing the linear system of equations is given~by
\[
\begin{pmatrix}
6 & -\frac{12}{d} & 0 & 0 \\[1mm]
0 & \frac{6(d-2)}{d} & 0 & 0 \\[1mm]
3 & \frac{3(d-k-2)}{d} & 0 & 0 \\[1mm]
0 & 0 & \frac{3(d-k-2)}{d} & 0 \\[1mm]
0 & 0 & 0 & \frac{6(d-2)}{d} \\[1mm]
0 & 0 & 0 & \frac{6(d-k-1)}{d}
\end{pmatrix} .
\]
This matrix has rank at most 4. Furthermore, when $k = d-2$, this matrix has rank 3. \big(Note that the rank also drops when $k=1$, but all components of $F^{(2)}$ collapse to a single function, which is uniquely determined by fixing the single function $A^{(2)}$.\big) So as expected, we can remove at most 4 traced degrees of freedom from \smash{$F^{(2)}$} using \smash{$A^{(2)}$}, except when $k = d-2$, in which case we can only remove 3 degrees of freedom. Indeed, for $k \neq d-2$, we have an obstruction given~by
\[
F^{(2)} =
\raisebox{2.5pt}{\scalebox{0.4}{\ydiagram{2,2}}}_{\,\circ}
 \; \oplus \;
 \, \raisebox{1pt}{\scalebox{0.4}{\ydiagram{2}}}_{\,\circ}
 \; \oplus \;
1.
\]
When $k = d-2$, we instead have an obstruction of the form
\[
F^{(2)} =
\raisebox{2.5pt}{\scalebox{0.4}{\ydiagram{2,2}}}_{\,\circ}
 \; \oplus \;
 \, \raisebox{1pt}{\scalebox{0.4}{\ydiagram{2}}}_{\,\circ}
 \; \oplus \;
 \, \raisebox{1pt}{\scalebox{0.4}{\ydiagram{1,1}}}_{\,\circ}
 \; \oplus \;
1.
\]
For the sake of uniformity, when $k \neq d-2$, we will choose to eliminate as much of the components of $F^{(2)}$ in $\raisebox{2.5pt}{\scalebox{0.2}{\ydiagram{3,1}}}$ and $\raisebox{1.5pt}{\scalebox{0.2}{\ydiagram{4}}}$ as possible. Certainly, as in the Riemannian case, we can easily remove the trace-free parts of these components, so we examine the traces.

First, observe that
\[\delta^{\alpha \beta} \delta^{\gamma \delta} F^{(2)}_{(\alpha \beta \gamma \delta)} = \frac{1}{3} \bigl(F^{(2)}_{\alpha \alpha \beta \beta} + 2 F^{(2)}_{\alpha \beta \beta \alpha}\bigr) .\]
So we demand that
\[F^{(2)}_{\alpha \alpha \beta \beta} + 2 F^{(2)}_{\alpha \beta \beta \alpha} + \frac{6(3d-4-2k)}{d} A_{\alpha \alpha \beta \beta}^{(2)} = 0 .\]
This has solutions so long as $(d,k) \neq (2,1)$, and thus we can always remove the full trace of \smash{$F^{(2)}_{(\alpha \beta \gamma \delta)}$} in that setting. In fact, when $(d,k) = (2,1)$, we may directly compute this full trace term (which is in fact the only possibly obstruction, as $k=1$) and see that it vanishes. However, as $A^{(2)}$ is not determined by the total trace of $F^{(2)}$, we may not proceed there.

Next we consider the trace-free symmetric part of the trace of the totally-symmetric component. By the same considerations as before, this component is given~by
\[\frac{1}{6}\bigl(F_{\gamma \gamma (\alpha \beta)_\circ} + 4 F_{\gamma (\alpha \beta)_\circ \gamma} + F_{(\alpha \beta)_\circ \gamma \gamma}\bigr) .\]
Thus, we demand that
\begin{align} \label{tfs-fully-symmetric}
F_{\gamma \gamma (\alpha \beta)_\circ} + 4 F_{\gamma (\alpha \beta)_\circ \gamma} + F_{(\alpha \beta)_\circ \gamma \gamma} + 18 A^{(2)}_{\gamma \gamma (\alpha \beta)_\circ} + \frac{6(3d-2k-8)}{d} A^{(2)}_{(\alpha \beta)_\circ \gamma \gamma} = 0 .
\end{align}
To fully determine both of the components of $A^{(2)}$ above, we can also demand that the trace-free symmetric part of the trace of the $\raisebox{2.5pt}{\scalebox{0.2}{\ydiagram{3,1}}}$ component of $F^{(2)}$ vanishes. This component of $F^{(2)}$ can be expressed as
\[\frac{1}{6}\bigr(F^{(2)}_{(\alpha \beta)_\circ \gamma \gamma} - F^{(2)}_{\gamma \gamma (\alpha \beta)_\circ}\bigl) .\]
So, we must have that
\begin{align} \label{tfs-pistol}
F^{(2)}_{(\alpha \beta)_\circ \gamma \gamma} - F^{(2)}_{\gamma \gamma (\alpha \beta)_\circ} - 6 A^{(2)}_{\gamma \gamma (\alpha \beta)_\circ} + 6 A^{(2)}_{(\alpha \beta)_\circ \gamma \gamma} =0 .
\end{align}
Combining equations~(\ref{tfs-fully-symmetric}) and~(\ref{tfs-pistol}), we find that the determinant of the system is $\frac{72}{d} (3d-k-4)$, which never vanishes for $1 \leq k \leq d-1$. Thus, we can completely remove the trace-free symmetric pieces of both $\raisebox{2.5pt}{\scalebox{0.2}{\ydiagram{3,1}}}$ and $\raisebox{1.5pt}{\scalebox{0.2}{\ydiagram{4}}}$.

Finally, we wish to remove the antisymmetric part of the trace of the $\raisebox{2.5pt}{\scalebox{0.2}{\ydiagram{3,1}}}$ component of $F^{(2)}$, which is given~by \smash{$F^{(2)}_{\gamma [\alpha \beta] \gamma}$}. So long as $k \neq d-2$, this piece can be removed as well. This accounts for all of the degrees of freedom in $A^{(2)}$, and thus we are left with a completely undetermined $\raisebox{2.5pt}{\scalebox{0.2}{\ydiagram{2,2}}}$ component of $F^{(2)}$, precisely as in the Riemannian setting.

Having fixed $A^{(2)}$ and in so doing, forced $F^{(2)} \in \Gamma(\ce M[-2] \times \raisebox{2.5pt}{\scalebox{0.2}{\ydiagram{2,2}}})$, it is now instructive to compute this obstruction.
We begin by noting that
\begin{align} \label{window-projection-bG}
\bG_{\alpha \beta} = \delta_{\alpha \beta} + \bigl(P_{\scalebox{0.2}{\ydiagram{2,2}}}\, F^{(2)}_{\alpha \beta \gamma_1 \gamma_2} \bigr) \sigma_{\gamma_1} \sigma_{\gamma_2} .
\end{align}
Furthermore, we have that
\[ P_{\scalebox{0.2}{\ydiagram{2,2}}}\, F^{(2)}_{\alpha \beta \gamma_1 \gamma_2} = \frac{1}{3} \bigl(2 F^{(2)}_{\alpha \beta \gamma_1 \gamma_2} + 2 F^{(2)}_{\gamma_1 \gamma_2 \alpha \beta} - 2 F^{(2)}_{\gamma_1 (\alpha \beta) \gamma_2} - 2F^{(2)}_{\gamma_2 (\alpha \beta) \gamma_1} \bigr) .\]
Now consider the following projection:
\[P_{\scalebox{0.2}{\ydiagram{2,2}}}\, N^A_{\gamma_1} N^B_{\gamma_2} \hd_A \hd_B \bG_{\alpha \beta} .\]
Using equation~(\ref{window-projection-bG}) and simplifying, we find that
\begin{align}
P_{\scalebox{0.2}{\ydiagram{2,2}}}\, N^A_{\gamma_1} N^B_{\gamma_2} \hd_A \hd_B \bG_{\alpha \beta} &{}\=
2F^{(2)}_{\alpha \beta \gamma_1 \gamma_2}
+ \frac{2}{9(d-2)} \delta_{\alpha \beta} \bigl(F^{(2)}_{\rho \gamma_1 \gamma_2 \rho} + F^{(2)}_{\rho \gamma_2 \gamma_1 \rho} - F^{(2)}_{\gamma_1 \gamma_2 \rho \rho} - F^{(2)}_{\rho \rho \gamma_1 \gamma_2}\bigr) \nonumber\\
&\quad{}+ \frac{2}{9(d-2)} \delta_{\gamma_1 \gamma_2} \bigl(F^{(2)}_{\rho \alpha \beta \rho} + F^{(2)}_{\rho \beta \alpha \rho} - F^{(2)}_{\alpha \beta \rho \rho} - F^{(2)}_{\rho \rho \alpha \beta}\bigr) \nonumber\\
&\quad{}+ \frac{1}{9(d-2)} \delta_{\gamma_1 \alpha} \bigl(F^{(2)}_{\gamma_2 \beta \rho \rho} + F^{(2)}_{\rho \rho \gamma_2 \beta} - F^{(2)}_{\rho \gamma_2 \beta \rho} - F^{(2)}_{\rho \beta \gamma_2 \rho}\bigr) \nonumber\\
&\quad{}+ \frac{1}{9(d-2)} \delta_{\gamma_2 \alpha} \bigl(F^{(2)}_{\gamma_1 \beta \rho \rho} + F^{(2)}_{\rho \rho \gamma_1 \beta} - F^{(2)}_{\rho \gamma_1 \beta \rho} - F^{(2)}_{\rho \beta \gamma_1 \rho}\bigr) \nonumber\\
&\quad{}+ \frac{1}{9(d-2)} \delta_{\gamma_1 \beta} \bigl(F^{(2)}_{\gamma_2 \alpha \rho \rho} + F^{(2)}_{\rho \rho \gamma_2 \alpha} - F^{(2)}_{\rho \gamma_2 \alpha \rho} - F^{(2)}_{\rho \alpha \gamma_2 \rho}\bigr) \nonumber\\
&\quad{}+ \frac{1}{9(d-2)} \delta_{\gamma_2 \beta} \bigl(F^{(2)}_{\gamma_1 \alpha \rho \rho} + F^{(2)}_{\rho \rho \gamma_1 \alpha} - F^{(2)}_{\rho \gamma_1 \alpha \rho} - F^{(2)}_{\rho \alpha \gamma_1 \rho}\bigr) . \label{conf-obs}
\end{align}

We may also compute the same projection using the fact that $\bG_{\alpha \beta} = N_{\alpha} \cdot N_{\beta}$. This computation is somewhat involved, but results in
\begin{align}
P_{\scalebox{0.2}{\ydiagram{2,2}}}\, N^A_{\gamma_1} N^B_{\gamma_2} \hd_A \hd_B \bG_{\alpha \beta} &{}\= -2\beta_{a \gamma_1 (\alpha} \beta^a_{\beta) \gamma_2} - \frac{2}{3} W_{\gamma_1 (\alpha \beta) \gamma_2} \nonumber\\
&\quad{}+ \frac{2}{3(d-2)} \bigl[
\delta_{\gamma_1 \gamma_2} \bigl(- \IIo^2_{\alpha \beta} + 2 \beta_{a \alpha \rho} \beta^a_{\rho \beta}\bigr)\! + \delta_{\alpha \beta} \bigl(- \IIo^2_{\gamma_1 \gamma_2} + 2 \beta_{a \gamma_1 \rho} \beta^a_{\rho \gamma_2}\bigr) \nonumber\\
&\quad{}
+ \delta_{\gamma_1 (\alpha} \bigl(\IIo^2_{\beta) \gamma_2 } - 2 \beta^a_{\beta) \rho } \beta_{a \rho \gamma_2} \bigr)
+ \delta_{\gamma_2 (\alpha} \bigl(\IIo^2_{\beta) \gamma_1 } - 2\beta^a_{ \beta) \rho} \beta_{a \rho \gamma_1} \bigr) \bigr] . \!\!\!\label{conf-obs-res}
\end{align}

Now, as traces are combined in non-trivial ways with the above expressions, it is fruitful to compute the irreducible components of $F^{(2)}|_{\Lambda}$ independently. To do so, first observe that
\[Tr_{\scalebox{0.2}{\ydiagram{2,2}}}^2 F^{(2)} := \delta^{\alpha \beta} \delta^{\gamma \delta} P_{\scalebox{0.2}{\ydiagram{2,2}}}\, F^{(2)}_{\alpha \beta \gamma \delta} = \frac{2}{3}\bigl(F^{(2)}_{\alpha \alpha \beta \beta} - F^{(2)}_{\alpha \beta \beta \alpha}\bigr) .\]
Evaluating this trace combination on equation~(\ref{conf-obs}), we have that
\[Tr_{\scalebox{0.2}{\ydiagram{2,2}}}^2 N^A_{\gamma_1} N^B_{\gamma_2} \hd_A \hd_B \bG_{\alpha \beta} = \frac{4(3d-2k-4)}{9(d-2)}\bigl(F^{(2)}_{\alpha \alpha \beta \beta} - F^{(2)}_{\alpha \beta \beta \alpha}\bigr) = \frac{2(3d-2k-4)}{3(d-2)} Tr_{\scalebox{0.2}{\ydiagram{2,2}}}^2 F^{(2)} .\]
But we can also take the double trace of equation~(\ref{conf-obs-res}). Doing so yields
\[Tr_{\scalebox{0.2}{\ydiagram{2,2}}}^2\, N^A_{\gamma_1} N^B_{\gamma_2} \hd_A \hd_B \bG_{\alpha \beta} \= -\frac{4(k-1)}{3(d-2)} \IIo^2_{\alpha \alpha} + \frac{2(3d-4k-2)}{3(d-2)} \beta_{a \alpha \beta} \beta^a_{\alpha \beta} - \frac{2}{3} W_{\alpha \beta \beta \alpha} .\]
Thus, we have that
\begin{align*}
Tr_{\scalebox{0.2}{\ydiagram{2,2}}}^2 F^{(2)} \= - \frac{2(k-1)}{3d-2k-4} \IIo^2_{\alpha \alpha} + \frac{3d-4k-2}{3d-2k-4} \beta_{a \alpha \beta} \beta^a_{\alpha \beta} - \frac{d-2}{3d-2k-4} W_{\alpha \beta \beta \alpha} .
\end{align*}

In general, this obstruction \smash{$\mathrm{Tr}_{\scalebox{0.2}{\ydiagram{2,2}}}^2 F^{(2)}$} depends in a non-trivial way on the parallelization of $N \Lambda$ through the appearance of $\beta$, however for special choices of $(d,k)$, we find that this obstruction is a true submanifold invariant that does not depend on the choice of orthonormal frame. This follows because $\IIo_{\alpha}$ and $n_{\alpha}$ both occupy tensor representations of ${\rm O}(k)$, unlike $\beta$. This occurs when $3d-4k-2 = 0$, i.e., for any integer $n \in \mathbb{N}$, then $\mathrm{Tr}_{\scalebox{0.2}{\ydiagram{2,2}}}^2 F^{(2)}$ is a true invariant so long as $(d,k) = (4n+2,3n+1)$.

We next consider the tracefree part of the trace of $F^{(2)}$. As above, it is useful to compute which component this is in terms of $F^{(2)}$ prior to projection. There we have that
\[\bigl[\bigl( \operatorname{t.f.} \circ \mathrm{Tr}_{\scalebox{0.2}{\ydiagram{2,2}}} \bigr) F^{(2)}\bigr]_{\alpha \beta} := \frac{1}{3} F^{(2)}_{(\alpha \beta)_\circ \rho \rho} + \frac{1}{3} F^{(2)}_{\rho \rho (\alpha \beta)_\circ} - \frac{2}{3} F^{(2)}_{\rho (\alpha \beta)_\circ \rho} ,\]
where $\operatorname{t.f.}$ is a projection operator to the trace-free part of a tensor.
Taking this trace of equations~(\ref{conf-obs}) and~(\ref{conf-obs-res}), we find that
\begin{gather*}
\bigl[\bigl( \operatorname{t.f.} \circ \mathrm{Tr}_{\scalebox{0.2}{\ydiagram{2,2}}} \bigr) F^{(2)}\bigr]_{\alpha \beta} \= - \frac{k-2}{3(3d-k-4)} \IIo^2_{(\alpha \beta)_\circ} - \frac{3d-2k-2}{3(3d-k-4)} \beta_{a \rho (\alpha} \beta^a_{\beta)_\circ \rho} \\ \hphantom{\bigl[\bigl( \operatorname{t.f.} \circ Tr_{\scalebox{0.2}{\ydiagram{2,2}}} \bigr) F^{(2)}\bigr]_{\alpha \beta} \=}{}
 - \frac{d-2}{3(3d-k-4)} W_{\rho (\alpha \beta)_{\circ} \rho} .
\end{gather*}
Unlike the double trace obstruction, there are no pairs $(d,k)$ such that this is a true extrinsic submanifold invariant. The final component of $F^{(2)}$, the trace-free part, can be read off from equations~(\ref{conf-obs}) and~(\ref{conf-obs-res}).

Now, given that $F^{(2)}$ occupies the window tensor representation, we have that $F^{(2)}_{\alpha \beta \gamma \delta} = F^{(2)}_{\gamma \delta \alpha \beta}$ and that \smash{$F^{(2)}_{(\alpha \beta \gamma) \delta} = 0$}. It follows that \smash{$F^{(2)}_{\alpha \beta \gamma \delta} + 2 F^{(2)}_{\gamma (\alpha \beta) \delta} = 0$}, and so
\[Tr_{\scalebox{0.2}{\ydiagram{2,2}}}^2 F^{(2)} = \frac{1}{3}F^{(2)}_{\alpha \alpha \beta \beta} \qquad \text{and} \qquad \bigl[\bigl( \operatorname{t.f.} \circ \mathrm{Tr}_{\scalebox{0.2}{\ydiagram{2,2}}} \bigr) F^{(2)}\bigr]_{\alpha \beta} = \frac{2}{3} F^{(2)}_{\rho \rho (\alpha \beta)_{\circ}} .\]
We now summarize all of these results in the following theorem.
\begin{Theorem} \label{ho-conformal-construction}
Let $\Lambda^{d-k} \hookrightarrow \bigl(M^d,\cc\bigr)$ be a parallelized submanifold embedding with $k \neq d-2$. Then there exists a family of defining densities $\sigma_{\alpha}$ for $\Lambda$ that agree modulo terms of order $\sigma^4$ such that
\[\bG_{\alpha \beta} = \delta_{\alpha \beta} + F^{(2)}_{\alpha \beta \gamma_1 \gamma_2} \sigma_{\gamma_1} \sigma_{\gamma_2} ,\]
where $F^{(2)} \in C^\infty M \times \raisebox{3pt}{\scalebox{0.2}{\ydiagram{2,2}}}$, and
\[\operatorname{t.f.} F^{(2)}_{\alpha \beta \gamma_1 \gamma_2} \= \operatorname{t.f.} \biggl(- \beta_{a \gamma_1 (\alpha} \beta^a_{\beta) \gamma_2} - \frac{1}{3} W_{\gamma_1 (\alpha \beta) \gamma_2} \biggr) ,\]
\[F^{(2)}_{\rho \rho (\alpha \beta)_\circ} \= - \frac{k-2}{2(3d-k-4)} \IIo^2_{(\alpha \beta)_\circ} - \frac{3d-2k-2}{2(3d-k-4)} \beta_{a \rho (\alpha} \beta^a_{\beta)_\circ \rho} - \frac{d-2}{2(3d-k-4)} W_{\rho (\alpha \beta)_{\circ} \rho} ,\]
and
\[F^{(2)}_{\alpha \alpha \beta \beta}\= - \frac{6(k-1)}{3d-2k-4} \IIo^2_{\alpha \alpha} + \frac{3(3d-4k-2)}{3d-2k-4} \beta_{a \alpha \beta} \beta^a_{\alpha \beta} - \frac{3(d-2)}{3d-2k-4} W_{\alpha \beta \beta \alpha} .\]
\end{Theorem}
\noindent
We call these defining densities \textit{associated defining densities} for $\Lambda^{d-k} \hookrightarrow \bigl(M^d,g\bigr)$.

We now consider the case where $k = d-2$. As observed above, we may not set \smash{$F^{(2)}_{\gamma [\alpha \beta] \gamma}$} to zero, and hence this term obstructs setting the Gram matrix to the identity. Indeed, we can compute this obstruction and find that
\[F^{(2)}_{\gamma [\alpha \beta] \gamma} \= \bar{\nabla}^a \beta_{a \alpha \beta} .\]
As noted in Section~\ref{gauge-fixed}, when $\Lambda$ is closed and $N \Lambda$ admits a~parallelization, there exists such a~parallelization such that this obstruction vanishes. Further, this is only conformally invariant for $k=d-2$: the divergence of a weight-0 1-form is only conformally invariant on surfaces. Regardless, we find that $\sigma_{\alpha}$ cannot be uniquely fixed to third order, as there is no constraint equation controlling \smash{$A^{(2)}_{[\alpha \beta] \rho \rho}$}. To that end, when $k = d-2$, we define an equivalence class of defining densities, given~by the following relation:
\[[\sigma_{\alpha}] \sim [\sigma_{\alpha} + A_{[\alpha \gamma]} \sigma_{\gamma} \sigma_{\rho} \sigma_{\rho}] \]
for any $A_{[\alpha \beta]} \in \Gamma(\ce M[-2] \times \raisebox{2.5pt}{\scalebox{0.2}{\ydiagram{1,1}}})$.

Note that despite the trouble that this sporadic obstruction provides, the remainder of the obstruction still occupies a window tensor representation and the projections to its components are unchanged from the $k \neq d-2$ case.

\begin{Remark}
As in Remark~\ref{curve-invariants}, because the orthonormal frame used to parallelize the normal bundle is geometrically determined for curves embedded in a conformal manifold, the obstructions at second order produced above are true submanifold invariants. The Bishop frame used here is of course different from the frame constructed in Fialkow's classical work~\cite{Fialkow2} on conformal curves, and would give rise to a different (but related) family of obstructions.
\end{Remark}

\subsection{Third order}
While there is an obstruction at second order, we can still ask which components at third order can be removed by a judicious choice of $\sigma_{\alpha}$. As the case where $k=d-2$ is qualitatively different from that of $k\neq d-2$, we consider these situations separately, beginning with the case where $k \neq d-2$. We will follow the method applied in the Riemannian case. To that end, our goal will be to show that there \textit{exists} such a defining density and then establish that this defining density is independent (to sufficiently higher order) on the initial choice of $F^{(1)}$.

We can now proceed to extend $\sigma_{\alpha}$ as in the Riemannian setting,
\[\sigma_{\alpha} \mapsto \tilde{\sigma}_{\alpha} := \sigma_{\alpha} + A^{(3)}_{\alpha \gamma_1 \gamma_2 \gamma_3 \gamma_4} \sigma_{\gamma_1} \sigma_{\gamma_2} \sigma_{\gamma_3} \sigma_{\gamma_4} ,\]
so that the Gram matrix takes the form
\begin{align}
\tilde{\bG}_{\alpha \beta} &{}= \delta_{\alpha \beta}+ F^{(2)}_{\alpha \beta \gamma_1 \gamma_2} \sigma_{\gamma_1} \sigma_{\gamma_2}\nonumber \\
&\quad{}+ \biggl[F^{(3)}_{\alpha \beta \gamma_1 \gamma_2 \gamma_3} + 8 A^{(3)}_{(\alpha \beta) \gamma_1 \gamma_2 \gamma_3} - \frac{24}{d} \delta_{\gamma_1 (\alpha} A^{(3)}_{\beta) \gamma \gamma \gamma_2 \gamma_3} \biggr] \sigma_{\gamma_1} \sigma_{\gamma_2} \sigma_{\gamma_3} .\label{3-order}
\end{align}
Now we could, in principle, perform the same representation theory analysis as in the second order case. Surely, we may force the trace-free part of $\tilde{F}^{(3)}$ to occupy a generalized window symmetry, as that calculation amounts to a rehashing of the Riemannian case. However, computations involving the representation theory of the remaining components quickly get out of hand. Instead, we are only interested when there are obstructions, special or otherwise. Because~$\tilde{F}^{(3)}$ has five Greek indices, there are no total trace terms, so the terms we will examine are terms that have one Greek index, i.e., are the result of a double trace.

First, observe that when $k=1$, there is only one maximal trace term, \smash{$F^{(3)}_{11111}$}, and only one term that contributes to the correction, \smash{$A^{(3)}_{11111}$}. As we are considering the case where $k \neq d-2$, this means that $d \neq 3$. So in this setting, it is easy to see that we may remove \smash{$F^{(3)}_{11111}$} in its entirety.

When $k\neq 1$, it is useful to consider the three distinct traces we may obtain at third order that can appear in equation~(\ref{3-order}). Tracing, we obtain
\begin{align*}
&\tilde{F}^{(3)}_{\gamma \gamma \rho \rho \alpha} = F^{(3)}_{\gamma \gamma \rho \rho \alpha }- \frac{8}{d} A^{(3)}_{\alpha \gamma \gamma \rho \rho} + \frac{8(d-2)}{d} A^{(3)}_{\gamma \gamma \rho \rho \alpha} ,\\
&\tilde{F}^{(3)}_{\gamma \rho \rho \gamma \alpha} = F^{(3)}_{\gamma \rho \rho \gamma \alpha }+ \frac{8(d-k-2)}{d} A^{(3)}_{\gamma \gamma \rho \rho \alpha} ,\\
&\tilde{F}^{(3)}_{\alpha \gamma \gamma \rho \rho} = F^{(3)}_{\alpha \gamma \gamma \rho \rho}+ \frac{4(d-k-3)}{d} A^{(3)}_{\alpha \gamma \gamma \rho \rho} + \frac{4(d-2)}{d} A^{(3)}_{\gamma \gamma \rho \rho \alpha} .
\end{align*}

The vanishing of these expressions clearly leads to an overdetermined system, so we must choose at least one trace that is unresolved by even the most judicious choice of $A^{(3)}$. In this case, we choose \smash{$F^{(3)}_{\gamma \rho \rho \gamma \alpha}$} to always be unfixed for reasons that will become clear later. In general, we can solve the remaining system. Suppose we let $\sigma_{\alpha} \mapsto \tilde{\sigma}_{\alpha}$ such that \smash{$\tilde{F}^{(3)}_{\gamma \gamma \rho \rho \alpha}$} vanishes (this is always possible) -- we do this by performing an arbitrary total-trace transformation of $\sigma_{\alpha}$ and then fixing \smash{$A^{(3)}_{\alpha \gamma \gamma \rho \rho}$} to depend on \smash{$F^{(3)}_{\gamma \gamma \rho \rho \alpha}$} and an arbitrary choice of \smash{$A^{(3)}_{\gamma \gamma \rho \rho \alpha}$} such that
\[{F}^{(3)}_{\gamma \gamma \rho \rho \alpha} \mapsto \tilde{F}^{(3)}_{\gamma \gamma \rho \rho \alpha} = 0 .\]
Under the same modification of $\sigma_{\alpha}$, we have that
\[F^{(3)}_{\alpha \gamma \gamma \rho \rho} \mapsto \tilde{F}^{(3)}_{\alpha \gamma \gamma \rho \rho} = F^{(3)}_{\alpha \gamma \gamma \rho \rho} + \frac{d-k-3}{2} F^{(3)}_{\gamma \gamma \rho \rho \alpha} + \frac{4(d-2)(d-k-2)}{d} A^{(3)}_{\gamma \gamma \rho \rho \alpha} .\]
Thus, for $k\neq d-2$, we can fix $A^{(3)}_{\gamma \gamma \rho \rho \alpha}$ so that $\tilde{F}^{(3)}_{\alpha \gamma \gamma \rho \rho} = 0$. It is entirely expected, however, that we would not be able to remove all of the maximal traces -- we expect that at least an entire generalized window's worth of obstructions will appear at third order, and that includes precisely one maximal trace term.

We now examine whether the construction provided so far is independent of the choice of~$F^{(1)}$ by mimicking the argument in Section~\ref{ho-Riem}. To do so, suppose that we had chosen two distinct extensions of $F^{(1)}|_{\Lambda}$, call them $F^{(1)}$ and $F^{(1)'}$. When $k\neq d-2$, these choices lead to two (naively) distinct defining densities $\sigma_{\alpha}$ and $\sigma_{\alpha}'$ that are determined modulo terms of fourth order. However, as their differentials must agree in any choice of metric representative along $\Lambda$, we must have that they differ at second order,
\[\sigma_{\alpha}' = \sigma_{\alpha} + M^{(1)}_{\alpha \gamma_1 \gamma_2} \sigma_{\gamma_1} \sigma_{\gamma_2} .\]
Notice that this procedure is just like modifying $\sigma_{\alpha}$ so that we may cancel $F^{(1)}$. The only distinction here is that $F^{(1)}$ is already established as vanishing, so we need only ensure that this uniquely fixes $M^{(1)}|_{\Lambda} = 0$. Indeed, so long as $k \neq d$, this is the case, by considerations in Section~\ref{conf-first-order}.

Similarly, the same style of argument applies to check that $\sigma_{\alpha}$ and $\sigma_{\alpha}'$ agree to third order. Again, we find that so long as $k \neq d-2$, $\sigma_{\alpha}$ is independent of the choice of $F^{(1)}$ modulo terms of order $\sigma^4$ -- the only pieces of $F^{(2)}$ that may be affected by a change in $\sigma$ at third order are precisely those terms that we forced to vanish in $F^{(2)}$. Hence, we must have that $M^{(2)}|_{\Lambda} = 0$.

While we do not proceed with the irreducible representation at third order (as the calculations get very messy but remain uninteresting), it is still clear that such a decomposition \textit{is} possible. Thus, one may choose $A^{(3)}$ to eliminate all but a preferred set of components of $F^{(3)}$. As~described above, one preferred component (when $k \neq d-2$) might be $F^{(3)}_{\gamma \rho \rho \gamma \alpha}$. Regardless of the choices made, so long as $\sigma_{\alpha}$ and $\sigma_{\alpha}'$ are constructed so that $F^{(3)}$ and $F^{(3)'}$ have the same vanishing components, $M^{(3)}$ is necessarily fixed to zero. This implies that the non-vanishing components of $F^{(3)}|_{\Lambda}$ are independent of the choice of $F^{(1)}$ and $\sigma_{\alpha}$ is fixed modulo terms of order $\sigma^5$.

We now turn our attention to the case where $k=d-2$. Rather than dealing with a unique $\sigma_{\alpha}$ for a given choice of $F^{(1)}$, we instead deal with a unique equivalence class $[\sigma_{\alpha}]$. We will attempt to fix the behavior of our equivalence class at fourth order, according to
\[[\sigma_{\alpha}] \mapsto \bigl[\sigma_{\alpha} + A^{(3)}_{\alpha \gamma_1 \gamma_2 \gamma_3 \gamma_4} \sigma_{\gamma_1} \sigma_{\gamma_2} \sigma_{\gamma_3} \sigma_{\gamma_4}\bigr] .\]
To enable us to work concretely with this equivalence class, we will pick an arbitrary representative $\sigma_{\alpha}$ and then parametrize the equivalence class by the arbitrary set of functions $A$. That is, we consider
\[\tilde{\sigma}_{\alpha} := \sigma_{\alpha} + A_{[\alpha \gamma]} \sigma_{\gamma} \sigma_{\rho} \sigma_{\rho} + A^{(3)}_{\alpha \gamma_1 \gamma_2 \gamma_3 \gamma_4} \sigma_{\gamma_1} \sigma_{\gamma_2} \sigma_{\gamma_3} \sigma_{\gamma_4} ,\]
where $A$ is treated as a free (function-valued) set of parameters. Then we have that
\begin{gather*}
\tilde{\bG}_{\alpha \beta} = \bG_{\alpha \beta} +\biggl( 8 A^{(3)}_{(\alpha \beta) \gamma_1 \gamma_2 \gamma_3} - \frac{24}{d} \delta_{\gamma_1 (\alpha} A^{(3)}_{\beta) \rho \rho \gamma_2 \gamma_3} + 2 \delta_{\gamma_2 \gamma_3} N_{(\alpha} \csdot \hd A_{\beta) \gamma_1}\\
\hphantom{\tilde{\bG}_{\alpha \beta} = \bG_{\alpha \beta} +\biggl(}{}
+ \frac{8}{d} N_{\gamma_3} \csdot \hd A_{\gamma_2 (\alpha} \delta_{\beta) \gamma_1} + \frac{4}{d} N_{\rho} \csdot \hd A_{\rho (\alpha} \delta_{\beta) \gamma_1} \delta_{\gamma_2 \gamma_3} \biggr) \sigma_{\gamma_1} \sigma_{\gamma_2} \sigma_{\gamma_3}.
\end{gather*}
We wish to examine the maximal trace terms again. These take the form
\begin{align*}
&F^{(3)}_{\gamma \gamma \rho \rho \alpha }- \frac{8}{d} A^{(3)}_{\alpha \gamma \gamma \rho \rho} + \frac{8(d-2)}{d} A^{(3)}_{\gamma \gamma \rho \rho \alpha} - \frac{2(d+2)}{3} N_{\gamma} \csdot \hd A_{\alpha \gamma} ,\\
& F^{(3)}_{\gamma \rho \rho \gamma \alpha }, \\
&F^{(3)}_{\alpha \gamma \gamma \rho \rho}- \frac{4}{d} A^{(3)}_{\alpha \gamma \gamma \rho \rho} + \frac{4(d-2)}{d} A^{(3)}_{\gamma \gamma \rho \rho \alpha} - \frac{d+2}{3} N_{\gamma} \csdot \hd A_{\alpha \gamma} .
\end{align*}
Evidently, there is no choice of $A^{(3)}$ that can eliminate \smash{$F^{(3)}_{\gamma \rho \rho \gamma \alpha}$}, regardless of the representative~that was chosen. Furthermore, so long as \smash{$F^{(3)}_{\gamma \gamma \rho \rho \alpha} \neq 2 F^{(3)}_{\alpha \gamma \gamma \rho \rho}$}, this system is inconsistent, and hence we may only choose to set one of these two maximal trace terms to zero. For consistency with the case where $k \neq d-2$, we can always fix $A^{(3)}$ so that \smash{$F^{(3)}_{\gamma \gamma \rho \rho \alpha} \mapsto \tilde{F}^{(3)}_{\gamma \gamma \rho \rho \alpha} = 0$} -- this choice is made so that we may compute a Willmore-like obstruction explicitly, as will be shown in the next section. In that case, we have that
\[\tilde{F}^{(3)}_{\alpha \gamma \gamma \rho \rho} = F^{(3)}_{\alpha \gamma \gamma \rho \rho} - \frac{1}{2} F^{(3)}_{\gamma \gamma \rho \rho \alpha} .\]
Importantly, $\tilde{F}^{(3)}_{\alpha \gamma \gamma \rho \rho}$ is \textit{independent of the choice of representative}. Furthermore (as in the argument for the case where $k \neq d-2$), for different choices of extensions of $F^{(1)}$, if both $[\sigma_{\alpha}]$ and $[\sigma_{\alpha}']$ are constructed so that \smash{$F^{(3)}_{\gamma \gamma \rho \rho \alpha}$} and \smash{$F^{(3)'}_{\gamma \gamma \rho \rho \alpha}$} vanish, then $M^{(3)}$ cannot affect \smash{$F^{(3)}_{\alpha \gamma \gamma \rho \rho}|_{\Lambda}$}. It thus follows that this obstruction is independent of the choice of $F^{(1)}$, and thereby is a true invariant of the parallelized submanifold embedding.

Observe that this part of $F^{(3)}$ is only an obstruction when $k=d-2$. This observation is key and suggests that this component of $F^{(3)}$ is a special invariant of the system. Furthermore, the uniqueness of this obstruction is fortuitous, as it is precisely this term that mimicks the obstruction density at order $d$ that occurs when this same procedure is applied to hypersurfaces~\cite{Will1}. In fact, for $d=3$ and $k=1$, this is precisely the celebrated Willmore invariant, as we will see in the following section. This suggests that for all $\Lambda^2 \hookrightarrow\bigl(M^d,\cc\bigr)$, we have that $F^{(3)}_{\alpha \gamma \gamma \rho \rho}$ represents a parallelized submanifold Willmore invariant.

\section{A holographic submanifold Willmore invariant} \label{Willmore-section}
As noted above, when $d-k=2$, $F^{(3)}$ has a special codimension-dependent obstruction that can be forced into one of its maximal traces, namely \smash{$F^{(3)}_{\alpha \gamma \gamma \rho \rho }$}. In this section, we provide a holographic formula for this obstruction and compute it explicitly.

Drawing inspiration from~\cite{Will1}, we can construct an extrinsically-coupled conformal Laplacian on $\Lambda$. Indeed, this conformal Laplacian exists on all tractor bundles with weight $\frac{2+k-d}{2}$. This is captured by the following straightforward result.
\begin{Theorem} \label{extr-lap}
Let $\Lambda^{d-k} \hookrightarrow \bigl(M^d,\cc\bigr)$ be a conformal submanifold embedding with some conformal defining map that satisfies $\bG_{\alpha \beta} = \delta_{\alpha \beta} + \mathcal{O}\bigl(\sigma^2\bigr)$. Then the operator
\[P_{2} := N_{\alpha} \csdot \hd N_{\alpha} \csdot D \colon \ \Gamma \biggl(\ct^\Phi M \biggl[\frac{2+k-d}{2} \biggr] \biggr)\Big|_{\Lambda} \rightarrow \Gamma \biggl(\ct^\Phi M \biggl[\frac{-2+k-d}{2} \biggr] \biggr)\Big|_{\Lambda} \]
is tangential and has principal symbol $-k \Delta^\top$, where $\Delta^\top := \bg^{ab} \nabla^\top_a \nabla^\top_b$ and $\nabla^\top$ is the tangential tractor-coupled Levi-Civita connection.
\end{Theorem}
\begin{proof}
Explicit computation in a choice of metric representative $g \in \cc$ yields
\begin{align*}
N_{\alpha} \csdot \hd N_{\alpha} \csdot D &\= [\nabla_{n_{\alpha}} + (w-1)\rho_{\alpha} ] \circ [(d+2w-2) (\nabla_{n_{\alpha}} + w\rho_{\alpha}) - \sigma_{\alpha} (\Delta + w J) ] \\
&\= -k \biggl(\Delta^\top + \frac{d-k-2}{2} (\nabla_{n_\alpha} \rho_{\alpha}) - \frac{d-k-2}{2} J - \frac{(d-k)(d-k-2)}{4} H_{\alpha} H_{\alpha} \biggr) .
\end{align*}
Evidently, this operator is tangential.
\end{proof}

\begin{Remark}
Note that this operator $P_2$ is a submanifold generalization of the extrinsic conformal Laplacian described in~\cite{Will1}.
\end{Remark}

We are now equipped to provide a holographic formula for \smash{$F^{(3)}_{\alpha \gamma \gamma \rho \rho}$}. From a straightforward computation, one finds that
\begin{align*}
\begin{split}
N_{\gamma} \csdot \hd N_{\gamma} \csdot \hd N_{\alpha}^A &{}\= -\frac{3}{d-2} X^A F^{(3)}_{\alpha \gamma_1 \gamma_1 \gamma_2 \gamma_2} + \frac{1}{d-2} X^A N_{\gamma} \csdot \hd K_{\gamma \alpha} - N_{\gamma} \csdot \hd B^A_{\alpha \gamma} \\
&\quad{}+ \frac{1}{d-2} N^A_{\gamma} K_{\gamma \alpha} + N^A_{\gamma_1} F^{(2)}_{\alpha \gamma_2 \gamma_2 \gamma_1} - \frac{1}{d-2} N^A_{\gamma_1} F^{(2)}_{\alpha \gamma_1 \gamma_2 \gamma_2} .
\end{split}
\end{align*}
Because we are interested in $F^{(3)}$, it is sensible to only consider the projection of the above display to the tractor bundle on $\Lambda$. As such, we have
\begin{align*} 
\top P_{2} N_{\alpha}^A \={} -3 X^A F^{(3)}_{\alpha \gamma_1 \gamma_1 \gamma_2 \gamma_2} + X^A N_{\gamma} \csdot \hd K_{\gamma \alpha} - (d-2) \top \bigl(N_{\gamma} \csdot \hd B^A_{\alpha \gamma}\bigr) ,
\end{align*}
where $\top$ also refers to the projection to the submanifold tractor bundle, such that \smash{$N^A_{\alpha} T^\top_A \= 0$} for any tractor $T$, i.e., $\top T^A := T^A - N^A_{\alpha} N_{\alpha} \csdot T$. From the principal symbol calculation of $P_2$ and the fact that \smash{$F^{(3)}_{\alpha \gamma_1 \gamma_1 \gamma_2 \gamma_2}$} appears as a coefficient for $X_A$, we have that the leading structure of this maximal-trace obstruction is $\bar{\Delta} H_{\alpha}$ -- as expected for a Willmore invariant. So, we may indeed use this obstruction to define a holographic Willmore invariant for a parallelized surface embedded in $d$ dimensions.

To explicitly compute this Willmore invariant requires more care. Indeed, if one were to calculate naively, one would realize that this computation involves tractor expressions such as
\[N_{\gamma} \csdot \hd F^{(2)}_{\omega [\gamma \alpha] \omega} .\]
For a given choice of $F^{(1)}$, these terms are uniquely determined, as $F^{(2)}$ is determined to all orders by $F^{(1)}$ and $A^{(1)}$ \big(which itself is fixed in terms of $F^{(1)}$\big). Furthermore, as $P_2$ is a tangential operator and \smash{$F^{(3)}_{\alpha \gamma \gamma \rho \rho}$} is a parallelized submanifold invariant, the dependence of such terms on the choice of $F^{(1)}$ must either cancel with one another or themselves vanish. Thus, such terms may simply be absorbed into $F^{(3)}$ by writing schematically
\[F^{(2)} = F^{(2c)} + F^{(2,1)} \sigma ,\]
where
\[F^{(2c)} := F^{(2)} - \sigma_{\alpha} N^A_{\alpha} \hd_A F^{(2)} .\]
Then $N_{\alpha} \csdot \hd F^{(2c)} \= 0$ and the terms involving $F^{(2,1)}$ are simply absorbed into $F^{(3)}$. While this may seem like we are changing the obstruction density, this is equivalent to simply choosing another $F^{(1)}$, as doing so simply reshuffles pieces of $F^{(2)}$ off $\Lambda$ into $F^{(3)}$. However, as \smash{$F^{(3)}_{\alpha \gamma \gamma \rho \rho}|_{\Lambda}$} is independent of $F^{(1)}$, this choice may be made with no ill-effects. (Note that this is the beginning of a solution to an extension problem which will be discussed in greater detail in Section~\ref{ext-prob-sec}.) Thus, for the sake of the computation that follows, we may drop terms involving $N_{\alpha} \csdot \hd F^{(2)}$.

To explicitly compute \smash{$F^{(3)}_{\alpha \gamma \gamma \rho \rho}$}, we must first handle $\top N_{\gamma} \csdot \hd B^A_{\alpha \gamma}$.
\begin{gather*}
2\top (N_{\gamma} \csdot \hd B^A_{\alpha \gamma}) \= \top \bigl[N_{\gamma} \csdot \hd \bigl(N_{\alpha} \csdot \hd N^A_{\gamma} - N_{\gamma} \csdot \hd N^A_{\alpha}\bigr) \bigr] \\
 \= \top \bigl(- \bigl(N_{\gamma} \csdot \hd\bigr)^2 N^A_{\alpha} + \bigl[N_{\gamma} \csdot \hd, N_{\alpha} \csdot \hd\bigr] N^A_{\gamma} + N_{\alpha} \csdot \hd N_{\gamma} \csdot \hd N^A_{\gamma} \bigr) \\
 \= \top \bigl(- \bigl(N_{\gamma} \csdot \hd\bigr)^2 N^A_{\alpha} + W_{\gamma \alpha}{}^A{}_{\gamma} + 2 B_{B \gamma \alpha} P^{BA}_{\gamma} + N_{\alpha} \csdot \hd N_{B \gamma} P^{AB}_{\gamma} \bigr) \\
 \= \top \biggl(- \bigl(N_{\gamma} \csdot \hd\bigr)^2 N^A_{\alpha} + W_{\gamma \alpha}{}^A{}_{\gamma} + 2 B_{B \gamma \alpha} P^{BA}_{\gamma} + N_{\alpha} \csdot \hd \biggl(B^A_{\gamma \gamma} + \frac{1}{2} \hd^A G_{\gamma \gamma} + \frac{1}{d-2} X^A K_{\gamma \gamma}\biggr) \biggr) \\
 \= \top \biggl(- \bigl(N_{\gamma} \csdot \hd\bigr)^2 N^A_{\alpha} + W_{\gamma \alpha}{}^A{}_{\gamma} + 2 B_{B \gamma \alpha} P^{BA}_{\gamma} + N_{\alpha} \csdot \hd \biggl( \frac{1}{2} \hd^A G_{\gamma \gamma} + \frac{1}{d-2} X^A K_{\gamma \gamma}\biggr) \biggr) \\
 \= \top \biggl(- \bigl(N_{\gamma} \csdot \hd\bigr)^2 N^A_{\alpha} + W_{\gamma \alpha}{}^A{}_{\gamma} + 2 B_{B \gamma \alpha} P^{BA}_{\gamma} + \frac{1}{2} N_{\alpha} \csdot \hd \hd^A G_{\gamma \gamma} + \frac{1}{d-2} X^A N_{\alpha} \csdot \hd K_{\gamma \gamma} \biggr) .
\end{gather*}
In the above computations, we have used the fact that $\hd_A X_B = h_{AB}$. Also, $W_{\gamma \alpha}{}^A{}_{\gamma}$ denotes the $W$-tractor contracted suitably with the tractor $N$ with the corresponding Greek indices. Now, to complete this calculation, we examine the term $\top N_{\alpha} \csdot \hd \hd^A G_{\gamma \gamma}$. To do so, we recall that we may drop terms of the form $N_{\alpha} \csdot \hd F^{(2)}$ along $\Lambda$. Thus, we have
\begin{gather*}
\top N_{\alpha} \csdot \hd \hd^A G_{\gamma \gamma} \= \top N_{\alpha} \csdot \hd \hd^A \bigl(F^{(2)}_{\gamma \gamma \beta_1 \beta_2} \sigma_{\beta_1} \sigma_{\beta_2} + F^{(3)}_{\gamma \gamma \beta_1 \beta_2 \beta_3} \sigma_{\beta_1} \sigma_{\beta_2} \sigma_{\beta_3}\bigr) \\
\qquad{}\= \top N_{\alpha} \csdot \hd \biggl( 2 \sigma_{(\beta_1} N^A_{\beta_2)} F^{(2)}_{\gamma \gamma \beta_1 \beta_2} - \frac{2}{d-2} X^A \bigl[ G_{\beta_1 \beta_2} F^{(2)}_{\gamma \gamma \beta_1 \beta_2} + 2 \sigma_{(\beta_1} N_{\beta_2)} \csdot \hd F^{(2)}_{\gamma \gamma \beta_1 \beta_2} \bigr] \biggr) \\
\quad\qquad{} + \top N_{\alpha} \csdot \hd \biggl(-\frac{6}{d-2} X^A \sigma_{(\beta_1} \bG_{\beta_2 \beta_3)} F^{(3)}_{\gamma \gamma \beta_1 \beta_2 \beta_3} + \mathcal{O}\bigl(\sigma^2\bigr) \biggr) \\
\qquad{}\= -\frac{6}{d-2} X^A \bG_{\alpha (\beta_1} \bG_{\beta_2 \beta_3)} F^{(3)}_{\gamma \gamma \beta_1 \beta_2 \beta_3} \\
\qquad{}\= -\frac{6}{d-2} X^A F^{(3)}_{\gamma \gamma \alpha \rho \rho} \\
\qquad{}\= 0 .
\end{gather*}
In this computation, we used the fact that $F^{(3)}_{\gamma \gamma \alpha \rho \rho}|_{\Lambda}$ is the only maximal trace that we have chosen to remove in general.

Combining the above two results with the formula for $\top \bigl(N_{\gamma} \csdot \hd\bigr)^2 N^A_{\alpha}$,
we have thus obtained a formula for the obstruction in terms of a tangential operator and lower-order terms,
\begin{gather}
X^A F^{(3)}_{\alpha \gamma \gamma \rho \rho} \= \top \biggl(-\frac{1}{6} P_{2} N^A_{\alpha} - \frac{d-2}{6} W_{\gamma \alpha}{}^A{}_{\gamma} - \frac{d-2}{3} B_{B \gamma \alpha} P^{BA}_{\gamma} \nonumber\\ \hphantom{X^A F^{(3)}_{\alpha \gamma \gamma \rho \rho} \= \top \biggl(}{}
 + \frac{1}{3} X^A N_{\gamma} \csdot \hd K_{\gamma \alpha} - \frac{1}{6} X^A N_{\alpha} \csdot \hd K_{\gamma \gamma} \biggr) .\label{F3-expression}
\end{gather}
To express the obstruction in terms of Riemannian quantities, we require a few more calculations. Note that both $P_2 N^A_{\alpha} \= \nabla^{\top b} \nabla^\top_b N^A_{\alpha}$ and $B_{B \gamma \alpha} P^{BA}_{\gamma} \= \beta^b_{\gamma \alpha} \nabla^\top_b N^A_{\gamma}$ both contain $\nabla^\top N_{\alpha}$, so we compute
\begin{align} \label{dtN}
\nabla^\top_a N_{B \alpha} \stackrel{\Lambda,g}{=} \begin{pmatrix}
0\\
\nabla^\top_a n_{b \alpha} - \bar{g}_{ab} H_{\alpha} \\
- \bar{\nabla}_a H_{\alpha} - \bar{g}_{ab} P^{bc} n^c_{\alpha}
\end{pmatrix}
 .
\end{align}
Further computation of $P_2 N$ is tedious but straightforward: a simple application of the projected tractor connection suffices.

Furthermore, computation of the last two terms can also be performed explicitly by noting that
\[N_{\alpha} \csdot \hd K_{\beta \gamma} \= (\nabla_{n_{\alpha}} + 2 H_{\alpha})\bigl[(\nabla_a n_{b \beta})\bigl(\nabla^a n^b_{\gamma}\bigr) - d \rho_{\beta} \rho_{\gamma} + 2 s_{(\beta} P^{ab} \nabla_a n_{b \gamma)}\bigr] .\]
We may thus calculate
\begin{gather}
N_{\alpha} \csdot \hd K_{\beta \gamma} \= 2\bigl(\nabla^a n^b_{(\beta}\bigr) \nabla_{|\alpha|} \nabla_a n_{b \gamma)} + 2d H_{(\beta} \nabla_{|\alpha|} \rho_{\gamma)} + 2 \delta_{\alpha (\beta} P^{ab} \nabla_a n_{b \gamma)} + 2 H_{\alpha} K_{\beta \gamma} \nonumber\\
\qquad{}\= 2\bigl(\II^{ab}_{(\beta} + n^a_{\rho} \beta^b_{\rho (\beta} + n^b_{\rho} \beta^a_{\rho (\beta} + n^a_{\rho} n^b_{\rho} H_{(\beta}\bigr) \nabla_{|\alpha|} \nabla_a n_{b \gamma)} + 2d H_{(\beta} \nabla_{|\alpha|} \rho_{\gamma)}
\nonumber\\
\qquad\quad{} + 2 \delta_{\alpha (\beta} P^{ab} \bigl(\II_{ab \gamma)} + n_{a |\rho} \beta_{b \rho| \gamma)} + n_{b |\rho} \beta_{a \rho| \gamma)} + n_{a |\rho} n_{b \rho|} H_{\gamma)}\bigr) + 2 H_{\alpha} K_{\beta \gamma} \nonumber\\
\qquad{}\= 2\bigl(\II^{ab}_{(\beta} + n^a_{\rho} \beta^b_{\rho (\beta} + n^b_{\rho} \beta^a_{\rho (\beta} +\bigl(g^{ab} - \bar{g}^{ab}\bigr) H_{(\beta}\bigr) \nabla_{|\alpha|} \nabla_a n_{b \gamma)} - 2 H_{(\beta} \nabla_{|\alpha|} \bigl(\Delta s_{\gamma)} + J s_{\gamma)}\bigr)
\nonumber\\
\qquad\quad{} + 2 \delta_{\alpha (\beta} P^{ab} \bigl(\II_{ab \gamma)} + n_{a |\rho} \beta_{b \rho| \gamma)} + n_{b |\rho} \beta_{a \rho| \gamma)} + n_{a |\rho} n_{b \rho|} H_{\gamma)}\bigr) + 2 H_{\alpha} K_{\beta \gamma} \nonumber\\
\qquad{}\= 2\bigl(\IIo^{ab}_{(\beta} + n^a_{\rho} \beta^b_{\rho (\beta} + n^b_{\rho} \beta^a_{\rho (\beta} \bigr) n^c_{|\alpha|} \nabla_c \nabla_a n_{b \gamma)}
\nonumber\\
\qquad\quad{} + 2 \delta_{\alpha (\beta} P^{ab} \bigl(\IIo_{ab \gamma)} + n_{a |\rho} \beta_{b \rho| \gamma)} + n_{b |\rho} \beta_{a \rho| \gamma)}\bigr) + 2 H_{\alpha} K_{\beta \gamma} \nonumber\\
\qquad{}\= 2n^c_{\alpha} \IIo^{ab}_{(\beta} \nabla_c \nabla_a n_{b \gamma)} + 2 n^c_{\alpha} n^a_{\rho} \beta^b_{\rho (\beta} \nabla_c \nabla_b n_{a \gamma)} + 2 n^c_{\alpha} n^b_{\rho} \beta^a_{\rho (\alpha} \nabla_c \nabla_a n_{b \gamma)}
\nonumber\\
\qquad\quad{} + 2 \delta_{\alpha (\beta} P^{ab} \bigl(\IIo_{ab \gamma)} + n_{a |\rho} \beta_{b \rho| \gamma)} + n_{b |\rho} \beta_{a \rho| \gamma)}\bigr) + 2 H_{\alpha} K_{\beta \gamma} \nonumber\\
\qquad{}\= 2n^c_{\alpha} \IIo^{ab}_{(\beta} \bigl(\nabla_a^\top \nabla_c n_{b \gamma)} + R_{cabd} n^d_{\gamma)}\bigr) + 2 n^c_{\alpha} n^a_{\rho} \beta^b_{\rho (\beta} \bigl(\nabla_b^\top \nabla_c n_{a \gamma)} + R_{cbad} n^d_{\gamma)} \bigr)\nonumber
\\
\qquad\quad{} + 2 n^c_{\alpha} n^b_{\rho} \beta^a_{\rho (\alpha} \bigl(\nabla^\top_a \nabla_c n_{b \gamma)} + R_{cabd} n^d_{\gamma)}\bigr) + 2 \delta_{\alpha (\beta} P^{ab} \bigl(\IIo_{ab \gamma)} + n_{a |\rho} \beta_{b \rho| \gamma)} + n_{b |\rho} \beta_{a \rho| \gamma)}\bigr) \nonumber\\
\qquad\quad{} + 2 H_{\alpha} K_{\beta \gamma} . \label{NDK}
\end{gather}

Further evaluating both equations~\eqref{dtN} and~\eqref{NDK} requires a few non-standard submanifold identities. Tracing the Gauss equation~(\ref{gauss-equation}) by applying $\bar{g} = g - n_{\alpha} \otimes n_{\alpha}$ to both sides yields the Fialkow--Gauss equation, namely
\begin{align*}
&(d-k-2) \bigl(P^\top_{ab} - \bar{P}_{ab} + H_{\alpha} \IIo_{ab \alpha} + \frac{1}{2} \bar{g}_{ab} H_{\alpha}^2 \bigr) \\
&\qquad{} = \IIo^2_{(ab)\alpha \alpha} - W_{\alpha ab \alpha}^\top - \frac{1}{2(d-k-1)} \bar{g}_{ab} \bigl(\IIo^2_{\alpha \alpha} + W_{\alpha \beta \beta \alpha}\bigr) .
\end{align*}
Similarly (when $d-k > 2$), tracing again with $\bar{g}$ yields a generalization of the \textit{theorema egregium},
\[J - P_{\alpha \alpha} = \bar{J} - \frac{d-k}{2} H_{\alpha} H_{\alpha} + \frac{1}{2(d-k-1)} \bigl(\IIo^2_{\alpha \alpha} - \bar{g}^{cd} W_{\alpha cd \alpha}\bigr) .\]
Similarly, tracing the Codazzi--Mainardi equation~(\ref{codazzi}) leaves us with
\[(d-k-1) P_{a \beta}^\top = \bar{\nabla} \csdot \IIo_{a \beta} - (d-k-1) \bigl(\bar{\nabla}_a H_{\beta} + \beta_{a \beta \alpha} H_{\alpha}\bigr) + \IIo^b_{a \alpha} \beta_{b \beta \alpha} + W_{\gamma \beta \gamma}{}^b \bar{g}_{ab} .\]
Conversely, we can take the tracefree part to obtain
\begin{gather*}
W_{abc \alpha}^\top - \frac{2}{d-k-1} \bigl(W_{\beta \alpha \beta [a}^\top \bar{g}_{b]c}\bigr) = 2 \bar{\nabla}_{[a} \IIo_{b]c \alpha} + 2\beta_{[a \alpha \beta} \IIo_{b]c \beta} \\ \hphantom{W_{abc \alpha}^\top - \frac{2}{d-k-1} \bigl(W_{\beta \alpha \beta [a}^\top \bar{g}_{b]c}\bigr) =}{}
- \frac{2}{d-k-1} \bar{g}_{c [a}\bigl( \bar{\nabla} \csdot \IIo_{b] \alpha} + \IIo_{b]c \beta} \beta^{c}_{\alpha \beta} \bigr).
\end{gather*}

We are now equipped to evaluate equations~\eqref{dtN} and~\eqref{NDK}. When $k=d-2$, we find that
\begin{align} \label{dtNeval}
\nabla^\top_a N_{B \alpha}
\stackrel{\Lambda,g}{=} \begin{pmatrix}
0\\
\IIo_{ab \alpha} + n_{b \beta} \beta_{a \beta \alpha} \\
-\bar{\nabla} \csdot \IIo_{a \alpha} - \beta^b_{\alpha \beta} \IIo_{ba \beta} + H_{\beta} \beta_{a \alpha \beta} - W_{\alpha \beta a \beta}^\top
\end{pmatrix}
\end{align}
and
\begin{align*}
\frac{1}{3} N_{\gamma} \csdot \hd K_{\gamma \alpha} - \frac{1}{6} N_{\alpha} \csdot \hd K_{\gamma \gamma}\stackrel{\Lambda,g}{=}& -\frac{2}{3} \bigl(\IIo^{ab}_{\beta} \bar{\nabla}_a \beta_{b \alpha \beta} + 2 \beta_{b \alpha \beta} \bar{\nabla} \csdot \IIo^b_{\beta}\bigr) - \frac{1}{3} \IIo^3_{\alpha \beta \beta} + \frac{1}{3} \beta^a_{\beta \gamma} \beta^b_{\beta \gamma} \IIo_{ab \alpha} \\
&{}- \frac{4}{3} \beta^a_{\alpha \beta} \beta^b_{\beta \gamma} \IIo_{ab \gamma} + \frac{2}{3} \beta^a_{\alpha \beta} W_{\beta \delta \delta a} + \frac{2}{3} \beta^a_{\beta \gamma} W_{\alpha \beta \gamma a} + \frac{1}{3} \IIo^{ab}_{\alpha} W_{\beta ab \beta} .
\end{align*}

Finally, we compute $\Delta^\top N_{B \alpha}$,
\begin{align}
\nonumber\Delta^\top N_{B \alpha} \stackrel{\Lambda,g}{=}& \begin{pmatrix}
0\\
\ast \\
\scalebox{0.75}{$
\begin{aligned}
-(\bar{\nabla}_a \bar{\nabla}_b + P^\top_{ab} + H_{\gamma} \IIo_{ab \gamma} ) \IIo^{ab}_{\alpha} - (\IIo^{ab}_{\gamma} \bar{\nabla}_a \beta_{b \alpha \gamma} + 2 \beta_{b \alpha \gamma} \bar{\nabla}_a \IIo^{ab}_{\gamma}) \\+ \beta_{a \alpha \beta} \beta_{b \beta \gamma} \IIo^{ab}_{\gamma} + 2 \beta_{a \alpha \beta} \bar{\nabla} \csdot \IIo^a_{\beta} + \beta^a_{\alpha \beta} W_{\beta \gamma \gamma a} - \bar{\nabla}^a W_{a \beta \beta \alpha}^\top
\end{aligned}
$}
\end{pmatrix}\\
&{}
- [\beta^a_{\alpha \beta} \beta_{a \gamma \beta} + \IIo^2_{\alpha \gamma} - \bar{\nabla} \csdot \beta_{\alpha \gamma}] N_{B \gamma}\label{LapTN}
\end{align}
Note that in the above computation, one must be careful to notice that the tractor in equation~(\ref{dtNeval}) is a \textit{bulk} tractor. After differentiation to obtain equation~(\ref{LapTN}), one must then project the resulting tractor to $\ct \Lambda$. For example, differentiating the middle slot of equation~(\ref{dtNeval}) yields a term in the middle slot of the form \smash{$\nabla^{\top a}\IIo_{ab \alpha} = \bar{\nabla} \csdot \IIo_{\alpha} - n_{b \gamma} \IIo_{ac \gamma} \IIo^{ac}_{\alpha}$}. Projecting the resulting tractor yields a term proportional to \smash{$H_{\gamma} \IIo^2_{\gamma \alpha}$} in the bottom slot, found in equation~(\ref{LapTN}).

With these identities in tow, we can explicitly compute, in terms of more familiar tensors, the obstruction \smash{$F^{(3)}_{\alpha \gamma \gamma \rho \rho}$} present in equation~(\ref{F3-expression}). Indeed, in a choice of scale $g \in \cc$, we find that
\begin{align*}
F^{(3)}_{\alpha \gamma \gamma \rho \rho} \={}& -\frac{d-2}{6} \bigl(\bar{\nabla}_a \bar{\nabla}_b + P^\top_{ab} + H_{\gamma} \IIo_{ab \gamma} \bigr) \IIo^{ab}_{\alpha} - \frac{1}{3} \IIo^{bc}_{\alpha} \bar{g}^{ad} W_{abcd} - \frac{1}{3} \IIo^{ab}_{\gamma} \IIo_{b \gamma}^c \IIo_{ac \alpha}
\\&- \frac{d-2}{6} \IVnn_{\alpha}- \frac{d+2}{6} \bigl[\IIo^{ab}_{\gamma} \bar{\nabla}_a \beta_{b \alpha \gamma} + 2 \beta_{b \alpha \gamma} \bar{\nabla}_a \IIo^{ab}_{\gamma} \bigr] + \frac{1}{3} \IIo^{ab}_{\alpha} \beta_{a \beta \gamma} \beta_{b\beta \gamma}
\\&- \frac{d+6}{6} \IIo^{ab}_{\gamma} \beta_{a \alpha \beta} \beta_{b \beta \gamma}+ \frac{2}{3} W_{\alpha \beta \gamma a} \beta^a_{\beta \gamma} + \frac{d+2}{6} W_{\beta \gamma \gamma a} \beta^a_{\alpha \beta} ,
\end{align*}
where
\[\IVnn_{\alpha} := C_{\alpha \beta \beta} + H_{\rho} W_{\alpha \beta \beta \rho} + \frac{1}{d-k-3} \bar{\nabla}^c W_{c \beta \beta \alpha}^\top \in \Gamma(\ce \Lambda[-3]\times \raisebox{1pt}{\scalebox{0.2}{\ydiagram{1}}})\]
is a novel conformal submanifold invariant for every $d-k \neq 3$. This can be verified by explicitly computing the conformal transformation $g \mapsto \Omega^2 g$. Furthermore, the combination
\[\IIo^{ab}_{\gamma} \bar{\nabla}_a \beta_{b \alpha \gamma} + \frac{2}{d-k-1} \beta_{b \alpha \gamma} \bar{\nabla}_a \IIo^{ab}_{\gamma} \in \Gamma(\ce \Lambda[-3] \times \raisebox{1pt}{\scalebox{0.2}{\ydiagram{1}}}) \]
is conformally invariant for all $d-k \neq 1$, and can be constructed from an appropriate combination of the tractor second fundamental form~\cite{Will1}, the tractor projection of $B$ to the submanifold, and the submanifold Thomas-$D$ operator. Finally, the operator
\[
L_{ab} := \bar{\nabla}_a \bar{\nabla}_b + P^\top_{ab} + H_{\gamma} \IIo_{ab \gamma} \colon \ \Gamma\bigl(\odot_\circ^2 T \Lambda[-3]\bigr) \rightarrow \Gamma(\ce \Lambda[-3])
\]
is invariant under local rescalings of the metric when $\Lambda$ is a surface, as is the case here.

Further simplification of the obstruction can be obtained by noting that in the setting of $d-k=2$, $\IIo_{\alpha}$ is a trace-free symmetric $2\times 2$
matrix. It follows from the Cayley--Hamilton theorem that $\IIo^a_{b \alpha} \IIo^b_{c \alpha}$ is proportional to the identity matrix, and so the $\IIo^3$ term vanishes, leaving us with
\begin{align} \label{surface-willmore}
F^{(3)}_{\alpha \gamma \gamma \rho \rho} \= -\frac{d-2}{6} L_{ab} \IIo^{ab}_{\alpha} - \frac{1}{3} \IIo^{bc}_{\alpha} \bar{g}^{ad} W_{abcd} - \frac{d-2}{6} \IVnn_{\alpha} + \text{terms containing } \beta .
\end{align}
Notably, this obstruction is similar to but distinct from the general codimension Willmore invariant for surfaces from~\cite{CGKTW}, expressed in our notation as
\begin{align*}
-L_{ab} \IIo^{ab}_{\alpha} + \IVnn_{\alpha} .
\end{align*}

Observe that when $d=3$ and $k=1$, all of $W$, $\IVnn_{\alpha}$, and $\beta$ vanish.
So for a surface embedded in a 3-volume, we have that
\[F^{(3)} \= -\frac{1}{6} \bigl(\bar{\nabla}_a \bar{\nabla}_b + P^\top_{ab} + H \IIo_{ab} \bigr) \IIo^{ab} ,\]
which matches~\cite{hypersurface_old} the famous Willmore invariant~\cite{willmore1965} up to a factor of $1/2$.

\begin{Remark}
An interesting question is whether this Willmore invariant is variational for higher codimension surface embeddings, and if so, what energy functional it comes from. We leave this for future work.
\end{Remark}

\section{The extension problem} \label{ext-prob-sec}
Recall from the previous section that it is sometimes useful to extend functions (or arrays thereof) on $\Lambda$ to a tubular neighborhood, at least formally. To that end, we consider these so-called \textit{extension problems} in both the Riemannian and conformal settings for submanifold embeddings. Such problems are of more general interest, see~\cite{GW}. A Riemannian problem can be stated as
\begin{Problem} \label{ext-problem}
Let $\bar{f} \in C^\infty \Lambda$ and let $s_{\alpha}$ be a canonical defining map for $\Lambda^{d-k} \hookrightarrow \bigl(M^d,g\bigr)$. Find a formal power series to as high order as possible for $f \in C^\infty M$ solving
\[\nabla_{n_{\alpha}} f = 0 , \qquad f|_{\Lambda} = \bar{f} .\]
\end{Problem}
However, unlike the hypersurface case, we will find that such a problem can not, in general, be solved. But, as formulated this problem depends on a choice of frame (and its extension off the submanifold) which makes this problem less interesting. Trivially, when $\bar{f}$ is constant, we may extend uniquely so that $f$ is constant as well. When $f$ is non-constant, we can modify this problem to become more geometrically meaningful and \textit{also} allow solutions. We will handle these two modifications first, and then we will consider the conformal extension problem.

\subsection{The symmetric extension problem}
Rather than solving Problem~\ref{ext-problem}, we instead weaken the constraint equations, so that the problem of interest takes the following form.

\begin{Problem} \label{sym-ext-problem}
Let $\bar{f} \in C^\infty \Lambda$ and let $s_{\alpha}$ be a canonical defining map for $\Lambda^{d-k} \hookrightarrow \bigl(M^d,g\bigr)$. For a given $m \in \mathbb{Z}_{\geq 1}$, find $f \in C^\infty M$ satisfying $f|_{\Lambda} = \bar{f}$ and
\[\nabla_{n_{\alpha}} f \= 0 , \qquad n^{a_1 a_2}_{(\alpha_1 \alpha_2)} \nabla_{a_1 a_2} f \= 0 , \qquad \ldots , \qquad n^{a_1 \cdots a_{m+1}}_{(\alpha_1 \cdots \alpha_{m})} \nabla_{a_1 \cdots a_{m}} f \= 0 .\]
\end{Problem}
\noindent
In the above and going forward, $n^{a_1 \cdots a_n}_{\alpha_1 \cdots \alpha_n} := n^{a_1}_{\alpha_1} \cdots n^{a_n}_{\alpha_n}$ and similarly $\nabla_{a_1 \cdots a_n} := \nabla_{a_1} \cdots \nabla_{a_n}$. To solve this problem requires a series of technical results.

\begin{Lemma} \label{normal-derivs-s}
Let $\Lambda \hookrightarrow (M,g)$ be a Riemannian submanifold embedding with a defining map $s_{\alpha}$ and let $m \in \mathbb{Z}_{\geq 2}$. Suppose that the defining map satisfies $G_{\alpha \beta} = \delta_{\alpha \beta} + \mathcal{O}\bigl(s^2\bigr)$ and
\[n^{a_1 \cdots a_{n-1}}_{(\alpha_1 \cdots \alpha_{n-1}} \nabla_{a_1 \cdots a_{n-1} b} G_{\gamma \delta)} \= 0 \= n^{a_1\cdots a_n}_{(\alpha_1 \cdots \alpha_n} \nabla_{a_1 \cdots a_n} G_{\gamma) \delta} \]
for all $2 \leq n \leq m$.
Then, for each $1 \leq p \leq m+1$, we have that
\[n^{a_1 \cdots a_p}_{(\alpha_1 \cdots \alpha_p)} \nabla_{a_1 \cdots a_p} s_{\beta} \= \begin{cases}
0, & p \neq 1, \\
\delta_{\alpha_1 \beta}, & p = 1 .
\end{cases}\]
\end{Lemma}
\begin{proof}
The case where $p=1$ follows trivially. We consider the case where $p \geq 2$ and prove via induction. When $p=2$, the identity reduces to the fact that \smash{$n \csdot \beta \= 0$}.

Now fix $p$, and suppose that the identity holds for all $k \leq p \leq m$. We wish to show that the identity holds for $k=p+1$. We then compute
\begin{align*}
 n^{a_{p}}_{\alpha_{p}} n^{a_{p+1}}_{\alpha_{p+1}} \nabla^{p-1} \nabla_{a_p} \nabla_{a_{p+1}} s_{\beta} \={}& n^{a_{p}}_{\alpha_{p}} \nabla^{p-1} \nabla_{a_{p}} G_{\alpha_{p+1} \beta} - n^{a_{p}}_{\alpha_{p}} n^{a_{p+1}}_{\beta} \nabla^{p-1} \nabla_{a_{p}} n_{a_{p+1} \alpha_{p+1}}\\
\={}& n^{a_{p}}_{\alpha_{p}} \nabla^{p-1} \nabla_{a_{p}} G_{\alpha_{p+1} \beta} -n^{a_{p}}_{\alpha_{p}} n^{a_{p+1}}_{\beta} \nabla^{p-1} \nabla_{a_{p+1}} n_{a_{p} \alpha_{p+1}}\\
\={}& n^{a_{p}}_{\alpha_{p}} \nabla^{p-1} \nabla_{a_{p}} G_{\alpha_{p+1} \beta} - \frac{1}{2} n^{a_{p+1}}_{\beta} \nabla^{p-1} \nabla_{a_{p+1}} G_{\alpha_{p} \alpha_{p+1}}\\
\={}& 0 .
\end{align*}
The first identity holds via the Leibniz identity and the induction hypothesis: moving $n$ past any number of covariant derivatives costs additional products, each of which has strictly fewer than $p+1$ symmetric normal derivatives on $s$. By the induction hypothesis, these additional terms vanish. The third equality follows by the same observation, as well as by noting that all such expressions are symmetric in $\alpha_{p}$ and $\alpha_{p+1}$. The final identity follows by the proposition hypothesis. The lemma follows.
\end{proof}

\begin{Corollary} \label{normal-derivs-n}
Let $\Lambda \hookrightarrow (M,g)$ be a Riemannian submanifold embedding with a defining map~$s_{\alpha}$ and let $m \in \mathbb{Z}_{\geq 2}$. Suppose that the defining map satisfies $G_{\alpha \beta} = \delta_{\alpha \beta} + \mathcal{O}\bigl(s^2\bigr)$ and
\[n^{a_1 \cdots a_{n-1}}_{(\alpha_1 \cdots \alpha_{n-1}} \nabla_{a_1 \cdots a_{n-1} b} G_{\gamma \delta)} \= 0 \= n^{a_1 \cdots a_n}_{(\alpha_1 \cdots \alpha_n} \nabla_{a_1 \cdots a_n} G_{\gamma) \delta} \]
for all $2 \leq n \leq m$.
Then, for every $1 \leq p \leq m$, we have that
\[n^{a_1 \cdots a_p}_{(\alpha_1 \cdots \alpha_p} \nabla_{a_1 \cdots a_p} n^b_{\beta)} \= 0 .\]
\end{Corollary}
\begin{proof}
We prove by induction. For $p = 1$, this follows simply because \smash{$n^a_{(\alpha} \nabla_a n^b_{\beta)} \= \beta^b_{(\alpha \beta)} \= 0$}. Now suppose the identity holds for all $k \leq p \leq m-1$. Now consider the case where $k = p+1$ and compute:
\begin{align*}
n^{a_1 \cdots a_p}_{(\alpha_1 \cdots \alpha_p} n^{a_{p+1}}_{\alpha_{p+1}} \nabla^{p} \nabla_{a_{p+1}} n^{b}_{\beta)} \={}& n^{a_1 \cdots a_p}_{(\alpha_1 \cdots \alpha_p} n^{a_{p+1}}_{\alpha_{p+1}} \nabla^{p} \nabla^{b} n_{a_{p+1} \beta)}
\= \frac{1}{2} n^{a_1 \cdots a_p}_{(\alpha_1 \cdots \alpha_p} \nabla^{p} \nabla^{b} G_{\alpha_{p+1} \beta)}
\= 0 .
\end{align*}
The second identity follows by commuting normal vectors to the right, which costs nothing via the Leibniz rule according to Lemma~\ref{normal-derivs-s}. The third identity follows from the hypothesis.
\end{proof}

\begin{Corollary} \label{normal-derivs-sn}
Let $\Lambda \hookrightarrow (M,g)$ be a Riemannian submanifold embedding with a defining map~$s_{\alpha}$ and let $m \in \mathbb{Z}_{\geq 2}$. Suppose that the defining map satisfies $G_{\alpha \beta} = \delta_{\alpha \beta} + \mathcal{O}\bigl(s^2\bigr)$ and
\[n^{a_1 \cdots a_{n-1}}_{(\alpha_1 \cdots \alpha_{n-1}} \nabla_{a_1 \cdots a_{n-1} b} G_{\gamma \delta)} \= 0 \= n^{a_1 \cdots a_n}_{(\alpha_1 \cdots \alpha_n} \nabla_{a_1 \cdots a_n} G_{\gamma) \delta} \]
for all $2 \leq n \leq m$.
Then, for each $1 \leq p \leq m+1$ and $1 \leq \ell$, we have that
\[n^{a_1 \cdots a_p}_{(\alpha_1 \cdots \alpha_p)} \nabla_{a_1 \cdots a_p} \bigl(s_{\beta_1} \cdots s_{\beta_{\ell}} n^{b_1}_{\beta_1} \cdots n^{b_{\ell}}_{\beta_{\ell}}\bigr) \= \begin{cases}
0, & p \neq \ell, \\
p! n^{b_1}_{(\alpha_1} \cdots n^{b_p}_{\alpha_p)}, & p = \ell .
\end{cases}\]
\end{Corollary}

\begin{proof}
For $p < \ell$, the result is trivial. Furthermore, for $p = \ell$, the result follows trivially by noting that each covariant derivative must act on a unique $s_{\beta}$, yielding $n_{\alpha} \csdot n_{\beta} = \delta_{\alpha \beta}$.

Finally, consider the case where $p > \ell$. Then, following Lemma~\ref{normal-derivs-s}, we may only have one derivative acting on each $s$. Thus, in each term in the Leibniz expansion there exists at least one normal vector that is acted upon by more than one normal derivative. However, note that each normal vector has had its Greek index $\beta_i$ has been replaced by a Greek index $\alpha_j$ for some~$j$ as a consequence of the action of the normal derivatives on $s$. Hence, we are examining terms of the form present in Corollary~\ref{normal-derivs-n}. Therefore, all such terms vanish and the result follows.
\end{proof}

We require one more technical lemma.

\begin{Lemma} \label{symmetric-deriv-G}
Let $\Lambda^{d-k} \hookrightarrow \bigl(M^d,g\bigr)$ be parallelized by $\{n_{\alpha}\}$ and let $s_{\alpha}$ be the corresponding canonical defining map.
Then
\[n^{a_1 \cdots a_{n-1}}_{(\alpha_1 \cdots \alpha_{n-1}} \nabla_{a_1 \cdots a_{n-1} b} G_{\gamma \delta)} \= 0 \= n^{a_1 \cdots a_n}_{(\alpha_1 \cdots \alpha_n} \nabla_{a_1 \cdots a_n} G_{\gamma) \delta} \]
for each $n \in \mathbb{Z}_{\geq 1}$.
\end{Lemma}
\begin{proof}
We will prove this by analyzing a fixed term $F^{(i)}s^i$ in the expansion of $G_{\alpha \beta}$ and using induction on $n$. The base case ($n=1$) trivially holds because $G = \delta + \mathcal{O}\bigl(s^2\bigr)$. Now fix some $n \geq 2$, and suppose that for each $1 \leq j \leq n$, we have that
\[n^{a_1 \cdots a_{j-1}}_{(\alpha_1 \cdots \alpha_{j-1}} \nabla_{a_1 \cdots a_{j-1} b} G_{\gamma \delta)} \= 0 \= n^{a_1 a_2 \cdots a_j}_{(\alpha_1 \alpha_2 \cdots \alpha_j} \nabla_{a_1 \cdots a_j} G_{\gamma) \delta} .\]
We would like to show that this identity holds in the case where $j=n+1$ as well. To do so, it is sufficient to show that, for each $i$, we have that
\begin{align} \label{decomp-window}
n^{a_1 \cdots a_{n+1}}_{(\alpha_1 \cdots \alpha_{n+1}} \nabla_{a_1 \cdots a_{n+1}} \bigl(F^{(i)}_{\gamma) \delta \gamma_1 \cdots \gamma_i} s_{\gamma_1} \cdots s_{\gamma_i}\bigr) \= 0 \= n^{a_1 \cdots a_{n}}_{(\alpha_1 \cdots \alpha_{n}} \nabla_{a_1 \cdots a_{n} b} \bigl(F^{(i)}_{\gamma \delta) \gamma_1 \cdots \gamma_i} s_{\gamma_1} \cdots s_{\gamma_i}\bigr) .\!
\end{align}
Certainly, for $n+1 < i$, we find that the identity holds as \smash{$\nabla^{n+1} s^i \= 0$}. If $n+1 = i$, then we must compute the appropriate symmetrization over $i+1$ indices of $F^{(i)}$. However, for either symmetrization, because $F^{(i)}$ occupies a generalized window representation according to Theorem~\ref{general-ho-ext}, this term necessarily vanishes.

Finally, consider the case where $n+1 > i$. The leading derivative term on $F^{(i)}$ takes the form $\nabla^{n-i+1} F^{(i)}$. However, such a term vanishes because $F^{(i)}$ (even away from $\Lambda$) occupies the generalized window representation, and this term symmetrizes over $i+1$ indices of $F^{(i)}$. On the other hand, subleading terms must be considered: those terms of the form $\bigl(\nabla^{\ell} s^i\bigr) \nabla^{n+1-\ell} F^{(i)}$, where $\ell > i$. Note that at least one $s$ is differentiated twice, but no $s$ is differentiated more than $n-i+2$ times. Note that via the induction hypothesis, we may use both Lemma~\ref{normal-derivs-s} and Corollary~\ref{normal-derivs-n} to examine these terms, as $n-i+2 \leq n$ (because $i \geq 2$).

Now, consider the term on the left-hand side of equation~(\ref{decomp-window}). All such terms of the form $\bigl(\nabla^\ell s^i\bigr) \nabla^{n+1-\ell} F^{(i)}$ for this term contain at least one multiplicand that vanishes due to Lemma~\ref{normal-derivs-s}. Thus, we need only consider the right-hand side of equation~(\ref{decomp-window}).

In that case, there are two types of terms in a Leibniz expansion that must be examined. The first type of term is one that contains a multiplicand of the form to which Lemma~\ref{normal-derivs-s} is applicable. These terms must therefore vanish. Hence, we must only examine terms of the second type, namely those terms where each of $i-1$ factors of $s$ are acted upon by $\nabla_n$ exactly once, and one remaining $s$ is acted upon by $\nabla_n^{\ell-i} \nabla_b$:
\[\bigl(n^{a_1 \cdots a_{\ell-i}}_{(\alpha_1 \cdots \alpha_{\ell-i}} \nabla_{a_1 \cdots a_{\ell-i}b }s_{|\gamma_1|}\bigr) \nabla^{n+1-\ell} F^{(i)}_{\gamma \delta) \gamma_1 \alpha_1 \cdots \alpha_{\ell-1}} .\]
But note that here $F^{(i)}$ is symmetrized over $i+1$ indices, and therefore must vanish.

Having shown that all terms in a Leibniz expansion for both the left and right-hand sides of equation~(\ref{decomp-window}) vanish for all $i$, we have that the lemma follows.
\end{proof}

We are now equipped to solve Problem~\ref{sym-ext-problem}.
\begin{Theorem} \label{sym-ext-prop}
Let $\Lambda \hookrightarrow (M,g)$ be a parallelized Riemannian submanifold embedding with a~canonical defining map $s_{\alpha}$, let $\bar{f} \in C^\infty \Lambda$, and let $m \in \mathbb{Z}_{\geq 2}$. Then there exists a unique $\tilde{f} \in C^\infty M$ modulo terms of order $m+2$ with $\tilde{f}|_{\Lambda} = \bar{f}$ satisfying
\[\nabla_{n_{\alpha}} \tilde{f} \= 0 , \qquad n^{a_1 a_2}_{(\alpha_1 \alpha_2)} \nabla_{a_1 a_2} \tilde{f} \= 0 , \qquad \ldots , \qquad n^{a_1 \cdots a_{m+1}}_{(\alpha_1 \cdots \alpha_{m+1})} \nabla_{a_1 \cdots a_{m+1}} \tilde{f} \= 0 .\]
\end{Theorem}

\begin{proof}
We begin by fixing $m$. Then let $f \in C^\infty M$ be any function satisfying $f|_{\Lambda} = \bar{f}$ and define
\begin{align}
\tilde{f} :={}& f - s_{\beta_1} n^{b_1}_{\beta_1} \nabla_{b_1} f + \frac{1}{2} s_{\beta_1} s_{\beta_2} n^{b_1}_{\beta_1} n^{b_2}_{\beta_2} \nabla_{b_1} \nabla_{b_2} f\nonumber \\
&- \frac{1}{6} s_{\beta_1} s_{\beta_2} s_{\beta_3} n^{b_1}_{\beta_1} n^{b_2}_{\beta_2} n^{b_3}_{\beta_3} \nabla_{b_1} \nabla_{b_2} \nabla_{b_3} f + \cdots .\label{sym-ext}
\end{align}
Now, we would like to show that this ansatz satisfies the conditions of the theorem. In particular, we must check that
\[n^{a_1 \cdots a_q}_{(\alpha_1 \cdots \alpha_q)} \nabla_{a_1 \cdots a_q} \tilde{f} \]
vanishes for all $1 \leq q \leq m+1$. Now we examine a generic term from $\tilde{f}$, say the term containing~$s^{\ell}$. That is, consider
\[\frac{(-1)^{\ell}}{\ell!} n^{a_1 \cdots a_q}_{(\alpha_1 \cdots \alpha_q)} \nabla_{a_1 \cdots a_q} \bigl(s_{\beta_1} \cdots s_{\beta_{\ell}} n^{b_1}_{\beta_1} \cdots n^{b_\ell}_{\beta_{\ell}} \nabla_{b_1 \cdots b_{\ell}} f\bigr) .\]
Trivially this expression vanishes when $q < \ell$. When $q \geq \ell$, we may rewrite this expression as
\[\frac{(-1)^{\ell}}{\ell!} n^{a_1 \cdots a_q}_{(\alpha_1 \cdots \alpha_q)} \sum_{r=0}^{q} {q \choose r} \bigl(\nabla_{a_1 \cdots a_r} s_{\beta_1} \cdots s_{\beta_{\ell}} n^{b_1}_{\beta_1} \cdots n^{b_\ell}_{\beta_{\ell}}\bigr)(\nabla_{a_{r+1} \cdots a_{q}} \nabla_{b_1 \cdots b_{\ell}} f) .\]
Now from Lemma~\ref{symmetric-deriv-G}, we have that the hypothesis of Corollary~\ref{normal-derivs-sn} holds for every $m$. Furthermore, as $q \geq \ell$, we have that this sum includes the case where $r = \ell$. Thus, we have that
\begin{gather*}
\frac{(-1)^{\ell}}{\ell!} n^{a_1 \cdots a_q}_{(\alpha_1 \cdots \alpha_q)} \sum_{r=0}^{q} {q \choose r} \bigl(\nabla_{a_1 \cdots a_r} s_{\beta_1} \cdots s_{\beta_{\ell}} n^{b_1}_{\beta_1} \cdots n^{b_\ell}_{\beta_{\ell}}\bigr)(\nabla_{a_{r+1} \cdots a_{q}} \nabla_{b_1 \cdots b_{\ell}} f) \\
\qquad{}= (-1)^{\ell} {q \choose \ell} n^{a_1 \cdots a_q}_{(\alpha_1 \cdots \alpha_q)} \nabla_{a_1 \cdots a_q} f .
\end{gather*}
It follows then that
\[n^{a_1 \cdots a_q}_{(\alpha_1 \cdots \alpha_q)} \nabla_{a_1 \cdots a_q} \tilde{f} \= \Biggl(\sum_{\ell = 0}^q (-1)^{\ell} {q \choose \ell} \Biggr) n^{a_1 \cdots a_q}_{(\alpha_1 \cdots \alpha_q)} \nabla_{a_1 \cdots a_q} f \= 0 ,\]
by the binomial theorem. Thus $\tilde{f}$ satisfies the required derivative constraints.

To complete the proof, it suffices to check that if we had initially chosen a distinct ${f' \in C^\infty M}$ and built the ansatz accordingly, we find that $\tilde{f}'-\tilde{f} = s_{\alpha_1} \cdots s_{\alpha_{m+2}} O_{\alpha_1 \cdots \alpha_{m+2}}$, for some smooth $O_{\alpha_1 \cdots \alpha_{m+2}}$. Now, note that because \smash{$f'|_{\Lambda} \= \bar{f}$}, we must have that $f' - f = s_{\alpha} T_{\alpha}$ for some smooth~$T_{\alpha}$. We must check that the first $m+1$ symmetric normal derivatives of this difference vanishes. However, from the above computation it is clear that $m+1$ symmetric normal deri\-vatives of \smash{$\widetilde{(s_{\alpha} T_{\alpha})}$} vanishes. And as the ansatz operator is linear, this shows that the difference of the ansatz must vanish to sufficiently high order. This completes the proof.
\end{proof}

\begin{Remark}
Note that the condition in the display of Theorem~\ref{sym-ext-prop} is independent of the choice of parallelization, and hence by uniqueness of the result, so is the extension.
\end{Remark}

\subsection{The restricted extension problem}

Alternatively, we may hope that we can extend $\bar{f}$ to force many normal derivatives to vanish as in Problem~\ref{ext-problem} if certain constraints are placed on the embedding map $\Lambda^{d-k} \hookrightarrow \bigl(M^d,g\bigr)$.
We find that, contrary to the case of hypersurface embeddings, which have solutions to all orders by the Picard--Lindel{\"o}f theorem, this extension problem encounters an obstruction. To see this, we proceed order by order.

Suppose that $f_0 \in C^\infty M$ is any function such that $f_0|_{\Lambda} = \bar{f}$, and let $f_1 = f_0 - s_{\alpha} \nabla_{n_{\alpha}} f_0$. Then a straightforward computation shows that
\[\nabla_{n_{\alpha}} f_1 = F_{\alpha \beta} s_{\beta} \]
for some array of functions $F$. Now suppose that we wish to proceed to a higher order by writing
\[f_2 = f_1 + s_{\alpha} s_{\beta} A_{\alpha \beta} .\]
Observe that here, we may consider only $A_{\alpha \beta} = A_{(\alpha \beta)}$. If $f_2$ is to solve Problem~\ref{ext-problem} to second order, then we must have that $A \propto F$. However, $F$ is not necessarily symmetric. While we can choose $A$ so that $F_{(\alpha \beta)}$ is cancelled, we cannot do the same for $F_{[\alpha \beta]}$. In fact, this is our first obstruction, and it can be calculated:
\begin{align*}
F_{[\alpha \beta]} \= -\nabla_{n_{[\alpha}} \nabla_{n_{\beta]}} f_1
\= -n^a_{[\alpha} \bigl(\nabla_a n^b_{\beta]}\bigr) \nabla_b f_1 - n^a_{\alpha} n^b_{\beta} \nabla_{[a} \nabla_{b]} f_1
\= -\beta^b_{\alpha \beta} \bar{\nabla}_b \bar{f} .
\end{align*}
Indeed, unless $\bar{f}$ is constant or $\beta = 0$, the submanifold extension problem cannot be solved beyond first order. The case where $\bar{f}$ is constant is uninteresting, as the solution is merely that $f = \bar{f}$. So going forward, we assume that $\beta = 0$ to proceed -- that is, that $\Lambda \hookrightarrow (M,g)$ has vanishing normal curvature and we choose a parallelization of $N \Lambda$ such that $\beta = 0$. Then, from the above considerations, we have that there exists $f_2$ such that
\[\nabla_{n_{\alpha}} f_2 = F_{\alpha \beta \gamma} s_{\beta} s_{\gamma} .\]
We wish to adjust $s$ such that $F$ vanishes along $\Lambda$, as before. So, we write $f_3 = f_2 + A_{\alpha \beta \gamma} s_{\alpha} s_{\beta} s_{\gamma}$. Indeed, by purely representation theory considerations, we have that $F_{\alpha \beta \gamma} \in C^\infty (M ,\raisebox{3pt}{\scalebox{0.2}{\ydiagram{2,1}}} \; \oplus \; \raisebox{2pt}{\scalebox{0.2}{\ydiagram{3}}})$, whereas $A_{\alpha \beta \gamma} \in C^\infty (M, \raisebox{1.5pt}{\scalebox{0.2}{\ydiagram{3}}})\,$. We now attempt to compute this obstruction as we did for the obstruction at first order.

To do so, a straightforward calculation shows that
\begin{align*}
F_{\alpha \beta \gamma} &\= \frac{1}{2} \nabla_{n_{\beta}} \nabla_{n_{\gamma}} \nabla_{n_{\alpha}} f_2
\= \frac{1}{2} n^b_{\beta} n^c_{\gamma} \bigl(\nabla_b \nabla_c n^a_{\alpha}\bigr) \nabla_a f_2 + \frac{1}{2} n^b_{\beta} n^c_{\gamma} n^a_{\alpha} \nabla_b \nabla_c \nabla_a f_2 \\
&\= \frac{1}{2} n^b_{\beta} n^c_{\gamma} \bigl(\nabla_b \nabla_c n^a_{\alpha}\bigr) \bar{\nabla}_a \bar{f} + \frac{1}{2} n^b_{\beta} n^c_{\gamma} n^a_{\alpha} \nabla_b \nabla_c \nabla_a f_2 \\
&\= \frac{1}{2} n^b_{\beta} n^c_{\gamma} \bigl(\nabla^{\top\, a} \nabla_b n_{c \alpha}\bigr) \bar{\nabla}_a \bar{f} + \frac{1}{2} R_{\gamma \alpha \beta}{}^a \bar{\nabla}_a \bar{f} + \frac{1}{2} n^b_{\beta} n^c_{\gamma} n^a_{\alpha} \nabla_b \nabla_c \nabla_a f_2 \\
&\= \frac{1}{2} R_{\gamma \alpha \beta}{}^a \bar{\nabla}_a \bar{f} + \frac{1}{2} n^b_{\beta} n^c_{\gamma} n^a_{\alpha} \nabla_b \nabla_c \nabla_a f_2 ,
\end{align*}
where the third identity follows from the known behavior of the normal derivative of $f_2$, and the last identity follows because $\beta = 0$. Of interest is the component of $F$ that occupies the $\raisebox{3pt}{\scalebox{0.2}{\ydiagram{2,1}}}$ representation -- the totally symmetric component is removable via a judicious choice of $A$,
\begin{gather*}
P_{\scalebox{0.2}{\ydiagram{2,1}}}\, F_{\alpha \beta \gamma} \= \frac{1}{3}(F_{\alpha \beta \gamma} + F_{\beta \alpha \gamma} - F_{\gamma \beta \alpha} - F_{\gamma \alpha \beta})
\= \bigl(P_{\scalebox{0.2}{\ydiagram{2,1}}}\, R_{\gamma \alpha \beta}{}^a\bigr) \bar{\nabla}_a \bar{f} \\ \hphantom{P_{\scalebox{0.2}{\ydiagram{2,1}}}\, F_{\alpha \beta \gamma} \=}{}
+ \frac{1}{6} \bigl(
n^b_{\beta} n^c_{\gamma} n^a_{\alpha} \nabla_b \nabla_c \nabla_a f_2
 + n^b_{\beta} n^c_{\gamma} n^a_{\alpha} \nabla_a \nabla_c \nabla_b f_2
- n^b_{\beta} n^c_{\gamma} n^a_{\alpha} \nabla_b \nabla_a \nabla_c f_2 \\ \hphantom{P_{\scalebox{0.2}{\ydiagram{2,1}}}\, F_{\alpha \beta \gamma} \=+ \frac{1}{6} \bigl(}{}
 - n^b_{\beta} n^c_{\gamma} n^a_{\alpha} \nabla_a \nabla_b \nabla_c f_2
 \bigr) \\ \hphantom{P_{\scalebox{0.2}{\ydiagram{2,1}}}\, F_{\alpha \beta \gamma}}{}
 \= \bigl(P_{\scalebox{0.2}{\ydiagram{2,1}}}\, R_{\gamma \alpha \beta}{}^a\bigr) \bar{\nabla}_a \bar{f} .
\end{gather*}
Thus, we see that while not a necessary condition, a sufficient condition to proceed to higher orders is to demand that $\bigl(M^d,g\bigr)$ is flat. In fact, this assumption allows us to solve the problem entirely.
\begin{Theorem} \label{ho-extension}
Suppose $\Lambda^{d-k} \hookrightarrow \bigl(M^d,g\bigr)$ is parallelizable, let $s_{\alpha}$ be a canonical defining map, and let $\bar{f} \in C^\infty \Lambda$. Then there exists $f \in C^\infty \Lambda$ such that $f|_{\Lambda} = \bar{f}$ and
\[\nabla_{n_{\alpha}} f = \mathcal{O}(s) .\]
Solutions are obstructed at second order by the invariant
\[\beta^b_{\alpha \beta} \bar{\nabla}_b \bar{f} .\]
Further, if the embedding has vanishing normal curvature and the embedding is parallelized by the rotation minimizing frame, then we may find $f \in C^\infty M$ satisfying
\[\nabla_{n_{\alpha}} f = \mathcal{O}\bigl(s^2\bigr) .\]
Finally, if the embedding has vanishing normal curvature, $\bigl(M^d,g\bigr)$ is flat, and the embedding is parallelized by the rotation minimizing frame, then there exists a unique $f \in C^\infty M$ satisfying
\[\nabla_{n_{\alpha}} f = 0 .\]
\end{Theorem}
\begin{proof}
The first two cases follow by the above considerations. Now suppose that $\mathcal{R} = 0 = R$. We prove by induction. Suppose that there exists $f \in C^\infty M$ such that $f|_{\Lambda} = \bar{f}$ and
\[\nabla_{n_{\alpha}} f = F_{\alpha \beta_1 \cdots \beta_{m-1}} s_{\beta_1} \cdots s_{\beta_{m-1}} .\]
First, we will show that $F$ is totally symmetric. By differentiating $m$ times, we find that
\[\frac{1}{(m-1)!} F_{\alpha \beta_1 \cdots \beta_{m-1}} = \nabla_{n_{\beta_1}} \cdots \nabla_{n_{\beta_{m-1}}} \nabla_{n_{\alpha}} f .\]
From Lemma~\ref{lower-order-derivs}, it follows that up to $m-1$ normal derivatives on $n^a_{\alpha}$ vanish. So,
\[\frac{1}{(m-1)!} F_{\alpha \beta_1 \cdots \beta_{m-1}} = n^a_{\alpha} n^{b_1}_{\beta_1} \cdots n^{b_{m-1}}_{\beta_{m-1}} \nabla_{b_1} \cdots \nabla_{b_{m-1}} \nabla_a f .\]
By hypothesis, covariant derivatives commute, and so $F$ is totally symmetric.

Now define $\tilde{f} = f + A_{\gamma_1 \cdots \gamma_{m}} s_{\gamma_1} \cdots s_{\gamma_{m}}$ for arbitrary $A$. Then
\begin{align*}
 \nabla_{n_{\gamma_1}} \cdots \nabla_{n_{\gamma_{m}}} \tilde{f} &\= \nabla_{n_{\gamma_1}} \cdots \nabla_{n_{\gamma_{m}}} (f + A_{\alpha_1 \cdots \alpha_{m}} s_{\alpha_1} \cdots s_{\alpha_{m}}) \\
&\= (m-1)! F_{\gamma_{m} (\gamma_1 \cdots \gamma_{m-1})} + m! A_{\gamma_1 \cdots \gamma_{m}} .
\end{align*}
Because $F$ is totally symmetric, there exists a judicious choice of $A$ such that we can totally eliminate $F$, and so we have that $\nabla_{n_{\alpha}} \tilde{f} = \mathcal{O}(s^m)$. This completes the induction.

Now observe that there is a unique extension that sets the symmetric normal derivatives (to a certain order) to zero, as given~by equation~(\ref{sym-ext}). As the extension constructed here has all normal derivatives vanishing (up to a certain order) it thus follows that it is indeed unique.
\end{proof}

We have thus determined sufficient conditions to solve Problem~\ref{ext-problem}.

\subsection{The conformal extension problem} 
Considering the success of the Riemannian symmetric extension problem, we consider the conformal equivalent.
\begin{Problem} \label{conformal-ext-problem}
Let $\Lambda \hookrightarrow (M,\cc)$ be a parallelized conformal submanifold embedding with an associated defining density $\sigma_{\alpha}$, let $\bar{f} \in \Gamma(\ce \Lambda[w])$, and let $m \in \mathbb{Z}_{\geq 2}$.
Find $f \in \Gamma(\ce M[w])$ satisfying $f|_{\Lambda} = \bar{f}$ and
\[N_{\alpha} \csdot D f \= 0 , \qquad N^{A_1 A_2}_{(\alpha_1 \alpha_2)} D_{A_1 A_2} f \= 0 , \qquad \ldots , \qquad N^{A_1 \cdots A_{m+1}}_{(\alpha_1 \cdots \alpha_{m+1}} D_{A_1 \cdots A_{m+1})}f \= 0 .\]
\end{Problem}

We provide the following ansatz to second order:
\begin{align*}
\tilde{f} :={}& f - \sigma_{\beta_1} N^{B_1}_{\beta_1} \hd_{B_1} f + \frac{1}{2} \sigma_{\beta_1} \sigma_{\beta_2} N^{B_1}_{(\beta_1} N^{B_2}_{\beta_2)_\circ} \hd_{B_1} \hd_{B_2} f \\&+ \frac{d-2k+2w-2}{2k(d-k+2w-2)} \sigma_{\beta_1} \sigma_{\beta_1} N^{B_1}_{\beta_2} N^{B_2}_{\beta_2} \hd_{B_1} \hd_{B_2}f
\\
&+ \frac{w}{d-2} \sigma_{\beta_1} \sigma_{\beta_2} \bigl(K_{\beta_1 \beta_2} - F^{(2)}_{\beta_1 \beta_2 \gamma \gamma}\bigr) f+\mathcal{O}(\sigma^3)
 ,
\end{align*}
where \smash{\raisebox{0.5pt}{$N^{B_1}_{(\beta_1} N^{B_2}_{\beta_2)_\circ}$} satisfies \smash{\raisebox{0.5pt}{$\delta_{\beta_1 \beta_2} N^{B_1}_{(\beta_1} N^{B_2}_{\beta_2)_\circ} \= 0$}}}. Formally, we are allowing ourselves to interchange between $D$ and $\hd$, noting that this will exclude some weights at the end of the construction.

By direct computation, one can check that indeed,
\[N \csdot D_{\alpha} \tilde{f} \= 0 \= N^{A_1}_{(\alpha_1} N^{A_2}_{\alpha_2)} \hd_{A_1} \hd_{A_2} \tilde{f} .\]
However, as the Thomas-$D$ operator is not a derivation, it becomes increasingly challenging to construct $\tilde{f}$ so that more symmetric normal Thomas-$D$s vanish on $\tilde{f}$. For brevity, we will halt the solution here, observing that we can find a canonical extension unique to second order in $f$ that solves Problem~\ref{conformal-ext-problem} to second order. Observe that this extension is only well defined when $w \neq 1-\frac{d-k}{2},1-\frac{d}{2}$.

For those sporadic weights, we must be more careful. When $w = 1-\frac{d}{2}$, we have trivially that \smash{$N_{\alpha} \csdot D f \= 0$}, but the condition is vacuous and hence we lose uniqueness -- so we are halted here. However, when $w = 1-\frac{d-k}{2}$, we may instead omit the offending term from our ansatz. In doing so, the ansatz no longer solves the problem to second order, but instead solves the following system:
\[N \csdot D_{\alpha} \tilde{f} \= 0 , \qquad N^{A_1 A_2}_{(\alpha_1 \alpha_2)_\circ} D_{A_1 A_2} \tilde{f} \= 0 .\]
Notably, we are obstructed from finding a solution that solves the full problem and instead can only solve the trace-free equivalent. This observation is consistent with the existence of the extrinsically-coupled conformal submanifold Laplacian operator described Theorem~\ref{extr-lap}.
In fact, a slight modification of the operator $P_{2}$ yields the problematic operator in the ansatz. To see this, we first perform a useful computation,
\begin{align*}
N_{\gamma} \csdot \hd N^A_{\beta} &{}= N^B_{\gamma} \hd^A N_{B \beta} \\
&{}= N_{B(\gamma} \hd^A N^B_{\beta)} + N_{B[\gamma} \hd^A N^B_{\beta]} \\
&{}= \frac{1}{2} \hd^A G_{\gamma \beta} + \frac{1}{d-2} X^A K_{(\gamma \beta)} + B^A_{\gamma \beta} \\
&{}= \frac{1}{2} \hd^A \bigl(\sigma_{\rho_1} \sigma_{\rho_2} F^{(2)}_{\gamma \beta \rho_1 \rho_2}\bigr) + \frac{1}{d-2} X^A K_{(\gamma \beta)} + B^A_{\gamma \beta} \\
&{}= \sigma_{(\rho_1} N^A_{\rho_2)} F^{(2)}_{\gamma \beta \rho_1 \rho_2} - \frac{2}{d-2} X^A \sigma_{(\rho_1} N_{\rho_2)} \csdot \hd F^{(2)}_{\gamma \beta \rho_1 \rho_2} - \frac{1}{d-2} X^A \bG_{\rho_1 \rho_2} F^{(2)}_{\gamma \beta \rho_1 \rho_2} \\
&\phantom{=} + \frac{1}{d-2} X^A K_{(\gamma \beta)} + B^A_{\gamma \beta} + \mathcal{O}\bigl(\sigma^2\bigr) .
\end{align*}
Consequently, we may write the $P_{2}$ in terms of the operator $P_2^\top := N^{B_1}_{\beta_1} N^{B_2}_{\beta_2} \hd_{B_1} \hd_{B_2}$ via
\begin{align*}
P_{2} \=k N_{\alpha} \csdot \hd N_{\alpha} \csdot \hd
\= k P_2^\top + k \bigl(N_{\alpha} \csdot \hd N^B_{\alpha}\bigr) \hd_B
\= k P_2^\top + \frac{k(k-d+2)}{2(d-2)} \bigl(K_{\alpha \alpha} - F^{(2)}_{\alpha \alpha \beta \beta}\bigr) .
\end{align*}
Notably, as $P_{2}$ is a tangential operator, so is $P_2^\top$. So we find that, when $\bar{f} \in \Gamma\bigl(\ce \Lambda\bigl[\frac{2+k-d}{2}\bigr]\bigr)$ satisfies \smash{$P_2^\top \bar{f} \= 0$}, the ansatz without the offending term solves the whole extension problem~to~order two -- otherwise, this obstruction genuinely prevents a full solution.

We may now summarize the results of this subsection in the following proposition.
\begin{Theorem} \label{conformal-ext-solution}
Let $\Lambda^{d-k} \hookrightarrow \bigl(M^d,\cc\bigr)$ be a parallelized conformal submanifold embedding with an associated defining density $\sigma_{\alpha}$.
Then, for any $\bar{f} \in \Gamma(\ce \Lambda[w])$ with $w \neq 1 - (d-k)/2, 1-d/2$, there exists a unique $f \in \Gamma(\ce M[w])$ modulo terms of order $\sigma^3$ such that
\[N \csdot D_{\alpha} f \= 0 , \qquad N^{A_1 A_2}_{(\alpha_1 \alpha_2)} D_{A_1 A_2} f \= 0 .\]
Furthermore, if \smash{$P_2^\top \bar{f} \= 0$}, there exists a unique solution \big(modulo terms of order $\sigma^3$\big) even when $w = 1- (d-k)/2$.
\end{Theorem}

\subsection*{Acknowledgements}

SB and J\v{S} are supported by the Czech Science Foundation (GACR) grant GA22-00091S. SB is also supported by the Operational Programme Research Development and Education Project No. CZ.02.01.01/00/22-010/0007541. SB and J\v{S} would like to thank A.~Rod Gover and Andrew Waldron for helpful comments and insights and the anonymous referees for relevant contributions to improve this work.

\pdfbookmark[1]{References}{ref}
\LastPageEnding

\end{document}